\documentclass[a4wide,11pt,onetabnum,reqno]{article}
\usepackage{bm}

\usepackage{mathrsfs}
\usepackage[nottoc]{tocbibind} 

\usepackage{fullpage}
\usepackage{textcomp}
\usepackage{enumerate}
\usepackage[table]{xcolor}
\usepackage[utf8]{inputenc}
\usepackage{amssymb, amsfonts}
\usepackage{caption}
\usepackage{subcaption}
\usepackage{booktabs}
\usepackage{siunitx}

\usepackage{amsbsy,amsthm} 
\usepackage{tikz}
\usetikzlibrary{patterns,calc,decorations.markings}
\usepackage{graphics}
\usepackage{caption}
\usepackage{amsmath}
\usepackage{algpseudocode}
\usepackage{color}
\usepackage{subcaption}
\usepackage{float,epsfig}
\usepackage{fixmath}

\usepackage{authblk}

\usepackage{pgfpages}
\usepackage{appendix}
\usepackage{hyperref} 
\usepackage{graphicx}
\usepackage{pst-all}
\usepackage{calligra}
\usepackage{cleveref}
\usepackage{algorithm}
\usepackage{multirow}
\usepackage{array}
\newcolumntype{H}{>{\setbox0=\hbox\bgroup}c<{\egroup}@{}}
\usepackage{accents}
\DeclareMathAlphabet{\mathpzc}{OT1}{pzc}{m}{it}
\DeclareMathAlphabet{\mathcalligra}{T1}{calligra}{m}{n}
\DeclareMathAlphabet{\mathpzc}{OT1}{pzc}{m}{it}
\DeclareMathAlphabet{\mathcalligra}{T1}{calligra}{m}{n}
\DeclareMathAlphabet{\mathpzc}{OT1}{pzc}{m}{it}
\DeclareMathAlphabet{\mathcalligra}{T1}{calligra}{m}{n}
\usepackage{blindtext}
\newtheorem{theorem}{\bf Theorem}[section]

\newtheorem{exam}[theorem]{\bf Example}
\newtheorem{remark}[theorem]{\bf Remark}
\newtheorem{lemma}[theorem]{\bf Lemma}
\newtheorem{assum}[theorem]{\bf Assumption}

\newcommand{\be}{\begin{equation}}
\newcommand{\ee}{\end{equation}}
\newcommand{\beno}{\begin{equation*}}
\newcommand{\eeno}{\end{equation*}}
\newcommand{\ba}{\begin{align}}
\newcommand{\ea}{\end{align}}
\newcommand{\bano}{\begin{align*}}
\newcommand{\eano}{\end{align*}}
\newcommand{\bea}{\begin{eqnarray}}
\newcommand{\eea}{\end{eqnarray}}
\newcommand{\beano}{\begin{eqnarray*}}
\newcommand{\eeano}{\end{eqnarray*}}

\renewcommand{\arraystretch}{3}
\numberwithin{equation}{section}

\def\arraystretch{1.3}

\def \cK{{\mathcal K}}

\def \Tc{{\mathcal T}}

\def \Pc{{\mathcal P}}

\def \Ac{{\mathcal A}}

\def \Ts{\mathscr{T}}

\def \Ac {\mathpzc{A}}
\def \Tc {\mathpzc{T}}

\def \max{{\rm max}}
\def \min{{\rm min}}
\def \sup{{\rm sup}}
\def \esssup{{\rm ess~sup}}

\begin{document}

\title{Space-Time Finite Element Approximation of Quasilinear Hyperbolic Equations Arising in Dynamic Strain-Limiting Elasticity}
\author[${\ast,a}$]{Ram Manohar}
\author[ ${\dagger,a}$]{Ananya G. Hegade}
\author[${\ddagger,a}$]{Harry Lee}
\author[{\textsection, $a$} ]{S. M. Mallikarjunaiah}
\affil[$a$]{Department of Mathematics \& Statistics,  Texas A\&M University-Corpus Christi, 6300 Ocean Dr., Corpus Christi, TX-78412, USA} \smallskip  \vspace{0.1cm} 
\affil[ \textsection]{Corresponding author E-mail: m.muddamallappa@tamucc.edu   \newline
\textit{E-mail addresses:~} {\texttt{$^\ast$ram.manohar@tamucc.edu},~ $^\dagger$\texttt{ahegde@islander.tamucc}, ~ and~  $^\ddagger$\texttt{harry.lee@tamucc.edu}}
}
\date{} 
\maketitle
\begin{abstract}
The numerical approximation of a class of quasilinear hyperbolic equations arising in dynamic strain-limiting elasticity—whose nonlinear constitutive law relates stress and strain nonlinearly—is investigated. Ellipticity may be lost by such problems in regions of large strain, leading to significant challenges in the design of stable and accurate numerical methods. To address these difficulties, a fully discrete continuous Galerkin finite element framework is developed, combining continuous linear finite elements for spatial discretization with the Hilber--Hughes--Taylor (\texttt{HHT}-$\alpha$) time integration scheme. A lifting technique is incorporated within the proposed formulation to treat non-homogeneous Dirichlet boundary conditions, and a consistent Newton iteration is employed for the efficient solution of the resulting nonlinear systems. Robustness is enhanced by the algorithmic dissipation introduced through the \texttt{HHT}-$\alpha$ method by suppressing nonphysical high-frequency oscillations near degenerate regions, while consistency and accuracy are preserved. Under suitable structural assumptions on the nonlinear constitutive coefficient, the weak formulation is established in appropriate energy spaces. Theoretical second-order convergence in the $L^2$-norm and first-order convergence in the $H^1$-norm, rapid nonlinear convergence with nearly mesh-independent Newton iterations, physically consistent wave-speed evolution, and stable energy dissipation are demonstrated through comprehensive numerical experiments. Furthermore, it is confirmed by the results that an accurate, stable, and computationally efficient framework for simulating nonlinear strain-limiting wave propagation is provided, offering a solid foundation for future extensions to multidimensional nonlinear elastodynamics, adaptive finite element methods, and fracture and damage mechanics.
\end{abstract}

\vspace{.1in}

\textbf{Key words.}
Quasilinear hyperbolic problem, continuous Galerkin FEM, $\texttt{HHT}$-$\alpha$-time integration.
	
\vspace{.1in}
\textbf{AMS subject classifications.} $65\mathrm{N}12$; $65\mathrm{N}15$; $65\mathrm{N}22$; $65\mathrm{N}30$; $65\mathrm{N}50$.
\section{Introduction}
A rigorous analysis of wave propagation is a fundamental part of modern engineering research and is the basis for many key technological developments in many fields. In the field of geophysical engineering, these analyses are essential for subsurface imaging and hydrocarbon exploration. Accurate modeling of seismic wave dispersion enables the identification of fluid reservoirs and geological faults. Likewise, in the realm of civil infrastructure, wave-based techniques are the essence of structural health monitoring and nondestructive evaluation. Here, ultrasonic wave scattering is used to detect nascent defects, e.g., micro-cracks or corrosion, in bridges, dams, and pipelines before they propagate to catastrophic failure. Furthermore, these principles are also applied in biomedical engineering, especially in elastography where shear wave propagation is studied to determine the mechanical properties of soft tissues for tumor detection. In traditional approaches, these fields are identified in the Classical Linearized Theory of Elasticity, where the material is assumed to be homogeneous, and the stress-strain relationship is strictly proportional. However, empirical evidence from complex materials such as rocks, concrete, and biological tissues shows that experimental responses often differ from these linear predictions. These deviations are fundamental changes in the behavior, which are often classified as \textit{classical physical nonlinearity, geometric nonlinearity, and hysteretic nonlinearity}. In these cases, the use of linearized models is insufficient to describe the complexity of real-world materials and requires the implementation of advanced nonlinear mathematical frameworks.
\par 
The central physical challenge is the behavior of materials near stress
concentrators such as crack tips. Classical linear elasticity, which is the
standard mathematical model of elastic wave propagation, predicts that
the strain field becomes singular at a crack tip, growing without bound.
This is physically inadmissible: no real material can sustain an
infinitely large strain. Yet this mathematical failure occurs precisely
at the location that is most important to engineers and material
scientists who seek to predict crack initiation and growth. Addressing
this failure is one of the primary motivations for the present work.
\par
In recent years, Rajagopal \cite{rajagopal2003implicit, rajagopal2007elasticity,rajagopal2011conspectus} introduced a novel implicit framework for constitutive modeling of elastic bodies that overcomes this limitation. In this framework, the linearized strain $\boldsymbol{\varepsilon}$ is expressed as a nonlinear function of the Cauchy stress $\mathbf{T}$, rather than the other way around. A specific and physically important subclass of these models, known as the strain-limiting models, ensures that the strain field remains uniformly bounded throughout the body regardless of how large the stress becomes. This effectively regularizes the crack-tip singularities of classical theory and produces physically realistic predictions near fractures \cite{rajagopal2011modeling, gou2015modeling,mallikarjunaiah2015direct,kulvait2013anti}. The present work extends this framework to the \textit{dynamic} setting of wave propagation, which has not been comprehensively studied before.
\par
A pioneering contribution in this direction is the theory of elasticity of Rajagopal \cite{rajagopal2003implicit, rajagopal2007elasticity, rajagopal2011conspectus}, which overcomes the limitations of the classical Cauchy and Green formulations via a rigorous thermodynamic framework. Unlike traditional approaches relating stress to strain, this approach involves implicit relations where the linearized strain ($\varepsilon$) relates to the Cauchy stress ($\mathbb{T}$) by a nonlinear function. 
\begin{equation}
\varepsilon =  \mathbf{f}(\mathbb{T})   
\end{equation}
An important subgroup of these models is strain-limiting theories, which are characterized by the requirement that the norm of the response function stays confined by a positive constant regardless of how high the stress is. This unique property effectively regularizes the mathematical singularities commonly found at crack tips in linear models, making this framework ideal for simulating fracture evolution in brittle materials and providing new insights into dynamic phenomena such as enhanced tip velocities and crack branching. While Rajagopal's founding work and following contributions by Malek and S\"{u}li have considerably advanced our knowledge of static boundary value issues, the dynamic extensions of these theories are mostly unknown. The current work on wave propagation in this context is limited to a specific subset of constitutive relations, leaving a full study of finite element analysis for dynamic wave propagation inside a body in a condition of anti-plane shear, using a constitutive model where
\begin{equation}
\varepsilon = \phi(|\mathbb{T}|)\frac{\mathbb{T}}{2}\;\;\; \text{with} \;\;\; \phi(r)=(1+ (\beta r)^\gamma)^{-1/\gamma}
\end{equation}
The strain-limiting subclass of these models is defined by the
requirement that the response function satisfies the uniform bound. This single property resolves the crack-tip singularity problem in the strain-limiting framework; the strain field near a crack tip remains bounded and physically admissible, regardless of how large the stress concentration becomes \cite{rajagopal2011modeling, gou2015modeling}. By inverting this relationship and substituting it into the balance of linear momentum, we derive the governing \textbf{second-order quasilinear hyperbolic partial differential equation} \cite{bustamante2015direct, kannan2014unsteady}: 
\begin{equation} 
\rho \frac{\partial^2 u}{\partial t^2} -\left(\frac{\nabla u}{[1-\beta ^\gamma |\nabla u|^\gamma]^{1/ \gamma}}\right )= 0
 \nonumber
\end{equation}
This motivates the formulation of a general  \textbf{second-order quasilinear hyperbolic partial differential equation}.

Our objective is to determine the displacement $u(x,t)$, defined for $(x,t) \in [\texttt{a},\,\texttt{b}] \times [0,\,\texttt{T}]$ a finite time horizon $\texttt{T} \in \mathbb{R}$ 
\begin{subequations}
\begin{align}
\rho\, \frac{\partial^2 u}{\partial t^2}-\frac{\partial}{\partial x}\Big(\mathpzc{N}(x, \big|\frac{\partial u}{\partial x}\big|)\, \frac{\partial u}{\partial x}  \Big) &=  f(x,t),~~~~~~ \text{in}~~ \quad (x,t)\in (\texttt{a},\, \texttt{b})\times (0,\,\texttt{T}], \label{hypc1.1a}\\
 u(\texttt{a},t)=0,  \quad u(\texttt{b},t)&=\psi(t),~~~~~~~~~ \text{for}~~ \quad  t \in (0,\texttt{T}],\\
 u(x,0)=\phi_1(x),  \quad \frac{\partial u}{\partial t}(x,0)&=\phi_2(x),~~~~~~~~ \text{for}~~ \quad  x \in (\texttt{a},\, \texttt{b}), \label{hypc1.1c}
\end{align}  
\end{subequations}
where the nonlinear term is defined by 
\begin{align*}
\mathpzc{N}\big(x, \big|\frac{\partial u}{\partial x} \big| \big):=\frac{1}{\big[1-\beta^\gamma \, \big|\frac{\partial u}{\partial x} \big|^\gamma \big]^{\frac{1}{\gamma}}}.
\end{align*}
To avoid the cumbersome notation, we will use the following convention from now onwards: $\frac{\partial u}{\partial t}:= \partial_t u ~ \text{or}~ \dot{u}$,  and $\frac{\partial u}{\partial x}:=\partial_xu$, while the $k-th$ order partial derivative of $u$ with respect to $x$ is represented by $\partial^{(k)}_xu$, respectively. 
\par
The mathematical and computational literature on strain-limiting elasticity has developed along four main directions. Rajagopal and his collaborators established the constitutive framework based on implicit constitutive relations and demonstrated its thermodynamic consistency while showing that classical linear elasticity is recovered as a limiting case \cite{rajagopal2003implicit, rajagopal2007elasticity, rajagopal2011conspectus, rajagopal2007response, rajagopal2009class}. Subsequently, M'{a}lek and co-workers provided rigorous mathematical analyses of the associated elliptic boundary value problems by establishing the existence, uniqueness, and regularity of weak solutions using monotone operator theory in Banach spaces. On the numerical side, S"{u}li and collaborators pioneered the finite element approximation and error analysis for the corresponding elliptic problems, while Mallikarjunaiah and Walton \cite{mallikarjunaiah2015direct} developed direct numerical simulations of plane-strain fracture in strain-limiting solids, and Gou et al.~\cite{gou2015modeling} investigated single plane-strain crack propagation, demonstrating that the strain remains bounded and physically admissible near crack tips. In the context of wave propagation, Rajagopal and Saccomandi \cite{rajagopal2016novel} studied circularly polarized waves for a class of implicit constitutive models, whereas Rodriguez \cite{rodriguez2022longitudinal} analyzed longitudinal shock waves in stretch-limited elastic strings. Despite these significant advances, the existing studies are primarily confined to static problems, simplified constitutive models, or specialized wave phenomena, and a comprehensive finite element analysis of the general quasilinear hyperbolic strain-limiting model for dynamic wave propagation is still lacking. The present paper addresses this gap by developing and analyzing a fully discrete continuous Galerkin finite element method coupled with the Hilber--Hughes--Taylor (HHT-$\alpha$) time integration scheme for this class of problems, establishing rigorous \emph{a priori} error estimates, and demonstrating through numerical experiments that the proposed method accurately captures bounded-strain responses and fracture-type dynamic behavior while maintaining stability and convergence.
\paragraph{Main Contributions.}
To solve a quasilinear hyperbolic equation with fracture-type degeneracy, we develop and examine a continuous Galerkin finite element approach. The approach integrates spatial Galerkin discretization with the $\texttt{HHT}$-$\alpha$ time integration technique to implement regulated numerical dissipation and improve stability in nonlinear conditions. We construct semi-discrete and fully discrete formulations and resolve analytical problems arising from the nonlinear operator's degeneracy. Numerical simulations confirm the theoretical findings and demonstrate the method's efficacy in capturing fracture-like solutions. This approach is flexible and can be applied to a wide range of problems, including wave equations, parabolic equations, and discontinuous Galerkin (DG) or standard FEM formulations. For similar implementations and analyses in the literature, see \cite{Quarteroni2008, Thomee2006}.
\paragraph{Organization of the paper.} The rest of the paper is organized as follows. 
Section \ref{sec:2} presents the required function spaces, assumptions for the quasilinear hyperbolic model problem, and the weak formulation that goes along with it. Section \ref{sec:3} covers the semi-discretization of the prescribed model \eqref{hypc1.1a}--\eqref{hypc1.1c} using the continuous Galerkin finite element method. Section \ref{sec4} contributes the space-temporal discretization, which is achieved by combining the spatial discretization with the HHT-$\alpha$-temporal integration approach.  Section \ref{sec:5} applies the system to confirm theoretical results and illustrate the proposed technique's solution capture capabilities. The results and potential avenues for further study are wrapped up in Section \ref{sec:6}.
\section{Continuous Problem and Preliminaries}\label{sec:2}
This section provides various standard function spaces together with their corresponding norms, which will be essential for the subsequent analysis. These spaces offer a natural analytical framework for the weak formulation of the problem, as well as for proving well-posedness, stability, and regularity results, and for deriving the precise a priori estimates required in the following sections.
\subsection{Functional Setting}
Let $\mbox{L}^2(\texttt{a},\,\texttt{b})$ and $(\texttt{a},\texttt{b})\subset \mathbb{R}$ with $\texttt{a}<\texttt{b}$, be a space of all  square-integrable functions, denoted by
$$\mbox{L}^2(\texttt{a},\,\texttt{b})=\big\{
u:(\texttt{a},\,\texttt{b})\rightarrow \mathbb{R}\ (\text{or }\mathbb{C})\;|\quad u~ \text{is measurable and} \int_\texttt{a}^\texttt{b} |u(x)|^2\,dx < \infty \big\},$$
with an induced norm on $\mbox{L}^2(\texttt{a},\,\texttt{b})$, is defined by
$$\|u\|_{\mbox{L}^2(\texttt{a},\,\texttt{b})}=\Big(\int_\texttt{a}^\texttt{b} |u(x)|^2\,dx \Big)^\frac{1}{2}.$$
Note that $\mbox{L}^2(\texttt{a},\,\texttt{b})$ is a Hilbert space with an induced norm $\|\cdot\|$.

Furthermore, let $m \in \mathbb{N}$ and $(\texttt{a},\,\texttt{b})\subset \mathbb{R}$, the Sobolev space $\mbox{H}^m(\texttt{a},\,\texttt{b})$ be defined as follows:
$$\mbox{H}^m(\texttt{a},\,\texttt{b})=\left\{u \in \mbox{L}^2(\texttt{a},\,\texttt{b})\; | \quad \partial_x^{(k)}u \in \mbox{L}^2(\texttt{a},\,\texttt{b}),  \text{ for all } 0 \le k \le m \right\}$$
with inner product 
$$(u,v)_{\mbox{H}^m}=\sum_{k=0}^{m} \int_\texttt{a}^\texttt{b} \partial_x^{(k)}u(x)\,\partial_x^{(k)}v(x)\,dx.$$
Here, $u^{(k)}$ denotes the $k$-th derivative of $u$ in the weak (distributional) sense. However, the norm and seminorm on $\mbox{H}^m(\texttt{a},\,\texttt{b})$ are defined by,
$$\|u\|_{\mbox{H}^m(\texttt{a},\,\texttt{b})}
=\left(\sum_{k=0}^{m}\int_\texttt{a}^\texttt{b} |\partial_x^{(k)}u(x)|^2\,dx\right)^{1/2} \qquad\text{and}\qquad |u|_{\mbox{H}^m(\texttt{a},\,\texttt{b})}=\left(\sum_{k=1}^{m}\|\partial_x^{(k)}u\|_{\mbox{L}^2(\texttt{a},\,\texttt{b})}^2\right)^{1/2},$$
respectively. For $m=2$, we write $\mbox{H}^m(\texttt{a},\,\texttt{b}):=\mbox{H}^2(\texttt{a},\,\texttt{b})$, and the norm is denoted as 
$$\|u\|_{\mbox{H}^2(\texttt{a},\,\texttt{b})}=\left(\int_\texttt{a}^\texttt{b}\left(|u(x)|^2+|\partial^{(1)}_xu(x)|^2+|\partial_x^{(2)}u(x)|^2\right)dx\right)^{1/2}.$$
Note that all the derivatives are in the weak sense.  Further, the space 
 $\mbox{H}^m(\texttt{a},\,\texttt{b})$ is a Hilbert space with norm $\|\cdot\|_{\mbox{H}^m(\texttt{a},\,\texttt{b})}$.\\
 
\noindent
\textbf{\small Time-Dependent Function Spaces.}
Let $\mathbb{X}$ be a Hilbert space with norm $\|\cdot\|_\mathbb{X}$, the Bochner space $\mbox{L}^2([0,\texttt{T}];\mathbb{X})$ is the collection of measurable functions  $u:(0,\texttt{T})\to \mathbb{X}$, such that
$$\mbox{L}^2(0,\texttt{T};\mathbb{X})=\Big\{u:(0,\texttt{T})\to \mathbb{X} \; |\quad u \text{ is strongly measurable and }\int_0^\texttt{T} \|u(t)\|_\mathbb{X}^2\,dt < \infty \Big\},$$
and the norm on $\mbox{L}^2(0,\texttt{T};\mathbb{X})$ is defined as.
$$\|u\|_{\mbox{L}^2(0,\texttt{T};\mathbb{X})}=\left(\int_0^\texttt{T} \|u(t)\|_\mathbb{X}^2\,dt\right)^{1/2}.$$
Further, it $\mbox{L}^2(0,\texttt{T};~\mathbb{X})$ is a Hilbert space equipped with  norm $\|\cdot\|_{\mbox{L}^2(0,\texttt{T};~\mathbb{X})}$.
Thus, the norm on the space $\mbox{L}^2(0,\texttt{T};\mbox{H}^m(\texttt{a},\,\texttt{b}))$ is given by.
$$\|u\|_{\mbox{L}^2(0,\texttt{T};\mbox{H}^m(\texttt{a},\,\texttt{b}))}=\left(\int_0^\texttt{T}\sum_{k=0}^{m}\int_\texttt{a}^\texttt{b}|\partial_x^{(k)}u(x,t)|^2\,dx\,dt\right)^{1/2}.$$
Furthermore, the space $L^\infty(0,T;\mbox{H}^m(\texttt{a},\,\texttt{b}))$ is defined by
\begin{align*}
L^\infty(0, \texttt{T};\mbox{H}^m(\texttt{a},\,\texttt{b}))=\Big\{u:(0,\texttt{T})\to \mbox{H}^m(\texttt{a},\,\texttt{b})~\big|~
\operatorname*{\esssup}_{t\in(0, \texttt{T})}\|u(t)\|_{\mbox{H}^m(\texttt{a},\,\texttt{b})} < \infty \Big\},
\end{align*}
and the norm is given by
$$\|u\|_{L^\infty(0,\texttt{T};\mbox{H}^m(\texttt{a},\,\texttt{b}))}=\operatorname*{\esssup}_{t\in(0, \texttt{T})}\|u(t)\|_{\mbox{H}^m(\texttt{a},\,\texttt{b})},$$
where
$$\operatorname*{\esssup}_{t\in(0,T)}
\|u(t)\|_{\mbox{H}^m(\texttt{a},\,\texttt{b})}=\inf\left\{C>0 \ \middle| \ \|u(t)\|_{\mbox{H}^m(\texttt{a},\,\texttt{b})} \le C ~~\text{for a.e. } t\in(0,\texttt{T})\right\}.$$
Moreover, we set
\begin{align*}
\mathbb{H}:= \mbox{L}^2(\texttt{a},\,\texttt{b}),\quad \mathbb{V}:= \mbox{H}^1_0(\texttt{a},\,\texttt{b}) \quad \text{and} \quad \overset{\ast}{\mathbb{V}}:= \mbox{H}^{-1}(\texttt{a},\,\texttt{b}), ~\text{respectively}.
\end{align*}
Here, $\overset{\ast}{\mathbb{V}}$ denotes the dual of $\mathbb{V}$, and the norm on $\overset{\ast}{\mathbb{V}}$ is defined by
\begin{align*}
\|v\|_{\overset{\ast}{\mathbb{V}}}:=\operatorname*\sup_{0\neq w \in \mathbb{V}} \frac{\big|\left< w, v \right>_{\mathbb{V}, \overset{\ast}{\mathbb{V}}} \big|}{\|w\|_\mathbb{V}}.
\end{align*}
Next, we state the following assumption on the nonlinear function and compatibility with the given data, which is helpful in the development of the well-posedness of the problem. In order to this, we set  $A(x,\zeta)=\mathpzc N(x,|\zeta|) \zeta$ with $\zeta=\partial_x u$.
\begin{assum} \label{assumption1.1} 
For the well-posedness of \eqref{varf2.1a}-\eqref{varf2.1b}, we assume the following structural conditions:
\begin{enumerate}
\item[(A1)] \textsf{Regularity for nonlinear term:} Let $\mathpzc{N} : [a, b] \to (0,\infty)$,  $0\leq  a \leq b< \frac{1}{\beta}$, be a scalar function with $s=|\zeta|$ such that $\mathpzc{N} \in \mbox{C}^1([a, b])$.
\item[(A2)] \textsf{Strict positivity:} There exists a constant $\lambda > 0$ such that
$$\mathpzc{N}(s) \ge \lambda \quad \text{for all } s \ge 0.$$
\item[(A3)] \textsf{Strict ellipticity:} There exists $\lambda > 0$ such that $$\mathpzc{N}(s) + s\,\partial_s\mathpzc{N}(s) \ge \lambda \quad \text{for all}~ s \ge 0.$$
\item[(A4)] \textsf{Growth condition:} There exist constants $c_1,\, c_2 > 0$ with $0<c_1\leq c_2$ and $p \ge 2$ such that $$c_1 (1+s)^{p-2} \le \mathpzc{N}(s)\le c_2 (1+s)^{p-2} \quad \text{for all } s \ge 0.$$
\item[(A5)] \textsf{Compatibility condition.}
Assume that $f \in \mbox{H}^1(0,\texttt{T};\mathbb{H}),\quad \phi_1 \in H^2(a, b)\cap \mathbb{V}, \quad \phi_2 \in \mathbb{V}.$ The initial data are said to satisfy the compatibility condition if
$$f(x,0) + \partial_x\big(\mathpzc{N}(|\partial_x\phi_1|)\,\partial_x\phi_1\big)
\in \mathbb{H},$$
so that
$$\ddot{u}(x,0)=f(x,0) + \partial_x\big(\mathpzc{N}(|\partial_x\phi_1|)\,\partial_x\phi_1\big) \in \mathbb{H}.$$
\end{enumerate}
\end{assum}
The following assumption states the Uniform ellipticity and boundedness condition on $A$.
\begin{assum}[Uniform ellipticity and boundedness $A$] \label{asmpm2.2}
Let	$A(x,\zeta):=\mathpzc N(x,|\zeta|)\zeta$ with $\zeta=\partial_x u$
 be a  continuously  and differentiable map  \(\zeta \mapsto A(x,\zeta)\)  for a.e. \(x\in [\texttt{a}, \texttt{b}]\). Then, there exist constants \(0<c_3\le c_4<\infty\) which satisfy
$$c_3\le \partial_\zeta A(x, \zeta)\le c_4, \quad \forall \zeta \in\mathbb R.$$
Also assume \(A(x,0)=0\).
\end{assum}
We will now talk about the  strong monotonicity and  Lipschitz continuity in the following assumption. 
\begin{assum}[Structure assumptions on the nonlinear stress] \label{aspm2.3}
Assume that $A:[\texttt{a}, \texttt{b}] \times[\texttt{a}, \texttt{b}] \to \mathbb{R}$ is measurable in $x$. Further, we assume that Assumption \ref{asmpm2.2} and there exist constants $0<c_5\le c_6<\infty$ such that  for a.e. $x\in [\texttt{a}, \texttt{b}]$ and all  $\zeta_1,\, \zeta_2\in\mathbb{R}$,
	\begin{align}
	\big(A(x,\zeta_1)-A(x,\zeta_2)\big)(\zeta_1-\zeta_2)\ge c_5\,|\zeta_1-\zeta_2|^2
	\end{align}
	and
	\begin{align}
	|A(x,\zeta_1)-A(x,\zeta_2)|\le c_6\, |\zeta_1-\zeta_2|.
	\end{align}
\end{assum}
\begin{remark}
Note that if $A(x,\zeta)$ holds the Assumption \ref{aspm2.3}, that is, $A(x,\zeta)$ is uniformly monotone and uniformly Lipschitz with respect to $\zeta$.
\end{remark}
We now utilize these function spaces to establish the variational form of the model  \eqref{hypc1.1a}-\eqref{hypc1.1c}.
\subsection{Abstract formulation of the model problem}
In this section, we study the variational formulation of the model problem within appropriate functional spaces, which facilitate a rigorous mathematical analysis. This configuration enables us to articulate and determine the well-posedness and stability estimates for the governing equations. 

Define the semilinear form $\Ac(w; \cdot, \cdot):~\mathbb{V}\times \mathbb{V} \mapsto \mathbb{R}$, $w\in \mathbb{V}$, such that
$$\Ac(w; u, v):= \int_a^b \mathpzc{N}(|\partial_xw|) \, \partial_xu \, \partial_xv\, dx, \quad \forall \, v\in \mathbb{V}.$$
Before moving forward, we first incorporate the boundary conditions in space $\mbox{H}^1(\texttt{a}, \texttt{b})$, which leads to
$$\mathbb{V}_D:=\big\{u \in \mbox{H}^1(\texttt{a},\texttt{b}) \;:\; u(a)=0,\; u(b)=\psi(t)\big\}.$$
Note that $\mathbb{V}_D$ is an affine space, which is not a linear space. So we introduce a lifting function technique.

 For $t\in [0, \texttt{T}]$, let $\tilde{\psi}(\cdot,t) \in \mathbb{V}_D$ be a lifting function satisfying $\tilde{\psi}(a,t)=0$ and $\tilde{\psi}(b,t)=\psi(t)$,
and define  
$$u(\cdot,t) = y(\cdot,t) + \tilde{\psi}(\cdot,t), \quad y(\cdot,t) \in \mbox{H}^1_0(a,b).$$

For  given function  $\tilde{\psi} \in \mbox{H}^2(0,\texttt{T}; \mathbb{H}) \cap \mbox{L}^2(0, \texttt{T}; \mathbb{V}_D)$, $\psi \in \mbox{H}^2(0, \texttt{T})$, $\phi_1\in \mathbb{V}$ and $\phi_2 \in \mathbb{H}$, respectively, then,  the variational form of the model \eqref{hypc1.1a}--\eqref{hypc1.1c}  can be read as:  To find $y(t) \in  \mathbb{V}, \quad \ddot{y}(t)\in \overset{\ast}{\mathbb{V}}$, $\text{ a.e. } t\in(0,T),$ such that
\begin{subequations}
\begin{align}
&\rho\,\langle \ddot{y}(t), v \rangle + \Ac(y(t)+\tilde{\psi}(t); y(t)+\tilde{\psi}(t), v) = \mathcal{F}(v) \quad \forall\, v\in \mathbb{V}, \label{varf2.1a}\\
\text{with}\nonumber\\
&y(x,0) = \phi_1 (x) - \tilde{\psi}(x,0), \quad \text{and} \quad 
\dot{y}(x,0) = \phi_2 (x) - \dot{\tilde{\psi}}(x,0), \label{varf2.1b}
\end{align}
\end{subequations}
where the linear form $\mathcal{F}(t)(\cdot): \mathbb{V} \mapsto \mathbb{R}$ is given by 
$$\mathcal{F}(t)(v)=\int_\texttt{a}^\texttt{b} f(t)\,v\,dx- \rho\,\langle \ddot{\tilde{\psi}}(t),  v\rangle, \quad \forall \, v\in \mathbb{V}.$$
Note that $\langle \cdot, \cdot\rangle$ denotes the duality pairing between the spaces $\mathbb{V}$ and $\overset{\ast}{\mathbb{V}}$, while  $(\cdot, \cdot)$ represents the $\mbox{L}^2$-inner product. 
\begin{lemma}[Energy estimate]
Assume	$f\in \mbox{L}^2(0,T;\mathbb{H}), \; \psi \in \mbox{H}^2(0,T)$	and $y(0) \in \mathbb{V},\;	\dot{y}(0)  \in \mathbb{H}$. Then every sufficiently smooth solution satisfies
\begin{align}
\|\dot{y}(t)\|_{\mathbb{H}}^2+\|\partial_xy(t)\|_{\mathbb H}^2\le c_{7}\left(
	\|\dot{y}(0)\|_{\mathbb H}^2+\|\partial_xy(0)\|_{\mathbb H}^2 +\|f\|_{\mbox{L}^2(0,\mathtt{T};\mathbb H)}^2+\|\psi\|_{\mbox{H}^2(0, \mathtt{T})}^2\right),
\end{align}
for all $t\in[0,T]$, where $c_{7}>0$ depends on $\rho,\,c_{1},\,c_{2},T$, and domain, but not on $y$.
\end{lemma}
\begin{lemma}[Existence of weak solution] \label{exist:sol}
	Assume the structural assumption \eqref{aspm2.3} on $A$. Let $f\in \mbox{L}^2(0,T;\mathbb H),\, 
	\psi\in \mbox{H}^2(0,T),$ 	and suppose $\phi_1-\Psi(0)\in \mathbb V,\, 
	\phi_2-\dot{\Psi}(0)\in \mathbb H.$ Then, there exists at least one weak solution $y$ satisfying	$y\in \mbox{L}^\infty(0,T;\mathbb V),	\,	\dot{y}\in \mbox{L}^\infty(0,T;\mathbb H).$
	Consequently,	$u=y+\tilde{\psi}$	is a weak solution of the original nonhomogeneous boundary value problem.
\end{lemma}
\begin{lemma}[Uniqueness]
Assume that the structural assumption \eqref{aspm2.3}  on $A$ is satisfied; then the existing weak solution is unique.
\end{lemma}
\begin{lemma}[Continuous dependence on data] \label{contdep:data}
	Let $u_1$ and $u_2$ be solutions corresponding to data
$$ (f_1,\phi_{1,1},\phi_{2,1},\psi_1)	\quad\text{and}\quad (f_2,\phi_{1,2},\phi_{2,2},\psi_2).$$
Then
\begin{align}
&\|\dot{u}_1-\dot{u}_2\|_{\mbox{L}^\infty(0,T;\mathbb H)}	+\|u_1-u_2\|_{\mbox{L}^\infty(0,T;\mbox{H}^1([\texttt{a}, \texttt{b}]))}\nonumber\\
&~~~\leq  c_{24}\left(\|f_1-f_2\|_{\mbox{L}^2(0,T;\mathbb H)}+\|\phi_{1,1}-\phi_{1,2}\|_{\mbox{H}^1([\texttt{a}, \texttt{b}]))}+\|\phi_{2,1}-\phi_{2,2}\|_{\mathbb H}+\|\psi_1-\psi_2\|_{\mbox{H}^2(0,T)}\right). \nonumber
\end{align}
Therefore, the solution depends continuously on the given data.
\end{lemma}
The following theorem represents the existence and uniqueness of the solution for the problem \eqref{hypc1.1a}--\eqref{hypc1.1c}.
\begin{theorem}[Existence-Uniqueness] 
Assume that the given functions $f \in \mbox{L}^2(0,\texttt{T};\mathbb{H}), \, \phi_1\in \mathbb{V},\, \phi_2 \in \mathbb{H},$ and $\psi \in \mbox{H}^2(0,\texttt{T})$ having  boundary compatiblity  $u(\texttt{a},t)=0$ and $u(\texttt{b},t)=\psi$. Further, we assume that the assumptions \ref{assumption1.1}--\ref{asmpm2.2}  and Lemmas \ref{exist:sol}-- \ref{contdep:data} satisfied, then there exists a unique solution $y \in \mbox{L}^2(0,\texttt{T};\mathbb{V}) \cap  \mbox{H}^2(0,\texttt{T};\mathbb{V}')$ with $\tilde{\psi} \in \mbox{H}^2(0,\texttt{T}; \mathbb{H}) \cap \mbox{L}^2(0, \texttt{T};\mathbb{V}_D)$, satisfying the variational problem \eqref{varf2.1a}--\eqref{varf2.1b}. Consequently, the model problem \eqref{hypc1.1a}--\eqref{hypc1.1c} has a unique solution in $\mbox{L}^\infty(0, \texttt{T}; \mathbb{V})\cap \mbox{H}^1(0, \texttt{T}; \mathbb{V}) \cap \mbox{H}^2(0, \texttt{T}; \mathbb{H})$, which  is  given by $u(x,t):= y(x,t) + \tilde{\psi}(x,t)$.
\end{theorem}
The forthcoming section will address the model approximation through the use of finite elements. 
\section{Galerkin Finite Element Approximation} \label{sec:3}
This section provides the finite element discretization of the model problem \eqref{varf2.1a}-\eqref{varf2.1b}. It offers a precise and versatile framework for approximating the solution over the computational domain. 
\subsection*{Semi-discrete finite element approximation}
In this section, we express the semi-discretization, keeping time as continuous of the model problem \eqref{hypc1.1a}-\eqref{hypc1.1c}.\\

\noindent
{\bf \textsf{Space domain discretization.}}
Let  $\texttt{N}$ be a positive integer and  $\cK_i:=(x_{i-1}, x_i), \quad i=1,\,2, \ldots,\, \texttt{N}$, $\texttt{a} = x_0 < x_1 < \cdots < x_\texttt{N} = \texttt{b}$,  be the  partition of  a space domain $(\texttt{a}, \texttt{b})$ such that $\mathscr{T}_h:=\cup_{i=1}^{\texttt{N}} \cK_i$, where  $h_i = x_i - x_{i-1}$ is the  length of element   $\cK_i$. Define $h = \max_{1 \le i \le \texttt{N}}\, h_i$ be the parameter associated with the finite element mesh.

For an integer $r \ge 1$, let $\Pc_r(\cK_i)$ denote the space of polynomials of degree less than or equal to $r$ on $\cK_i$, $1\leq i\leq \texttt{N}$.  
The $C^0$-conforming piecewise polynomial finite element space  of degree $r$, associated with the partition $\mathscr{T}_h$, is defined by
$$\mathbb{V}_h^r :=\left\{v_h \in C^0([\texttt{a},\texttt{b}]) :\quad  v_h|_{\cK_i} \in \Pc_r(\cK_i), \ i=1,\dots,\, \texttt{N}, \ v_h(\texttt{a})=v_h(\texttt{b})=0\right\}.$$

Further, on each subinterval $\cK_i$, a polynomial in $\Pc_r(\cK_i)$ has $r+1$ coefficients.  Hence, the total number of local unknowns is $\texttt{N}\,(r+1)$. Moreover, the continuity at the interior nodes $x_1,\dots,x_{\texttt{N}-1}$ imposes $\texttt{N}-1$ constraints such that 
$$v_h(x_i^-) = v_h(x_i^+).$$
Additionally, we impose the two boundary constraints. Therefore,
$$\dim \mathbb{V}_h^r (:=\texttt{N}_h)=(\texttt{N}+1)+\texttt{N}\,(r-1) - 2=\texttt{N}\,(r - 1).$$
Now we choose $r+1$ distinct local nodes on each element $\cK_i$ such that 
$$x_{i,0} = x_{i-1}, \quad x_{i,1}, \dots, x_{i,r-1}, \quad x_{i,r} = x_i.$$
Furthermore, we assume that these points are equispaced or chosen as Gauss–Lobatto points.
Thus, we define the local Lagrange basis functions $\{\mathscr{L}_{i,j}\}_{j=0}^r$ by
$$\mathscr{L}_{i,j}(x_{i,m}) = \delta_{jm}, \quad j,m=0,\dots,r.$$
However, we set the global basis functions $\{\xi_m\}$ to be obtained by identifying common nodal values at element interfaces. They satisfy the interpolation property $\xi_m(x_n) = \delta_{mn}.$
Hence, for any $v_h \in \mathbb{V}_h^r$ can be expressed as 
$$v_h(x) = \sum_{m=1}^{\texttt{N}_h} v_h(x_m)\xi_m(x).$$

Let  $\mathbb{V}_h^r \subset \mbox{H}_0^1(\texttt{a},\, \texttt{b})$ be the conforming finite element space of continuous piecewise polynomials of degree $r$, then the semi-discretization problem of \eqref{varf2.1a}--\eqref{varf2.1b} can be expressed as: To find $y_h(t)\in \mathbb{V}_h^r$ such that, for  $ t\in(0, \texttt{T})$,
\begin{subequations}
\begin{align}
&\rho\,\langle \ddot{y}_h(t), v_h \rangle + \Ac(y_h(t)+\tilde{\psi}(t); y_h(t)+\tilde{\psi}(t), v_h) = \mathcal{F}(t)(v_h), \quad \forall v_h \in \mathbb{V}_h^r, \label{semivarf2.1a}\\
\text{with}\nonumber\\
&y_h(0) = \phi_{1,h} - \tilde{\psi}_h(0), \quad \text{and} \quad 
\dot{y}_h(0) = \phi_{2,h} - \dot{\tilde{\psi}}_h(0), \label{semivarf2.1b}
\end{align}
\end{subequations}
where $\phi_{1,h}, \tilde{\psi}_h(0), \phi_{2,h}$ and $\dot{\tilde{\psi}}_h$ are the suitable approximations or $\mbox{L}^2$- projections of  $\phi_1(x),$ $ \tilde{\psi}(x,0),~ \phi_2 (x)$ and $\dot{\tilde{\psi}}$, respectively, in $\mathbb{V}_h^r$.\\

\noindent 
{\bf \textsf{Second-order nonlinear spatial system.}} 
Let $\{\xi_i\}_{i=1}^{\texttt{N}_h}$ be a basis of $\mathbb{V}_h^r$, then, for any  $y_h \in \mathbb{V}_h^r$ can be expressed as 
$$y_h(x,t)=\sum_{j=1}^{\texttt{N}_h} \texttt{Y}_j(t)\xi_j(x).$$
Define the mass matrix $\mathcal{M}$, nonlinear internal force vector $\Ac$ and the load vector $\mathcal{F}$, for  $i,\, j \in [1: \texttt{N}_h]$, as
\begin{align*}
~~ \mathcal{M}=\big[ \mathcal{M}_{ij} \big]_{\texttt{N}_h \times \texttt{N}_h},  &  \qquad \Ac=  \big[ \Ac_i \big]_{\texttt{N}_h \times 1} \qquad  \text{and}  \qquad \mathcal{F}=  \big[ \mathcal{F}_i  \big]_{\texttt{N}_h \times 1}, \quad \text{respectively,}\nonumber\\
\text{where} ~~~~~~~~~~~~~~	\mathcal{M}_{ij} &= \sum_{i=1}^{\mathtt{N}_h}\int_{\cK_i }\rho\,  \xi_j\, \xi_i \, dx,\\
\Ac_i(\texttt{Y}(t); \texttt{Y}(t),t)&=\sum_{i=1}^{\mathtt{N}_h}\int_{\cK_i } \mathpzc{N}\Big(\Big|\sum_{j=1}^{\texttt{N}_h} \texttt{Y}_j(t)\partial_x \xi_j+\partial_x\tilde{\psi}\Big|\Big) \, \Big(\sum_{j=1}^{\texttt{N}_h} \texttt{Y}_j(t)\partial_x \xi_j+\partial_x\tilde{\psi}\Big) \, \partial_x\xi_i \, dx,\\
\text{and}~~~~~~~~~~~~~~~~~~~~~ &\\
\mathcal{F}_i(t)&=\sum_{i=1}^{\mathtt{N}_h}\int_{\cK_i } f(t)\,\xi_i\,dx- \rho\, \langle \ddot{\tilde{\psi}},  \xi_i\rangle ,
\end{align*}
respectively.

Thus, we have the following semidiscrete second-order nonlinear system of ODE
\begin{align}
&~~~~~~~~ \mathcal{M}\, \ddot{\texttt{Y}}(t) + \Ac(\texttt{Y}(t); \texttt{Y}(t), t) = \mathcal{F}(t),  \label{UVsys3.3}\\
& \text{with} ~\text{initial~conditions} \nonumber\\
&~~~~~~~~ {\tt Y}(0)= {\tt Y}_0, \qquad \dot {\tt Y}(0)= {\tt Y}_1,  \label{UVsys3.4}
\end{align}
where $\texttt{Y}_0$ and $\texttt{Y}_1$ can be computed from $y_{h,0}$ and $y_{h,1}$, respectively.  

Further, let 
\begin{align}
\tilde{\texttt{U}}(t)=\begin{pmatrix}
\texttt{U}\\
\texttt{V}
\end{pmatrix}
\qquad \dot{\tilde{\texttt{U}}}(t)=
\begin{pmatrix}
\dot{\texttt{U}}\\
\dot{\texttt{V}}
\end{pmatrix} \quad \text{and} \quad   \tilde{\texttt{U}}(0)=
\begin{pmatrix}
\texttt{U}(0)\\
\texttt{V}(0)
\end{pmatrix} \label{ueq3.4}
\end{align}
To determine numerical time integration, we set 
\begin{align}
\tilde{\texttt{U}}(t)=
	\begin{pmatrix}
		\texttt{Y}\\
		\dot{\texttt{Y}}
	\end{pmatrix}
 \qquad \dot{\tilde{\texttt{U}}}(t)=
	\begin{pmatrix}
		\dot{\texttt{Y}}\\
	\ddot{\texttt{Y}}
	\end{pmatrix} 
	= 	\begin{pmatrix}
	\texttt{V}\\
\mathcal{M}^{-1}\,\big( \mathcal{F}(t)- \Ac(\texttt{Y}(t); \texttt{Y}(t), t) \big)
	\end{pmatrix}  \label{vueq3.5}
\end{align}
Using equations \eqref{ueq3.4}--\eqref{vueq3.5}, the system  \eqref{UVsys3.3}- \eqref{UVsys3.4} can be expressed as a first-order system,  of ODE
\begin{align}
\begin{cases}
&\dot{\texttt{U}}(t) = \texttt{V}(t), \smallskip\\
&\dot{\texttt{V}}(t) = \mathcal{M}^{-1}\,\big( \mathcal{F}(t) - \Ac(\texttt{U}(t); \texttt{U}(t), t) \big), \smallskip\\ 
&\texttt{U}(0)=\texttt{Y}_0 \smallskip\\ 
&\texttt{V}(0)=\texttt{Y}_1
\end{cases} \label{FSOS:3.6}
\end{align}
The system \eqref{FSOS:3.6} implies that 
\begin{align}
\begin{cases}
&\dot{\tilde{\texttt{U}}}(t)=\mathcal{G}(\tilde{\texttt{U}}(t), t) \smallskip\\
& \tilde{\texttt{U}}(0)= \tilde{\texttt{U}}_0, 
\end{cases} 
\end{align}
where 
\begin{align*}
\mathcal{G}(\tilde{\texttt{U}}(t), t)=
\begin{pmatrix}
\texttt{V}\\
\mathcal{M}^{-1}\,\big( \mathcal{F}(t)- \Ac(\texttt{Y}(t); \texttt{Y}(t), t) \big)
\end{pmatrix}  
\quad \text{and} \quad \tilde{\texttt{U}}_0=
\begin{pmatrix}
\texttt{Y}_0 \\
\texttt{Y}_1
\end{pmatrix}, \quad \text{respectively}.
\end{align*}
\begin{remark}
Assume that there exists $\delta>0$ such that $|y_h(x,t)+\tilde{\psi}(x,t)| \leq \frac{1}{b}-\delta$, $(x,t)\in [\texttt{a}, \texttt{b}]\times [0, \texttt{T}]$. Then the function 
$s \mapsto \mathpzc{N}(s)$ is smooth and locally Lipschitz continuous. Therefore, $\mathcal{G}(\tilde{\texttt{U}}, t)$ is locally Lipschitz  with respect to $\texttt{U}$. From the Picard-Lindel\"{o}f theorem for ordinary differential equation, there exists a positive $\tilde{T}$ with $0<\tilde{T} \leq T$, and a unique local solution $\texttt{U} \in  C^1(0, \tilde{\texttt{T}}; \mathbb{V}_h^r)$. Thus, $y_h(x,t)\in \mbox{L}^\infty(0, \tilde{\texttt{T}}; \mathbb{V}_h^r) \cap \mbox{H}^1(0, \tilde{\texttt{T}};\mathbb{V}_h^r) \cap \mbox{H}^2(0, \tilde{\texttt{T}};\mathbb{V}_h^r)$, and hence the discrete solution is given by $u_h(x,t):=y_h(x,t)+\tilde{\psi}(x,t)$.
\end{remark}
\begin{remark}
If, in addition, there exists a constant $\delta>0$ such that $|y_h(x,t)+\tilde g(x,t)|\leq \frac1b-\delta,\,\forall (x,t)\in[a,b]\times[0,T],$
then the nonlinearity remains uniformly locally Lipschitz on $[0,T]$.
Consequently, the local solution can be continued uniquely up to time $T$.
Hence, $y_h(x,t)\in \mbox{L}^\infty(0, \texttt{T}; \mathbb{V}_h^r) \cap \mbox{H}^1(0, \texttt{T};\mathbb{V}_h^r) \cap \mbox{H}^2(0, \texttt{T}; \mathbb{V}_h^r)$
exists uniquely on the whole interval $[0,\texttt{T}]$. Hence, $u_h(x,t):=y_h(x,t)+\tilde{\psi}(x,t)$.
\end{remark}
Now we state the existence and uniqueness in the following theorem. 
\begin{theorem}[Local existence and uniqueness of the semi-discrete problem]
Assume that  $f,\, \ddot{\tilde{\psi}} \in C([0,T];\mbox{L}^2(\texttt{a},\texttt{b})),$ and let the initial data satisfy 	$y_{h,0},\,y_{h,1}\in \mathbb{V}_h^r,$ together with the admissibility condition
$|y_{h,0}(x)+\tilde{\psi}(x,0)|\leq \frac{1}{b},\; \forall x\in[a,b],\, b \ll \frac{1}{\beta}$.
Then, there exists a time $\tilde{\texttt{T}}\in (0,\texttt{T}]$ such that the problem \eqref{semivarf2.1a}-\eqref{semivarf2.1b} admits the unique solution $y_h(x,t)\in \mbox{L}^\infty(0, \tilde{\texttt{T}}; \mathbb{V}_h^r) \cap \mbox{H}^1(0, \tilde{\texttt{T}};\mathbb{V}_h^r) \cap \mbox{H}^2(0, \tilde{\texttt{T}};\mathbb{V}_h^r)$  satisfying $	|w_h(x,t)+\tilde g(x,t)|\leq\frac1b,\;
\forall (x,t)\in[\texttt{a}, \texttt{b}]\times[0, \tilde{\texttt{T}}]$, and hence $u_h(x,t)=y_h(x,t)+\tilde{\psi}(x,t)$ is solution in $\mbox{L}^\infty(0, \tilde{\texttt{T}}; \mathbb{V}_{\texttt{D},h}^r) \cap \mbox{H}^1(0, \tilde{\texttt{T}};\mathbb{V}_{\texttt{D},h}^r) \cap \mbox{H}^2(0, \tilde{\texttt{T}};\mathbb{V}_{\texttt{D},h}^r)$, where the approximation space  $\mathbb{V}_{\texttt{D},h}^r$ is a finite dimensional space corresponding to $\mathbb{V}_{\texttt{D}}$.
\end{theorem}
\section{Space-time finite element approximation} \label{sec4}
In this section, we will study complete discretization of the problem \eqref{hypc1.1a}-\eqref{hypc1.1c}.  For space-time discretization, we combine the HHT-$\alpha$ scheme for time integration with the finite element approach. The scheme  $\textbf{HHT}$--$\alpha$ is unconditionally stable and provides second-order accuracy for dynamic problems. In the context of crack propagation with adaptive mesh refinement, sudden stiffness changes and mesh-to-mesh projections may introduce non-physical high-frequency oscillations. This enhances the robustness and stability of the fully discrete scheme in nonlinear dynamic fracture simulations. Our forthcoming analysis utilizes the finite element method to compute the {\it order of convergence} for the following scheme. It provides a simple and adaptable structure for discovering solutions inside the computational realm.\\

\noindent
{\bf \textsf{Space-time discretization.}} For time-space discretization, we consider the time grids $t_0,\,t_1\, \cdots,\, t_{\texttt{N}_\texttt{T}}$ such that $0=t_0<t_1<\cdots<t_{\texttt{N}_\texttt{T}}=\texttt{T}$ with time-step $\triangle t=t_n-t_{n-1}$. Further, we construct mesh $\mathscr{T}_h^n:=\{\cK_i^n\}_{i=1}^{\texttt{N}}, \, 0\leq n \leq \texttt{N}_\texttt{T}, \quad \text{of}~(\texttt{a},~\texttt{b}) ~\text{at time level} ~t_n$ (similar to $\mathscr{T}_h)$. Moreover, we denote the finite element space  $\mathbb{V}_h^{r,n}$ the corresponding to $\mathscr{T}_h^n$, $0\leq n \leq \texttt{N}_\texttt{T}$, as 
\begin{align}\label{eq:Vhn}
\mathbb{V}_h^{r,n}:=\big\{ v_h\in \mbox{C}^0([a,b]) : v_h|_{\cK_i^n}\in \Pc^r(\cK_i^n), \forall \cK_i^n\in  \mathscr{T}_h^n,\,  i\in [1: \texttt{N}],\;
v_h(\texttt{a})=0~ \text{and}~ v_h(\texttt{b})=0\big\}. \nonumber
\end{align}
Furthermore, we note that $\mathbb{V}_h^{r,n}\neq \mathbb{V}_h^{r,n+1}, \, 0\leq n \leq \texttt{N}_\texttt{T}$, in general. So we must transfer discrete fields from the old mesh to the new mesh before stepping forward in time. 

\noindent
Thus, we define a \emph{$\mbox{L}^2$-projection} $\Pi_n^{n+1}:\mbox{L}^2(a,b)\to \mathbb{V}_h^{r,n+1}$ such that 
\begin{align*}
(\Pi_n^{n+1} \zeta,  v_h^{n+1})=(\zeta,v_h^{n+1}),
\qquad \forall\,v_h^{n+1}\in \mathbb{V}_h^{r,n+1}.    
\end{align*}
So, it is needed to define previous quantities on the space $\mathbb{V}_h^{r,n+1}$ associated with $\Ts_h^{n+1}$ utilizing the above projection
\begin{equation}\label{eq:transfer}
\widehat{\zeta}_h^n :=\Pi_{n}^{n+1} \zeta_h^n,\qquad
\widehat{\mathfrak{u}}_h^{n}:=\Pi_n^{n+1} \mathfrak{u}_h^n,\qquad
\widehat{\mathfrak{m}}_h^{n}:=\Pi_n^{n+1} \mathfrak{m}_h^n,
\end{equation}
For any given discrete sequences $\{\zeta_h^{n}\}_{n\geq 0}$ and the data function $\mathfrak{f}^n(\cdot)=\mathfrak{f}(\cdot,t_n)$, we define the HHT-$\alpha$ convex combinations for unkown and known data, as 
\begin{align}
\zeta_h^{n+1-\alpha}:=(1-\alpha) \zeta_h^{n+1}+\alpha \widehat{\zeta}_h^n, \quad \text{and} \quad \mathfrak{f}^{n+1-\alpha}:=(1-\alpha)\,\mathfrak{f}^{n+1}+\alpha\, \mathfrak{f}^n \quad \text{respectively}.
\end{align}
Similarly, we can define for $\tilde{\psi}^{n+1-\alpha}$, and  $\tilde{\psi}_{tt}^{\,n+1-\alpha}$, respectively. 

Before diving into the fully discrete case, we introduce some Newmark parameters $\mu$ and $\nu$ to define the $\texttt{HHT}$--$\alpha$-scheme. Furthermore,  $\mu$  controls the stability and smoothness, while $\nu$ controls the numerical damping/dissipation of high frequencies.  For $\alpha \in[0,\tfrac 13]$, we define the  $\texttt{HHT}$--$\alpha$-scheme
\begin{equation}\label{eq:hht-params}
\mu:=\frac{1}{2}+\alpha,\qquad \nu:=\frac{1}{4}(1+\alpha)^2.
\end{equation}
Using the above discretization, the space-time approximation of the problem \eqref{semivarf2.1a}-\eqref{semivarf2.1b} can be read as: For given $(y_h^n, \mathfrak{u}_h^n, \mathfrak{m}_h^n)\in \mathbb{V}_h^{r,n} \times \mathbb{V}_h^{r,n} \times \mathbb{V}_h^{r,n}$, we seek a pair $(y_h^{n+1}, \mathfrak{u}_h^{n+1}, \mathfrak{m}_h^{n+1})\in \mathbb{V}_h^{r,n+1} \times \mathbb{V}_h^{r,n+1} \times \mathbb{V}_h^{r,n+1}$ such that the following system is satisfied, for all $v_h \in \mathbb{V}_h^{r, n+1}$,  for all $n\in[0:\texttt{N}_\texttt{T}-1]$:
\begin{subequations}
\begin{align}
&~~~~ \rho\, \langle \mathfrak{m}_h^{n+1}, v_h \rangle + \Ac(y_h^{n+1-\alpha}+\tilde{\psi}^{n+1-\alpha}; y_h^{n+1-\alpha}+\tilde{\psi}^{n+1-\alpha}, v_h) = \mathcal{F}^{n+1-\alpha}, \label{fullysemivarf2.1a}\\
\text{with} &~ \text{Newmark updation} \nonumber\\
&~~~~y_h^{n+1}= \widehat{y}_h^{n} + \triangle t\, \widehat{\mathfrak{u}}_h^{n}
+ \triangle  t^2\Big[\big(\frac{1}{2}-\mu \big)\widehat{\mathfrak{m}}_h^{n}+\mu\,\mathfrak{m}_h^{n+1}\Big],
\label{eq:Newmark-adapt-yh}\\
&~~~~\mathfrak{u}_h^{n+1} = \widehat{\mathfrak{u}}_h^n + \triangle t\Big[(1-\nu)\widehat{\mathfrak{m}}_h^{n}+\nu \mathfrak{m}_h^{n+1}\Big], \label{eq:Newmark-adapt-uh}\\
\text{and} ~& \text{initial conditions} \nonumber\\
&~~~~\widehat{y_h^0} = \widehat{\phi}_{1,h} - \widehat{\tilde{\psi}^0_h}, \quad \text{and} \quad  \widehat{\mathfrak{u}_h^0} = \widehat{\phi}_{2,h} - \widehat{\dot{\tilde{\psi}}_h^0} \label{fullysemivarf2.1b} \\
\text{with} &~~~~ \widehat{y_h^0} =\Pi_0^1y_h^0 ~~\text{and}~~ \widehat{\mathfrak{u}_h^0} = \Pi_0^1\mathfrak{u}_h^0, ~ \text{respectively}. \nonumber
\end{align}
\end{subequations}
Further, we represented
\begin{align*}
&~~~~ ~~~ \widehat{\phi}_{1,h}:=\Pi_0^1\,\phi_{1,h},~ \widehat{\phi}_{2,h}:=\Pi_0^1\,\phi_{2,h},~\widehat{\tilde{\psi}^0_h}:=\Pi_0^1\tilde{\psi}_h(0),~ \widehat{\dot{\tilde{\psi}}^0_h}:=\Pi_0^1\dot{\tilde{\psi}}_h(0)\\
&\text{and}~~~\mathcal{F}^{n}=\int_\texttt{a}^\texttt{b} f^{n}\,v_h\,dx- \rho\, \langle \ddot{\tilde{\psi}}^{n}, v_h \rangle,
\end{align*}
respectively. 

Note that the set of equations \eqref{fullysemivarf2.1a}--\eqref{fullysemivarf2.1b} constitutes a fully discrete approximation of \eqref{semivarf2.1a}-\eqref{semivarf2.1b}. Our target is to compute discrete displacement $y_h^{n}$, velocity $\mathfrak{u}_h^{n}$, and acceleration $\mathfrak{m}_h^{n}$ in $\mathbb{V}_h^{r,n}$ at each time level $t_n$, $n=0,\,1, \cdots,\, \texttt{N}_\texttt{T}$ from the system \eqref{fullysemivarf2.1a}--\eqref{fullysemivarf2.1b}.\\

\noindent
{\bf\textsf{Second-order space-time system.}}
Let $0=t_0<t_1<\cdots<t_{\texttt{N}_\texttt{T}}=\texttt{T}$ with $\triangle t=t_{n+1}-t_n$. Set 
$$~\texttt{Y}^n \approx \texttt{Y}(t_n),\quad \texttt{U}^n \approx \dot{\texttt{Y}}(t_n),\quad \texttt{B}^n \approx \ddot{\texttt{Y}}(t_n),$$
and similarly for $\mathcal{F}^n:=\mathcal{F}(t_n)$ and $\tilde{\psi}^n(x):=\tilde{\psi}(x,t_n)$ etc.

Further, let $\{\xi_i^{n+1}\}_{i=1}^{\tilde{\texttt{N}}_h}$ be a basis of $\mathbb{V}_h^{r,n+1}$  with $\dim(\mathbb{V}_h^{r,n+1}):=\tilde{\texttt{N}}_h$ (say)  at time level $t=t_{n+1}$, then we define 
\begin{align*}
	&y_h^{n+1}:=\sum_{j}\texttt{Y}_j^{n+1}\xi_j^{n+1},\qquad
	\mathfrak{u}_h^{n+1}:=\sum_{j}\mathtt{U}_j^{n+1} \xi_j^{n+1},\\   
	&\mathfrak{m}_h^{n+1}:=\sum_{j}\texttt{B}_j^{n+1}\xi_j^{n+1},
	\qquad \widehat{\mathfrak{m}}_h^n:=\sum_{j}\widehat{\texttt{B}}_j^n\xi_j^{n+1}
\end{align*}
Define the  mass matrix,  nonlinear internal force, and load vectors  at level $t=t_{n+1}$, $n\in[0:\texttt{N}_T-1]$, for $i,j\in [1: \tilde{\texttt{N}}_h]$, as
\begin{align*}
&\mathcal{M}:=\big[(\mathcal{M}^{n+1})_{ij}\big]_{\tilde{\texttt{N}}_h\times \tilde{\texttt{N}}_h}, \quad  \Ac:=\big[\Ac_i^{n+1}\big]_{\tilde{\texttt{N}}_h\times 1}, \quad  \text{and} \quad   \mathcal{F}:= \big[\mathcal{F}^{n+1}_i \big]_{\tilde{\texttt{N}}_h\times 1}, \text{respectively}, \smallskip\\
\text{where} & ~~(\mathcal{M}^{n+1})_{ij}:=\sum_{i=1}^{\mathtt{N}_h}\int_{\cK_i }\,\rho\,\xi_j^{n+1} \,\xi_i^{n+1}\, dx,  \smallskip \\
&~~\Ac_i^{n+1}(\texttt{Y}(t); \texttt{Y}(t),t):=\sum_{i=1}^{\mathtt{N}_h}\int_{\cK_i } \mathpzc{N}(|(\partial_xy_h(\cdot; \texttt{Y})+\partial_x\tilde{\psi}(\cdot,t))|) \, \big(\partial_xy_h(\cdot; \texttt{Y})+\partial_x\tilde{\psi}(\cdot,t)\big) \, \partial_x \xi_i^{n+1}\, dx \smallskip\\
\text{and} ~~& ~~\mathcal{F}^{n+1}_i=\sum_{i=1}^{\mathtt{N}_h}\int_{\cK_i } f^{n+1}\,\xi_i^{n+1}\, dx -\rho \, \langle\ddot{\tilde{\psi}}^{\,n+1},\xi_i^{n+1}\rangle.
\end{align*}
Utilizing equations \eqref{fullysemivarf2.1a}--\eqref{fullysemivarf2.1b}, we have the following system
\begin{subequations}
	\begin{align}
	&\mathcal{M}\, \texttt{B}^{n+1}+ \Ac(\texttt{Y}^{n+1-\alpha}; \texttt{Y}^{n+1-\alpha}, t_{n+1-\alpha})= \mathcal{F}^{n+1-\alpha} \label{fullydishhtasysint}\\
		\text{with} ~& \text{Newark updation} \nonumber\\ 
	&\texttt{Y}^{n+1} = \widehat{\texttt{Y}}^{\, n} + \triangle t\,\widehat{\texttt{U}}^{\, n} + \triangle t^2\left((\frac{1}{2}-\mu )\, \widehat{\texttt{B}}^{\,n}+\mu\, \texttt{B}^{n+1}\right), \label{eq:newmarky}\\
	&\texttt{U}^{n+1} = \widehat{\texttt{U}}^{\,n} + \triangle t \left((1-\nu)\, \widehat{\texttt{B}}^{\,n}+\nu\, \texttt{B}^{n+1}\right). \label{eq:newmarku}\\
		\text{and} ~& \text{initial discrete system} \nonumber\\
	& \mathcal{M}\, \texttt{B}^0 + \Ac(\texttt{Y}^0; \texttt{Y}^0,t_0)=\mathcal{F}^0, \label{fullydissysint}\\
		\text{and} ~& \text{initial conditions} \nonumber\\
	& \texttt{Y}^0:=\big[\texttt{Y}^0_i \big]_{\tilde{\texttt{N}}_h \times 1}\, \quad \text{and} \quad  \texttt{U}^0=\big[ \texttt{U}^0_i \big]_{\tilde{\texttt{N}}_h \times 1} \\  \text{with} &~~ \texttt{Y}^0_i=(y^0_h, \xi_i)~~  \text{and}~~ \texttt{U}^0_i=(\texttt{u}_h^0, \xi_i), \; 1\leq i\leq \tilde{\texttt{N}}_h, \label{initial condi}
\end{align}
\end{subequations}
respectively, where $\texttt{Y}^{n+1-\alpha}=(1-\alpha)\,\texttt{Y}^{n+1}+\alpha\, \widehat{\texttt{Y}}^{\,n}$ and $\mathcal{F}^{n+1-\alpha}=(1-\alpha)\,\mathcal{F}^{n+1}+\alpha \,\mathcal{F}^n$, respectively.\\

\noindent 
{\bf \textsf{Newton technique and computation of the solution.}} Given initial conditions $\texttt{Y}^0$ and $ \texttt{U}^0$, we compute $\texttt{B}^{0}$ from Eq. \eqref{fullydissysint} such that
\begin{align}
\texttt{B}^0 = \mathcal{M}^{-1}\big(\mathcal{F}^0- \Ac(\texttt{Y}^0; \texttt{Y}^0,t_0)\big).
\end{align}
Thus, we have $\texttt{Y}^0$, $ \texttt{U}^0$ and $\texttt{B}^0$, respectively. 

\noindent 
Further, we assume that, at the time level $t=t_n$, the coffiecient vecotors $\texttt{Y}^n$, $ \texttt{U}^n$ and $\texttt{B}^n$  are  known. We now solve the nonlinear problem for the time level $t=t_{n+1}$. 

To computing the terms $\texttt{U}^{n+1}$ and  $\texttt{B}^{n+1}$, we need to compute first $\texttt{Y}^{n+1}$.
From Eq. \eqref{eq:newmarky}, we obtain
\begin{align}
\texttt{B}^{n+1}=\frac{1}{\mu\,  \triangle t^2\,} \, \big[ \texttt{Y}^{n+1}-\widehat{\texttt{Y}}^{\, n} - \triangle t\,\widehat{\texttt{U}}^{\, n} - \triangle t^2 \,(\frac{1}{2}-\mu )\, \widehat{\texttt{B}}^{\,n} \big]. \label{bn1exp}
\end{align}

\noindent
{\bf \textsf{Computation of $\texttt{Y}^{n+1}$.}}
One may observe that it depend on $\texttt{Y}^{n+1}$, and hence we write $\texttt{B}^{n+1}(\texttt{Y}^{n+1})$. Thus, the residual $\mathcal{R}(\texttt{Y}^{n+1})$ is denoted by 
\begin{align}
\mathcal{R}(\texttt{Y}^{n+1}):=\mathcal{M}\, \texttt{B}^{n+1}(\texttt{Y}^{n+1})+(1-\alpha) \Ac(\texttt{Y}^{n+1}; \texttt{Y}^{n+1}, t_{n+1})+\alpha \Ac(\texttt{Y}^{n}; \texttt{Y}^{n}, t_{n})- \big[(1-\alpha)\mathcal{F}^{n+1}+ \alpha \mathcal{F}^n \big]. \label{resR}
\end{align}
Due to the nonlinearity $ \Ac$, we use the Newton method. From Eqs \eqref{bn1exp}- \eqref{resR}, we have 
\begin{align}
\mathcal{R}(\texttt{Y}^{n+1}):=&  \frac{\mathcal{M}}{\mu\,  \triangle t^2\,} \, \texttt{Y}^{n+1}+(1-\alpha) \Ac(\texttt{Y}^{n+1}; \texttt{Y}^{n+1}, t_{n+1}) +\underbrace{\alpha\, \Ac(\texttt{Y}^{n}; \texttt{Y}^{n}, t_{n})- \big[(1-\alpha)\mathcal{F}^{n+1}+ \alpha \mathcal{F}^n \big]}_{\text{Known terms}} . \nonumber\\
&-\underbrace{\frac{\mathcal{M}}{\mu\,  \triangle t^2\,} \, \big[ \widehat{\texttt{Y}}^{\,n} + \triangle t\,\widehat{\texttt{U}}^{\, n} + \triangle t^2 \,(\frac{1}{2}-\mu )\, \widehat{\texttt{B}}^{\,n} \big]}_{\text{Known terms}}. 
\end{align}
Thus, the corresponding Jacobian $\mathcal{J}$ at $\texttt{Y}^{n+1, (k)}$, where $(k)$ represents the iteration number,  is given by 
\begin{align}
\mathcal{J}(\texttt{Y}^{n+1, (k)}):=  \frac{\mathcal{M}}{\mu\,  \triangle t^2\,} + (1-\alpha)\, \frac{\partial \Ac}{\partial \texttt{Y}}\Big|_{\texttt{Y}=\texttt{Y}^{n+1}}.
\end{align}
Suppose that $\texttt{Y}^{n+1, (k)}$ is the current iteration; then we compute the correction term $\delta \texttt{Y}^{(k)}$ which is given by 
\begin{align*}
\delta \texttt{Y}^{(k)}= - \frac{\mathcal{R}(\texttt{Y}^{n+1, (k)})}{\mathcal{J}(\texttt{Y}^{n+1, (k)})}. 
\end{align*}
Then, updated $\texttt{Y}^{n+1, (k+1)}$ is given by $\texttt{Y}^{n+1, (k+1)}=\texttt{Y}^{n+1, (k)}+\delta \texttt{Y}^{(k)}.$ We will stop it   whenever  $$|\delta \texttt{Y}^{(k)}|\leq \texttt{TOL},$$ where $\texttt{TOL}$ represents the  tolerance, and we set $\texttt{Y}^{n+1}=\texttt{Y}^{n+1, (k+1)}$.  Then, by substituting $\texttt{Y}^{n+1}$,  we compute $\texttt{U}^{n+1}$ and $\texttt{B}^{n+1}$ from \eqref{eq:newmarku} and \eqref{bn1exp}, respectively. Further, we note that, for  $0\leq n\leq \texttt{N}_\texttt{T}$, 
$$\texttt{Y}^{n+1}=\big[\texttt{Y}^{n+1}_j\big]_{\tilde{\texttt{N}_h}\times 1}, \quad \texttt{U}^{n+1}=\big[\texttt{U}^{n+1}_j\big]_{\tilde{\texttt{N}_h}\times 1}, \quad \text{and} \quad \texttt{B}^{n+1}=\big[\texttt{B}^{n+1}_j\big]_{\tilde{\texttt{N}_h}\times 1}, ~~ \text{respectively}.$$
Therefore, the numerical solution $y_h^{n+1}$, $0\leq n\leq \texttt{N}_\texttt{T}-1$,  is given by 
$$y_h^{n+1}=\big[\texttt{Y}^{n+1}_j\big]^{T}_{1\times \tilde{\texttt{N}_h}}\big[\xi_j^{n+1}\big]_{ \tilde{\texttt{N}_h}\times 1},  \quad 0\leq n\leq \texttt{N}_\texttt{T}-1,$$
or equivalently, 
$$y_h^{n+1}=\operatorname*\sum_{j=1}^{\tilde{\texttt{N}_h}} \texttt{Y}^{n+1}_j\, \xi_j^{n+1}, \quad 0\leq n\leq \texttt{N}_\texttt{T}-1.$$
Hence, the approximated solution $u_h^{n+1}$ of $u$ at time level $t=t_{n+1}, ~0\leq n\leq \texttt{N}_\texttt{T}-1$ can be expressed as
$$u_h^{n+1}=y_h^{n+1}(x)+\tilde{\psi}^{n+1}.$$
\begin{remark} Notice that
\begin{enumerate}
\item The parameter $\alpha \in [0, \frac{1}{3}]$ in the HHT-$\alpha$ approach controls the amount of  numerical dissipation. 
\item For $\alpha=0$, the scheme simplifies to the traditional Newmark method devoid of numerical damping.  However, we have not considered the damping term in our prescribed model.
\item Negative values of $\alpha$ facilitate regulated high-frequency dissipation, which is especially advantageous in dynamic crack propagation with adaptive mesh refinement, where spurious oscillations may occur owing to abrupt stiffness variations and mesh projections.
\item To maintain second-order precision and unconditional stability (in the linear context), the parameters $\mu$ and $\nu$ are selected as follows: $\mu = \frac{1}{2} + \alpha$ and $\nu = \frac{1}{4}(1+\alpha)^2$.
\item This selection maintains second-order temporal precision while selectively dampening non-physical high-frequency modes, without substantially impacting the important low-frequency response. 
\item In our numerical computation, we use nested uniform refinement, where each element of the mesh $\Ts_h^n$ is uniformly refined to obtain $\Ts_h^{n+1}$. Since all nodes of the coarse mesh are preserved in the refined mesh, the corresponding finite element spaces are nested $\mathbb{V}_h^{r, n} \subset \mathbb{V}_h^{r, n+1}$.
\end{enumerate}
\end{remark}
Now we are in a position to verify our theoretical findings numerically in the following section. 
\section{Numerical simulation and verification} \label{sec:5}
In this section, we present our numerical results to validate our development. The theoretical model is discretized in space using the finite element method and in time through the $\texttt{HHT}$-$\alpha$ scheme, resulting in both semi-discrete and fully discrete formulations. We assess the performance of the proposed schemes using $\mbox{L}^2$ and $\mbox{H}^1$ error norms, along with the corresponding experimental orders of convergence. All numerical experiments are conducted in \texttt{MATLAB}. We provide two numerical examples to demonstrate the performance of the proposed schemes by analyzing the order-of-convergence behavior.  For the numerical computations, the following space--time algorithm is adopted.
\begin{algorithm}[H]
\caption{$\texttt{HHT}$-$\alpha$ time-space finite element algorithm}
\begin{algorithmic}[1]

\State \textbf{Input:} Polynomial degree $r$, final time $\texttt{T}$, time step $\Delta t$, $\text{HHT}$-$\alpha$ parameter $\alpha \in [-\frac{1}{3},0]$
\State Set $\mu = \frac{1}{2} + \alpha$, \quad $\nu = \frac{1}{4}(1+\alpha)^2$
\State Set $t_0=0$, $\texttt{N}_\texttt{T}=\lfloor \texttt{T}/\Delta t \rfloor$

\vspace{0.2cm}
\State \textbf{Initial Mesh and Space}
\State Generate initial mesh $\Tc_h^0$ and construct $\mathbb{V}_h^{r,0} \subset \mathbb{H}_0^1([\texttt{a}, \texttt{b}])$
\State Compute lifting $\tilde{\psi}(\cdot,t)$
\State Define lifted initial data:
$$y_0 = \phi_1 - \tilde{\psi}(\cdot,0), \qquad  y_1 = \phi_2 - \tilde{\psi}_t(\cdot,0)$$
\State Project onto $\mathbb{V}_h^{r,0}$: $y_h^0 = \Pi^1_0 \phi_1, \qquad \texttt{u}_h^0 = \Pi^1_0 \phi_2$
\State Compute $\texttt{m}_h^0 \in \mathbb{V}_h^{r,0}$ from discrete equilibrium at $t_0$

\vspace{0.2cm}
\For{$n=0,1,\dots,\texttt{N}_\texttt{T}-1$}

    \State $t_{n+1} = t_n + \Delta t$

    \vspace{0.1cm}
    \State \textbf{(a) Adaptive Mesh Update}
    \State Refine/coarsen mesh near crack tip to obtain $\Tc_h^{n+1}$
    \State Construct new space $\mathbb{V}_h^{r,n+1} \subset \mbox{H}_0^1([\texttt{a}, \texttt{b}])$

    \vspace{0.1cm}
    \State \textbf{(b) Mesh-to-Mesh Projection}

   $$ \widehat{y}_h^{\,n} := \Pi_n^{n+1} y_h^n, \quad \widehat{\texttt{u}}_h^{\,n} := \Pi_n^{n+1} \texttt{u}_h^n, \quad \widehat{\texttt{m}}_h^{\,n} := \Pi_n^{n+1} \texttt{m}_h^n$$

    \vspace{0.1cm}
    \State \textbf{(c) Newmark Predictors}
    $$y_h^{\,n+1}= \widehat{y}_h^{\,n} + \Delta t\, \widehat{\texttt{u}}_h^{\,n} + \Delta t^2\, \Big(\frac{1}{2}-\mu \Big) \widehat{\texttt{m}}_h^{\,n},
   \qquad 
    \texttt{u}_h^{\,n+1} = \widehat{\texttt{u}}_h^{\,n}+ \Delta t\, (1-\nu) \widehat{\texttt{m}}_h^{\,n}$$

    \vspace{0.1cm}
    \State \textbf{(d) Solve \texttt{HHT}-$\alpha$-Equilibrium for $\texttt{m}_h^{n+1}$}
    \State Find $\texttt{m}_h^{n+1} \in \mathbb{V}_h^{k,n+1}$ such that
    $$(\texttt{m}_h^{n+1},v_h) + \Ac\big(y_h^{n+1-\alpha}; v_h,t_{n+1-\alpha}\big)  =  \mathcal{F}(v_h,t_{n+1-\alpha}),  \quad \forall v_h \in \mathbb{V}_h^{r,n+1},$$
    \State (Solve nonlinear system using Newton iterations.)

    \vspace{0.1cm}
    \State \textbf{(e) Newmark Update}
    \begin{align*}
     y_h^{n+1} =\widehat{y}_h^{\,n}+ \mu \Delta t^2 \texttt{m}_h^{n+1} \qquad
     \texttt{u}_h^{n+1} =\widehat{u}_h^{\,n} + \nu \Delta t\, \texttt{m}_h^{n+1}     
    \end{align*}

    \vspace{0.1cm}
    \State \textbf{(f) Recover Physical Solution:} $u_h^{n+1}= y_h^{n+1} + \tilde{\psi}_h^{\,n+1}$
\EndFor
\end{algorithmic}
\end{algorithm}
\noindent
{\bf \textsf{Errors and Order of Convergence.}}
Let $u(t_n)$ denote the exact solution at time $t_n$ and $u_h^n$ the corresponding numerical approximation. We cmpute the  $\mbox{L}^2$ and $\mbox{H}^1$-errors, for $n \in [0:\texttt{N}_\texttt{T}]$,  by
\begin{equation}
	\texttt{E}_{\mbox{L}^2}(h_\ell) = \| u(t_n) - u_{h_\ell}^n \|_{\mbox{L}^2([\texttt{a}, \texttt{b}])}, \quad \text{and}
	\quad \texttt{E}_{\mbox{H}^1}(h_\ell) = \| u(t_n) - u_{h_\ell}^n \|_{\mbox{H}^1([\texttt{a}, \texttt{b}])}, ~ \text{respectively}.
\end{equation}
Further, the experimental order of convergence (EOC)  for $\mbox{L}^2$ and  $\mbox{H}^1$-errors are computed using the following formula 
\begin{equation}
	\text{EOC}_{\mbox{L}^2} 
	= \frac{\log \left( \texttt{E}_{\mbox{L}^2}(h_\ell) / \texttt{E}_{\mbox{L}^2}(h_{\ell+1}) \| \right)}{\log \left(h_\ell / h_{\ell+1} \right)} \quad \text{and} \quad \text{EOC}_{\mbox{H}^1} 
	= \frac{\log \left( \texttt{E}_{\mbox{H}^1}(h_\ell) / \texttt{E}_{\mbox{H}^1}(h_{\ell+1}) \| \right)}{\log \left(h_\ell / h_{\ell+1} \right)},
\end{equation}
respectively.  Here $h_\ell$ and $h_{\ell+1}$ are two successive mesh sizes.   Further, we will use the following formulae to obtain the following physical quantities.\\

\noindent 
{\bf \textsf{Discrete strain, stress, and local wave speed.}} We will figure out the discrete strain using the formula
\begin{equation*}
	\varepsilon_h(x, t_n)= \partial_x u_h^n(x), \qquad  x\in (\texttt{a}, \texttt{b}), \, n\in[0 :\texttt{N}_\texttt{T}].	\label{eq:discrete_strain}
\end{equation*}
Next, the discrete stress is obtained by using the nonlinear constitutive law, which is demonstrated by 
\begin{equation*}
	\sigma_h(x,, t_n)= \frac{\varepsilon_h(x, t_n)} {\sqrt{\,1-\beta^2\bigl(\varepsilon_h(x, t_n )\bigr)^2\,}}, \qquad x\in (\texttt{a}, \texttt{b}),  \, n\in[0:\texttt{N}_\texttt{T}],
	\label{eq:discrete_stress}
\end{equation*}
provided that $1-\beta^2\bigl(\varepsilon_h(x, t_n)\bigr)^2>0$.

Moreover,  the discrete local tangent wave speed is given by
\begin{equation*}
	c_{\mathrm{loc},h}(x, t_n)	=	\rho^{-1/2}\left(1-\beta^2\bigl(\varepsilon_h(x, t_n)\bigr)^2\right)^{-3/4}, \qquad x\in  (\texttt{a}, \texttt{b}),  \, n\in[0:\texttt{N}_\texttt{T}], \label{eq:discrete_wave_speed}
\end{equation*}
which can be easily followed using the derivative of the nonlinear constitutive relation.

Additionally, the discrete stored elastic-energy density and kinetic-energy are determined using the following formulas, for $ x\in (\texttt{a}, \texttt{b}),\, n\in[0:\texttt{N}_\texttt{T}]$, 
\begin{align*}
	\psi_h(x, t_n)~=~\frac{1-\sqrt{1-\beta^2(\varepsilon_h(x, t_n))^2}}{\beta^2}, \quad \text{and} \quad \mathpzc{K}_h(x, t_n)~=~\frac{\rho}{2}\left(v_h^n(x)\right)^2,\quad v_h^n=\partial_t u_h^n,
\end{align*}
where $v_h^n$ denotes the discrete velocity. 
 For convenience, we use the following abbreviations and parameters in our computation, which are described in tables ~\ref{computationalnotation} and \ref{parameterslist}, respectively:
\begin{table}[H]
	\centering
	\caption{Notation used in the computational results.}
	\label{computationalnotation}
	\renewcommand{\arraystretch}{1.20}
	\begin{tabular}{ll}
		\hline
		\textbf{Symbol} & \textbf{Description} \\
		\hline
		$N_e$ & Number of finite elements \\
		$N_{\mathrm{dof}}$ & Number of degrees of freedom \\
		$h$ & Mesh size \\
		$\Delta t$ & Time-step size \\
		$N_t$ & Total number of time steps \\
		$\overline{N}_{\mathrm{Newton}}$ & Average number of Newton iterations per time step \\
		$N_{\mathrm{Newton}}^{\max}$ & Maximum number of Newton iterations during the simulation \\
		$t_{\mathrm{CPU}}$ & Total CPU time \\
		$\|e\|_{L^\infty_{\mathrm{nodal}}}$ & Maximum nodal error \\
		$c_{\min}^{\texttt{T}}$ & Minimum local wave speed at the final simulation time $\texttt{T}$ \\
		$c_{\max}^{\texttt{T}}$ & Maximum local wave speed at the final simulation time $\texttt{T}$ \\
		$c_{\max}^{ST}$ & Maximum local wave speed over the entire space--time domain \\
		$D_{\min}^{ST}$ & Minimum constitutive denominator over the entire simulation \\
		\hline
	\end{tabular}
\end{table}

\begin{table}[H]
	\centering
	\renewcommand{\arraystretch}{1.25}
	\caption{Parameters used in the numerical simulations.}
	\label{parameterslist}
	\begin{tabular}{lll} 
		\hline
		\textbf{Parameter}  & \textbf{Value} \\
		\hline
		Material density $(\rho)$  & $1.0$ \\
		Strain-limiting parameter  $(\beta)$    & $0.25$ \\
		HHT-$\alpha$ parameter  $(\alpha)$ &
		$0.10$  \\
		Newmark parameter $(\mu)$  &
		$0.60$  \\
		Newmark parameter $(\nu)$  &
		$0.3025$  \\
		Final time  $(\mathtt{T})$ &
		$0.25$  \\
		Newton tolerance for convergence & $10^{-10}$  \\
		Maximum Newton iterations & $30$  \\
		Gaussian quadrature points per element &
		$2$ \\
		\hline
	\end{tabular}
\end{table}
Our objective is to validate the proposed theoretical framework through numerical experiments incorporating the above features. We will consider two different examples. The first example is the nonlinear case with zero boundary, and the second example is the nonlinear case with nonzero boundary.
\begin{exam}[Homogeneous Dirichlet   Boundary] \label{exple5.1}
Consider the model
\begin{align*}
\partial_{tt}u(x,t)-\partial_x\!\left(\frac{ \partial_x u}{[1-\beta^2\,|\partial_x u(x,t)|^2]^{\frac{1}{2}}}\right) = f(x,t),
\quad x\in(0,1), \ t\in(0,\texttt{T}),
\end{align*}
with homogeneous Dirichlet boundary conditions $u(0,t)=u(1,t)=0$ and  initial conditions  $u(x,0)=\sin(\pi x) $ and $ u_t(x,0)=0.$ Further, we  $u(x,t)$ be the  exact solution given by  $$u(x,t)=\sin(\pi x)\cos(\pi t),$$ then the load fuction $f(x,t)$ is computed accordingly. 
\end{exam}
To analyze the above problem, the numerical simulations were performed using the following parameters. The material density was taken as $\rho=1.0$, while the strain-limiting constitutive parameter was set to $\beta=0.25$, which satisfies the admissibility condition $\beta<1/\pi$ required for the manufactured solution. Time integration was carried out using the HHT-$\alpha$ method with $\alpha=0.10$, corresponding to the Newmark parameters $\mu=0.60$ and $\nu=0.3025$. 

The computations were performed over the time interval $[0,\texttt{T}]$ with the final time $\texttt{T}=0.25$. Note that $0=t_{0}<t_{1}<\cdots<t_{\texttt{N}_\texttt{T}}=\texttt{T}$ denote a uniform partition of the time interval $[0,\texttt{T}]$. The time step is defined by $\Delta t=\frac{\texttt{T}}{N},$
where $N$ denotes the total number of time steps. In the present implementation, the time step is chosen automatically according to the spatial mesh size. Specifically, if $h$ denotes the maximum element size of the finite element mesh, then the number of time steps is selected as $N=\left\lceil\frac{\texttt{T}}{h}\right\rceil,$
and the corresponding time step is computed as
$$\Delta t=\frac{\texttt{T}}{N}.$$
This choice yields $\Delta t\approx h$ except for element $N_e=2$, thereby balancing the temporal and spatial discretization errors while ensuring that the final time $\texttt{T}=0.25$ is reached exactly. 

At each time step, we solved the nonlinear system using Newton's method with a convergence tolerance $10^{-10}$ and a maximum of iterations $30$. All computations use  meshes $N_e=2^j$, $j=1,\,2,\ldots,\, 10$. 
Numerical integration was carried out using a two-point Gaussian quadrature rule on each element. 

\noindent
Table~\ref{Convergence_errors} presents the convergence history of the proposed finite element approximation on a sequence of uniformly refined meshes consisting of elements $N_e=2^j$, $j=1,\,2,\ldots,\, 10$. The corresponding errors in the $\mbox{L}^2$ and $\mbox{H}^1$-norms, together with their experimentally observed orders of convergence (EOC) are reported. Both errors are  decreasing monotonically on uniform refined meshes, which demonstrate the stability and the consistency of the proposed HHT-$\alpha$ finite element scheme. In particular, the $\mbox{L}^2$error decreases from $9.1718\times10^{-2}$  to  $4.2050\times10^{-8}$, while the $\mbox{H}^1$-error decreases from  $6.9456\times10^{-1}$ to $1.3912\times10^{-3}$, confirming continuous improvement of the numerical approximation with mesh refinement.

\begin{figure}[H]
	\centering
	\begin{subfigure}[b]{0.48\textwidth}
		\centering
		\includegraphics[width=\linewidth,height=0.20\textheight]{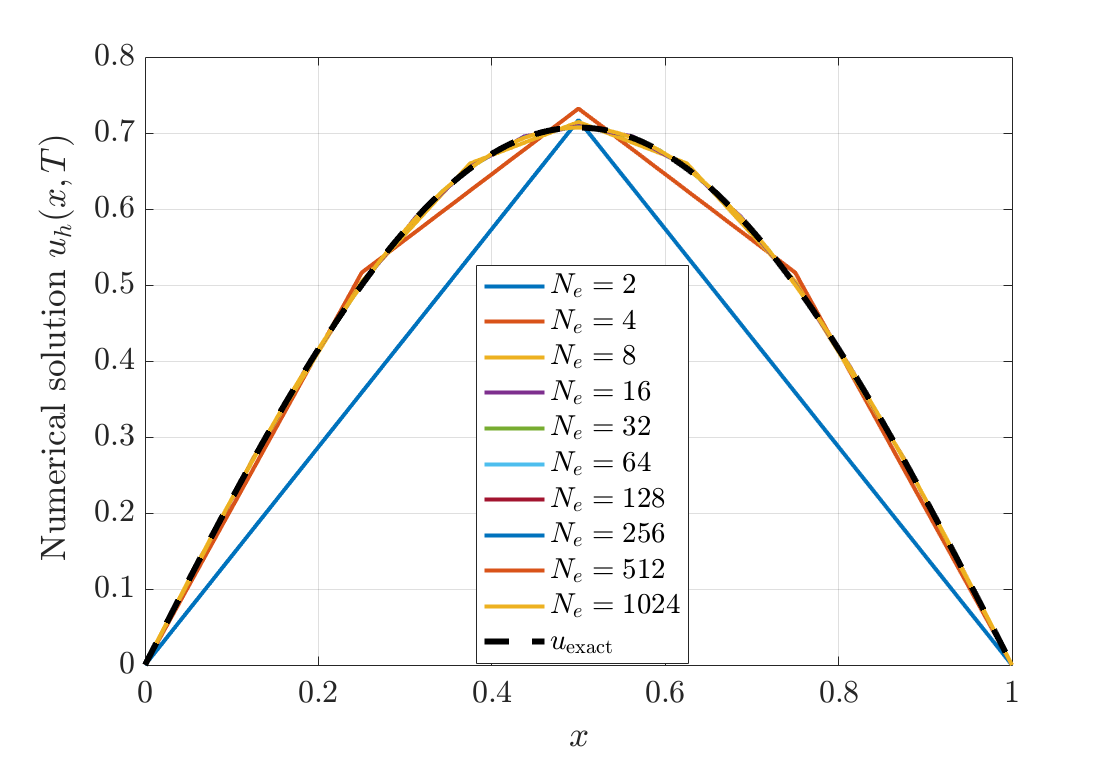}
		\caption{Approximate displacement at Final time}
		\label{solution_profilesa}
	\end{subfigure}
	\hfill
	\begin{subfigure}[b]{0.48\textwidth}
		\centering
		\includegraphics[width=\linewidth,height=0.20\textheight]{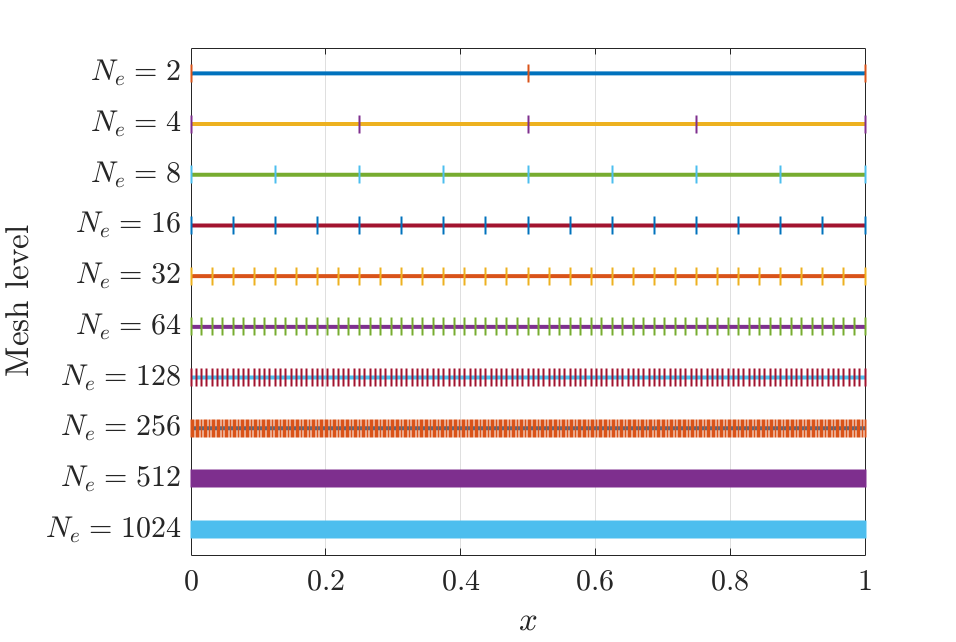}
		\caption{Uniform  finite-element mesh levels}
		\label{meshes}
	\end{subfigure}
	
	\vspace{0.3cm}
	\begin{subfigure}[b]{0.48\textwidth}
		\centering
		\includegraphics[width=\linewidth,height=0.20\textheight]{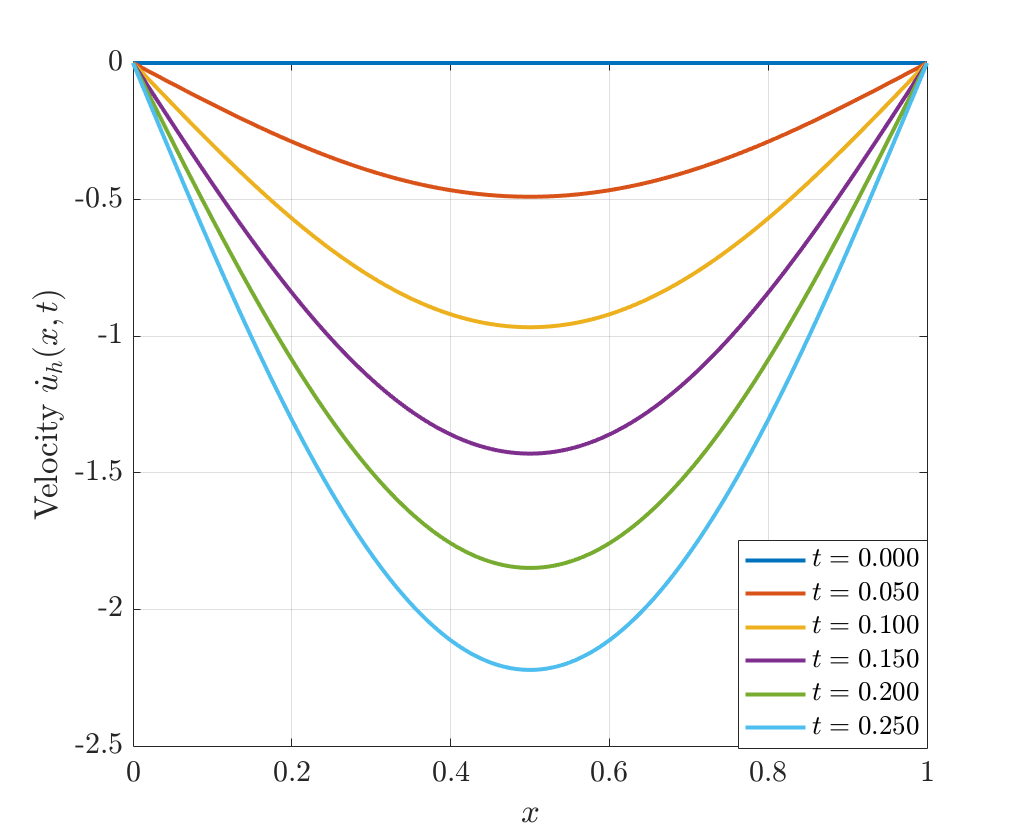}
		\caption{Approximate velocity}
		\label{velocityprofilesc}
	\end{subfigure}
	\hfill
	\begin{subfigure}[b]{0.48\textwidth}
		\centering
		\includegraphics[width=\linewidth,height=0.20\textheight]{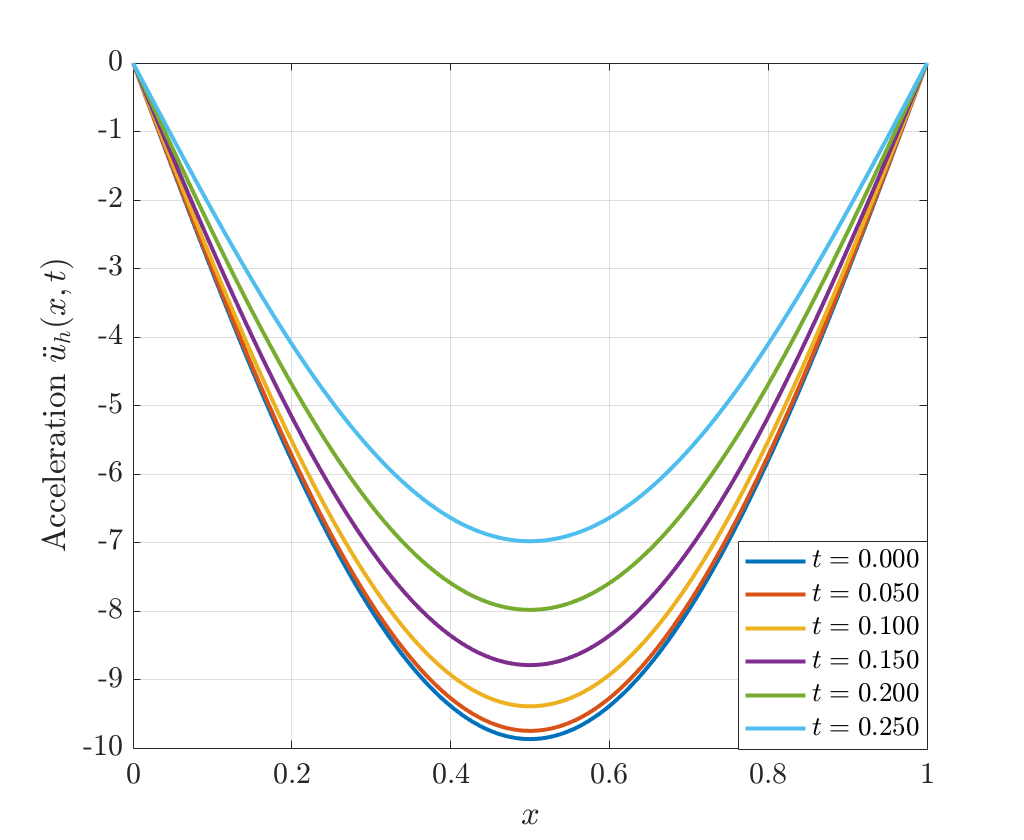}
		\caption{Approximate acceleration}
		\label{accelerationsprofilesd}
	\end{subfigure}
	\caption{Profiles of approximate solution,  displacement, velocity and acceleration}
	\label{profilesdva}
\end{figure}

It is evident from Table~\ref{Convergence_errors} that the experimental order of convergence (EOC) rates closely match the theoretical predictions for linear finite elements. The $\mbox{L}^2$-error exhibits an EOC approaching $2$, while the $\mbox{H}^1$-error consistently achieves an EOC of approximately $1$ throughout the refinement process. After the first refinement level, the convergence rates quickly become stable, with the $\mbox{L}^2$-EOC reaching $2.0124$ and the $\mbox{H}^1$-EOC remaining practically at $1.0000$ on the finer meshes. This behavior shows that the computations have reached the asymptotic convergence regime where the discretization error is dominated by the mesh size and not the preasymptotic effects.
\begin{table}[H]
	\centering
	\caption{Convergence history of the finite element approximation.}
	\label{Convergence_errors}
	\renewcommand{\arraystretch}{1.30}
	\setlength{\tabcolsep}{7pt}
	\begin{tabular}{cccccc}
		\hline
		Mesh  ($N_e$) ~&~
		$\|u-u_h\|_{\mbox{L}^2}$ ~&~
		EOC &
		~&~	$\|u-u_h\|_{\mbox{H}^1}$ ~&~
		EOC \\
		\hline
		$2$    & $9.1718\times10^{-2}$ & --      &    
		& $6.9456\times10^{-1}$ & --      \\
		$4$    & $9.7436\times10^{-3}$ & 3.2347  & 
		& $3.5710\times10^{-1}$ & 0.9597  \\
		$8$    & $1.6292\times10^{-3}$ & 2.5803  & 
		& $1.7840\times10^{-1}$ & 1.0012  \\
		$16$   & $2.8454\times10^{-4}$ & 2.5174  & 
		& $8.9093\times10^{-2}$ & 1.0017  \\
		$32$   & $5.5909\times10^{-5}$ & 2.3475  & 
		& $4.4526\times10^{-2}$ & 1.0007  \\
		$64$   & $1.2241\times10^{-5}$ & 2.1914  & 
		& $2.2260\times10^{-2}$ & 1.0002  \\
		$128$  & $2.8583\times10^{-6}$ & 2.0985  & 
		& $1.1129\times10^{-2}$ & 1.0001  \\
		$256$  & $6.9041\times10^{-7}$ & 2.0496  & 
		& $5.5647\times10^{-3}$ & 1.0000  \\
		$512$  & $1.6966\times10^{-7}$ & 2.0248  & 
		& $2.7823\times10^{-3}$ & 1.0000  \\
		$1024$ & $4.2050\times10^{-8}$ & 2.0124  & 
		& $1.3912\times10^{-3}$ & 1.0000  \\
		\hline
	\end{tabular}
\end{table}
Tables~\ref{Mesh_solver} and~\ref{Error_wave} summarize the computational performance and numerical behavior of the proposed finite element scheme over a sequence of uniformly refined meshes. As shown in Table~\ref{Mesh_solver}, the mesh size and time-step size are reduced proportionally, while the number of time steps increases accordingly to maintain a consistent spatial--temporal discretization. The nonlinear Newton solver exhibits robust convergence throughout all simulations, requiring at most four iterations per time step on the coarsest meshes and only two iterations on the finer meshes. As expected, the total CPU time increases with mesh refinement due to the larger number of degrees of freedom and time steps. The maximum nodal error reduces rapidly with mesh refinement as shown in Table~\ref{Error_wave}, confirming the convergence of the proposed method. Furthermore, the computed local wave-speed quantities stabilize as the mesh is refined: the minimum wave speed at the final time converges to approximately $1.0000$, while the maximum wave speeds at the final time and over the entire space--time domain approach $1.3186$ and $2.0534$, respectively. At the same time, the minimum constitutive denominator remains strictly positive throughout the simulations, indicating that the admissibility condition is satisfied and the nonlinear constitutive relation remains well-defined for all mesh levels.
\begin{table}[H]
	\centering
	\caption{Mesh characteristics, nonlinear solver performance, and computational cost.}
	\label{Mesh_solver}
	\renewcommand{\arraystretch}{1.15}
	\setlength{\tabcolsep}{5pt}
	\resizebox{\textwidth}{!}{
		\begin{tabular}{ccccccccc}
			\hline
			$N_e$ &
			$N_{\mathrm{dof}}$ &
			$h$ &
			$\Delta t$ &
			$N_t$ &
			$\overline{N}_{\mathrm{Newton}}$ &
			$N_{\mathrm{Newton}}^{\max}$ &
			$t_{\mathrm{CPU}}$ (s) \\
			\hline
			2    & 3    & $5.0000\times10^{-1}$ & $2.5000\times10^{-1}$ &   1 & 4.00 & 4 & 0.009 \\
			4    & 5    & $2.5000\times10^{-1}$ & $2.5000\times10^{-1}$ &   1 & 4.00 & 4 & 0.010 \\
			8    & 9    & $1.2500\times10^{-1}$ & $1.2500\times10^{-1}$ &   2 & 4.00 & 4 & 0.015 \\
			16   & 17   & $6.2500\times10^{-2}$ & $6.2500\times10^{-2}$ &   4 & 3.25 & 4 & 0.026 \\
			32   & 33   & $3.1250\times10^{-2}$ & $3.1250\times10^{-2}$ &   8 & 3.00 & 3 & 0.052 \\
			64   & 65   & $1.5625\times10^{-2}$ & $1.5625\times10^{-2}$ &  16 & 3.00 & 3 & 0.089 \\
			128  &129   & $7.8125\times10^{-3}$ & $7.8125\times10^{-3}$ &  32 & 2.00 & 2 & 0.156 \\
			256  &257   & $3.9063\times10^{-3}$ & $3.9063\times10^{-3}$ &  64 & 2.00 & 2 & 0.479 \\
			512  &513   & $1.9531\times10^{-3}$ & $1.9531\times10^{-3}$ & 128 & 2.00 & 2 & 1.854 \\
			1024 &1025  & $9.7656\times10^{-4}$ & $9.7656\times10^{-4}$ & 256 & 2.00 & 2 & 7.345 \\
			\hline
	\end{tabular}}
\end{table}
\begin{table}[htbp]
	\centering
	\caption{Nodal error and local wave-speed statistics.}
	\label{Error_wave}
	\renewcommand{\arraystretch}{1.15}
	\setlength{\tabcolsep}{6pt}
	\begin{tabular}{ccccccc}
		\hline
		$N_e$ ~&~
		$\|e\|_{L^\infty_{\mathrm{nodal}}}$ ~&~
		$c_{\min}^{T}$ ~&~
		$c_{\max}^{T}$ ~&~
		$c_{\max}^{ST}$ ~&~
		$D_{\min}^{ST}$\\
		\hline
		2    & $9.9825\times10^{-3}$ & 1.1087 & 1.1087 & 1.2408 & 0.7500 \\
		4    & $2.5061\times10^{-2}$ & 1.0364 & 1.2617 & 1.6818 & 0.5000 \\
		8    & $7.4392\times10^{-3}$ & 1.0090 & 1.3031 & 1.9368 & 0.4142 \\
		16   & $2.0386\times10^{-3}$ & 1.0022 & 1.3152 & 2.0223 & 0.3910 \\
		32   & $5.3190\times10^{-4}$ & 1.0006 & 1.3178 & 2.0455 & 0.3851 \\
		64   & $1.3557\times10^{-4}$ & 1.0001 & 1.3184 & 2.0514 & 0.3836 \\
		128  & $3.4191\times10^{-5}$ & 1.0000 & 1.3186 & 2.0529 & 0.3833 \\
		256  & $8.5831\times10^{-6}$ & 1.0000 & 1.3186 & 2.0533 & 0.3832 \\
		512  & $2.1501\times10^{-6}$ & 1.0000 & 1.3186 & 2.0534 & 0.3832 \\
		1024 & $5.3804\times10^{-7}$ & 1.0000 & 1.3186 & 2.0534 & 0.3832 \\
		\hline
	\end{tabular}
\end{table}
Further, Figure~\ref{profilesdva} illustrates the finite element solution and its temporal derivatives computed using the proposed HHT-$\alpha$ finite element method for the manufactured solution. Figure~\ref{profilesdva} (a) compares the finite element solutions at the final simulation time obtained on a sequence of uniformly refined meshes with the exact manufactured solution. The numerical solutions corresponding to the number of elements $N_e=2^j$, $j=1,\,2,\ldots,\, 10$ (cf., Figure \ref{profilesdva} (b))  are virtually indistinguishable from the analytical solution over the entire computational domain. 

Figure~\ref{profilesdva} (a) depicts the displacement profiles at several time instants ranging from $t=0$ to $t=0.25$. The displacement preserves the sinusoidal spatial distribution throughout the simulation while its amplitude decreases monotonically with time according to the factor $\cos(\pi t)$. Moreover, the homogeneous Dirichlet boundary conditions are satisfied exactly, with the displacement remaining zero at both ends of the domain for all times.
The corresponding velocity profiles are presented in Figure~\ref{profilesdva} (c). The analytical velocity is $u_t(x,t) = -\pi \sin(\pi x) \sin(\pi t)$ because
the velocity is initially zero and then becomes more negative with time. The maximum magnitude occurs at the center of the domain, whereas the velocity vanishes at the boundaries, in complete agreement with the analytical solution. The computed velocity remains smooth throughout the simulation, indicating the temporal accuracy of the HHT-$\alpha$ scheme.

Figure~\ref{profilesdva} (d) shows the acceleration profiles corresponding to $ u_{tt}(x,t)=-\pi^2\sin(\pi x)\cos(\pi t).$ The acceleration has the same shape in space as the displacement, but the opposite sign. Its magnitude decays continuously over time from 1 to $\cos(\pi/4)$ as $\cos(\pi t)$ decays. The computed acceleration profiles are smooth, symmetric relative to the center of the domain, and free of spurious oscillations, showing the robustness of the proposed time integration method.  Overall, the numerical results shown in Figure~\ref{profilesdva} confirm the capability of the presented HHT-$\alpha$ finite element formulation to accurately retrieve the exact displacement, velocity, and acceleration fields over the simulation time. The numerical solution shows the expected symmetry, satisfies the prescribed boundary conditions, and shows a smooth temporal evolution. 

\begin{figure}[htpb]
	\centering
\includegraphics[width=\linewidth,height=0.45\textheight]{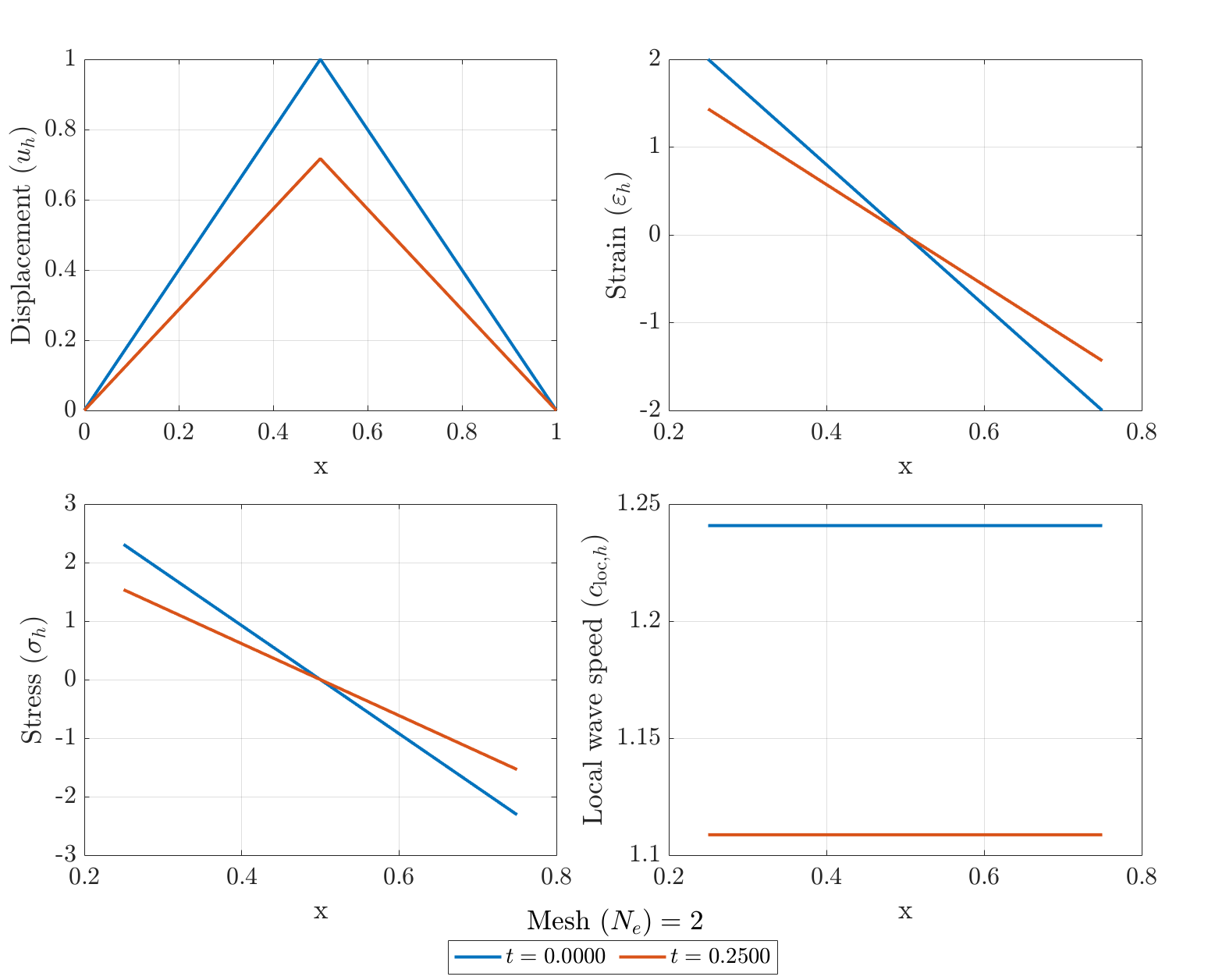} \\
\vspace{0.2cm}
\includegraphics[width=\linewidth,height=0.45\textheight]{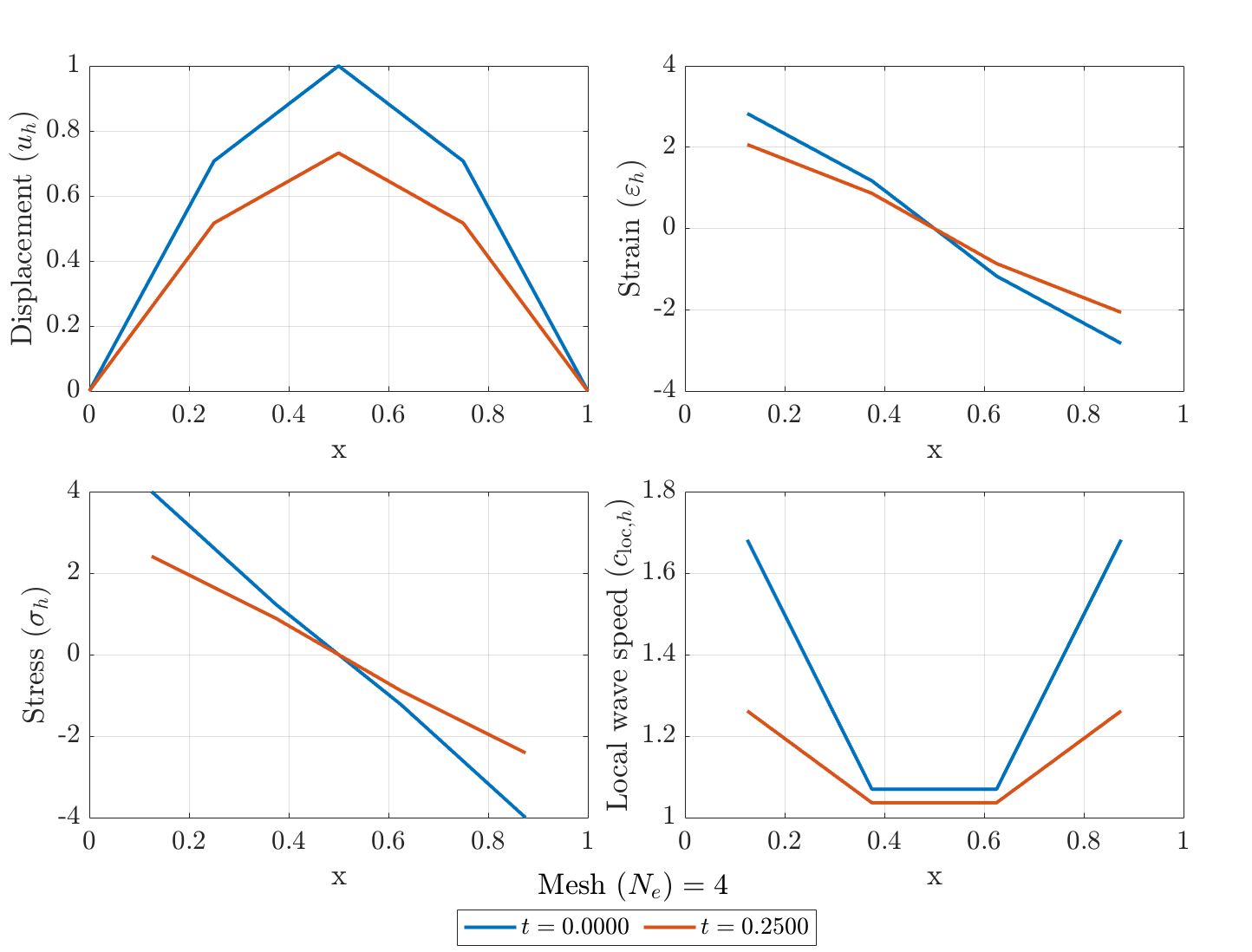}
\caption{Evolution of the discrete displacement, strain, stress, and local tangent wave speed at selected time instants for meshes with $N_e=2$ and $N_e=4$.}
\label{mechdatamesh2n4}
\end{figure}
\begin{figure}[htpb]
	\centering
	\includegraphics[width=\linewidth,height=0.45\textheight]{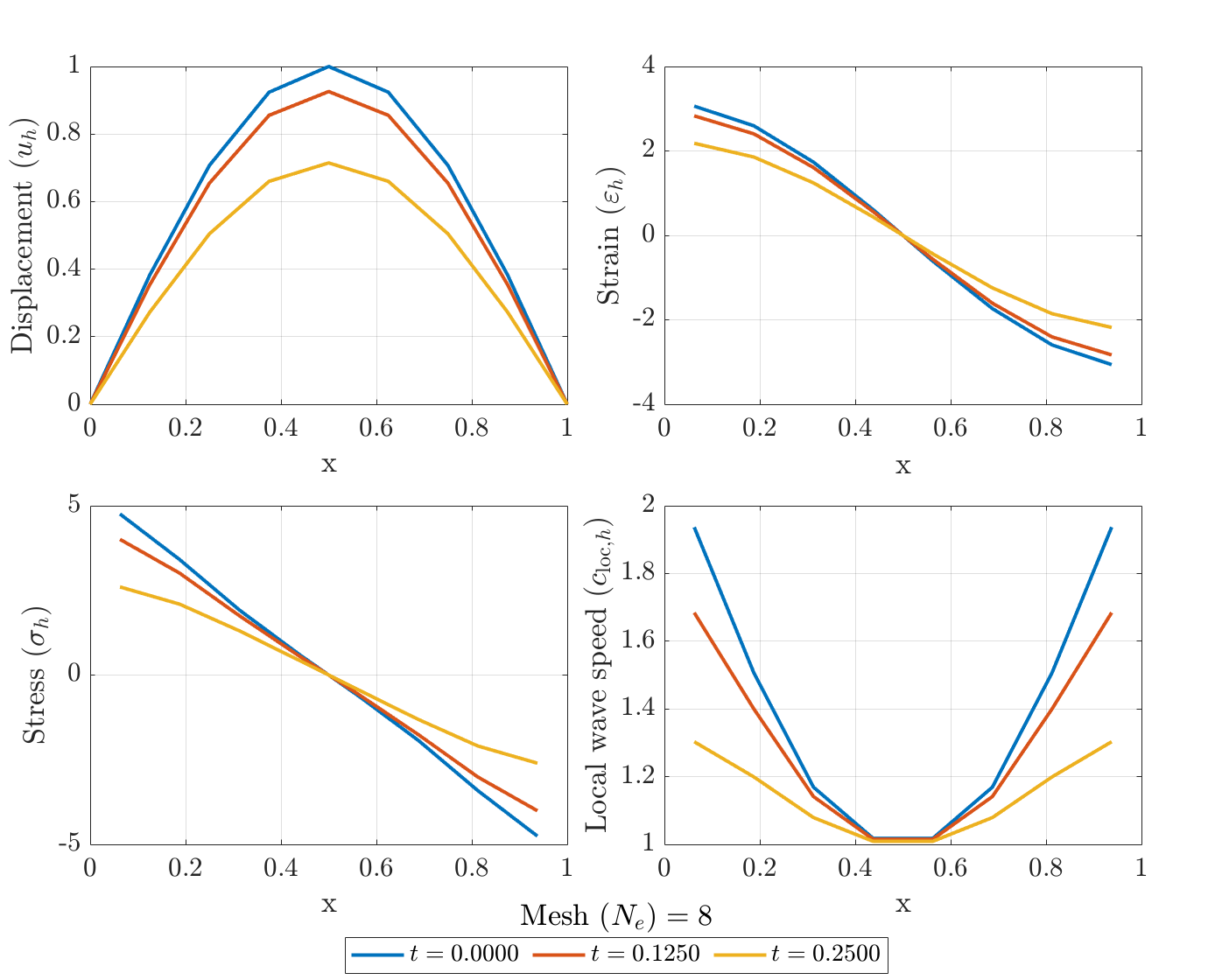} \\
	\vspace{0.3cm}
	\includegraphics[width=\linewidth,height=0.45\textheight]{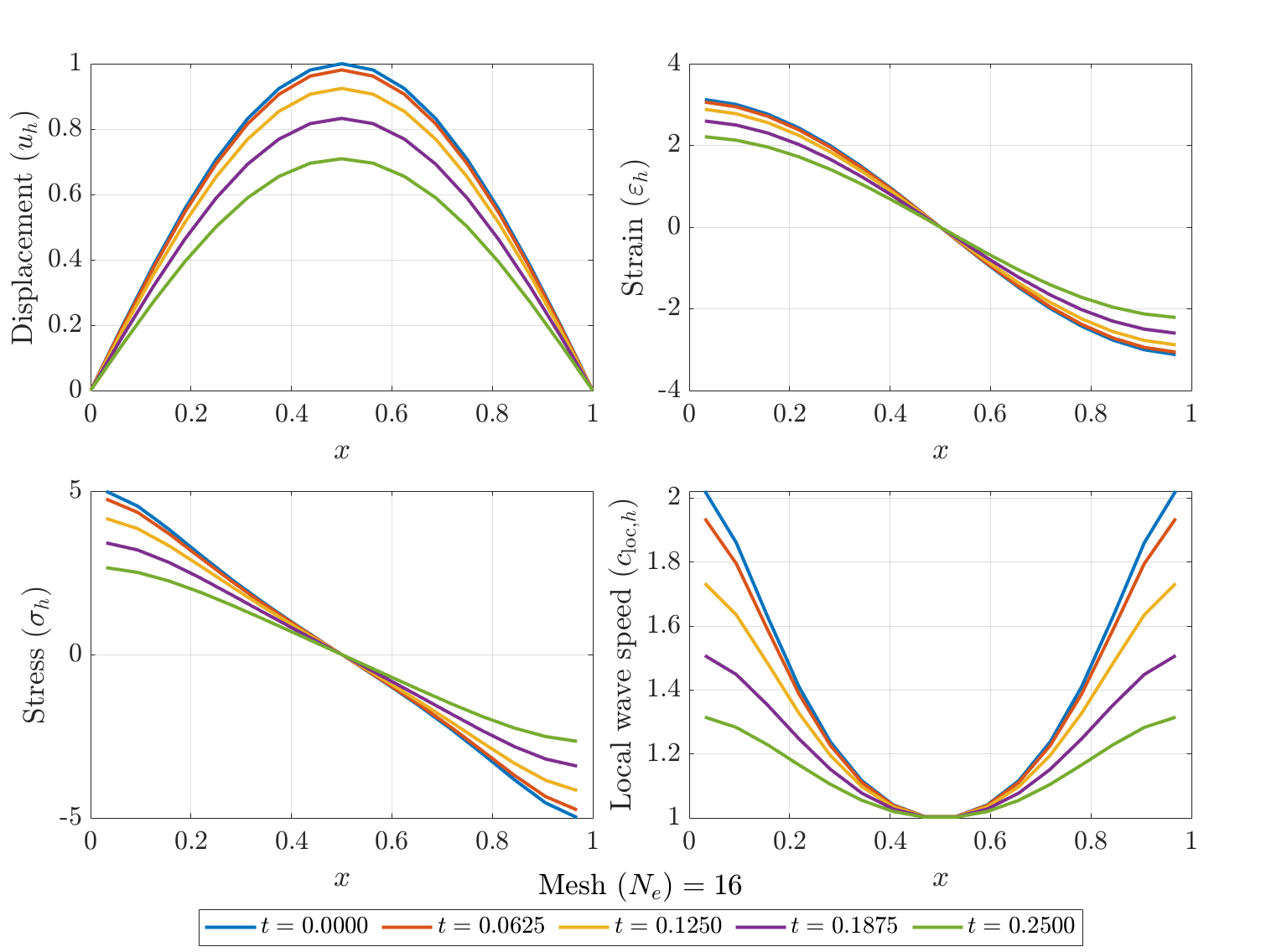}
\caption{Evolution of the discrete displacement, strain, stress, and local tangent wave speed at selected time instants for meshes with $N_e=8$ and $N_e=16$.}
	\label{mechdatamesh8n16}
\end{figure}
\begin{figure}[htpb]
	\centering
	\includegraphics[width=\linewidth,height=0.45\textheight]{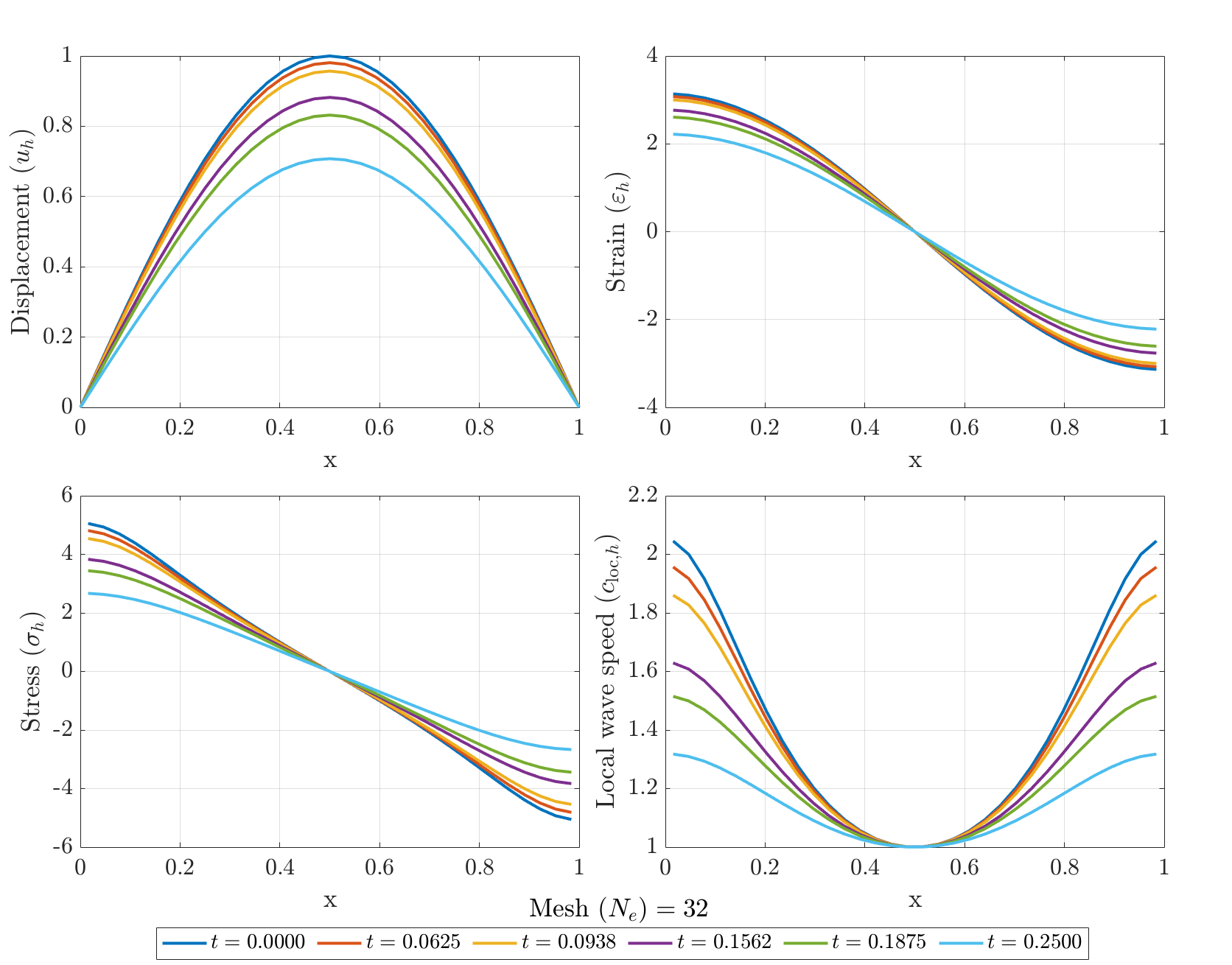} \\
	\vspace{0.3cm}
	\includegraphics[width=\linewidth,height=0.45\textheight]{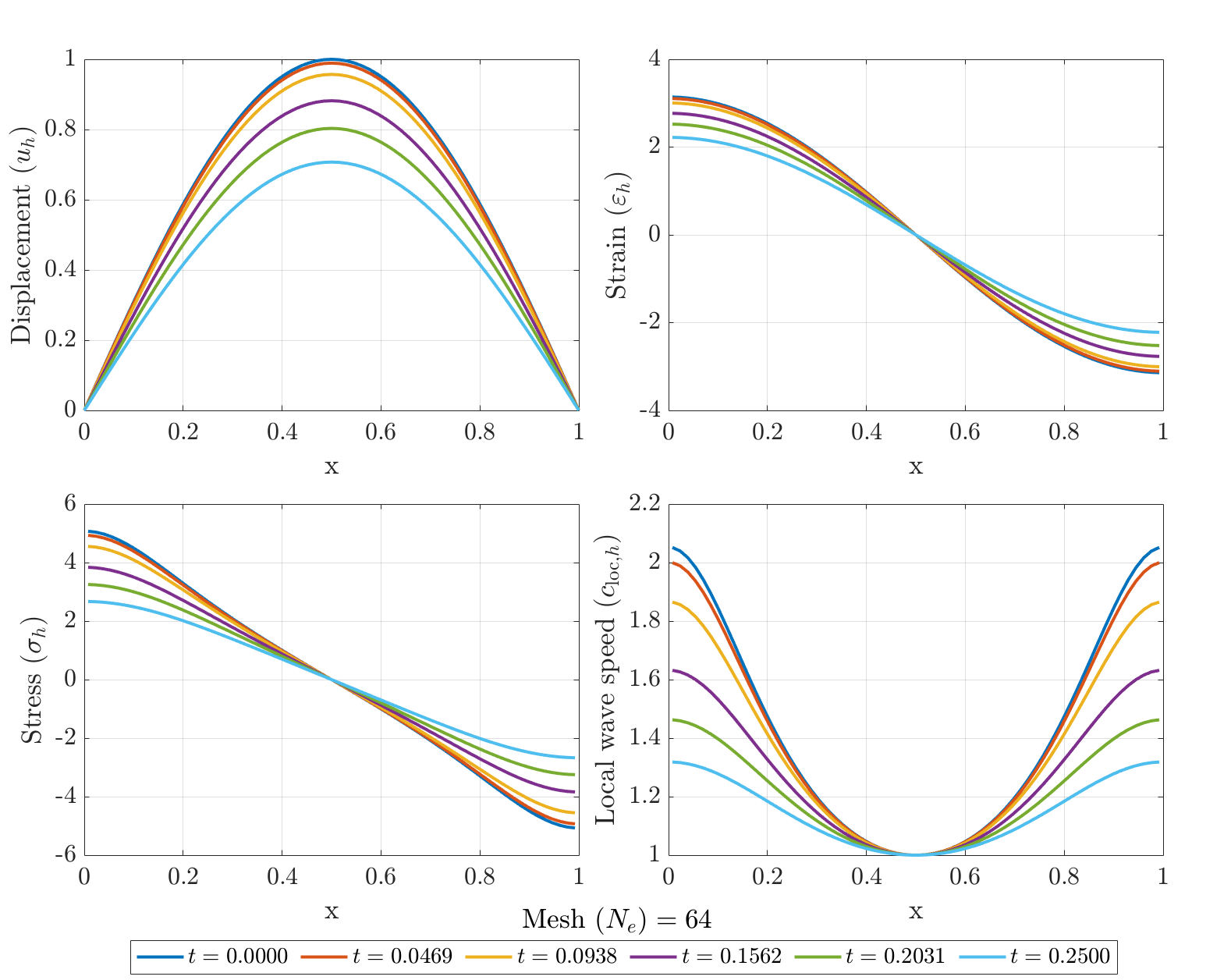}
\caption{Evolution of the discrete displacement, strain, stress, and local tangent wave speed at selected time instants for meshes with $N_e=32$ and $N_e=64$.}
	\label{mechdatamesh32n64}
\end{figure}
\begin{figure}[htpb]
	\centering
	\includegraphics[width=\linewidth,height=0.45\textheight]{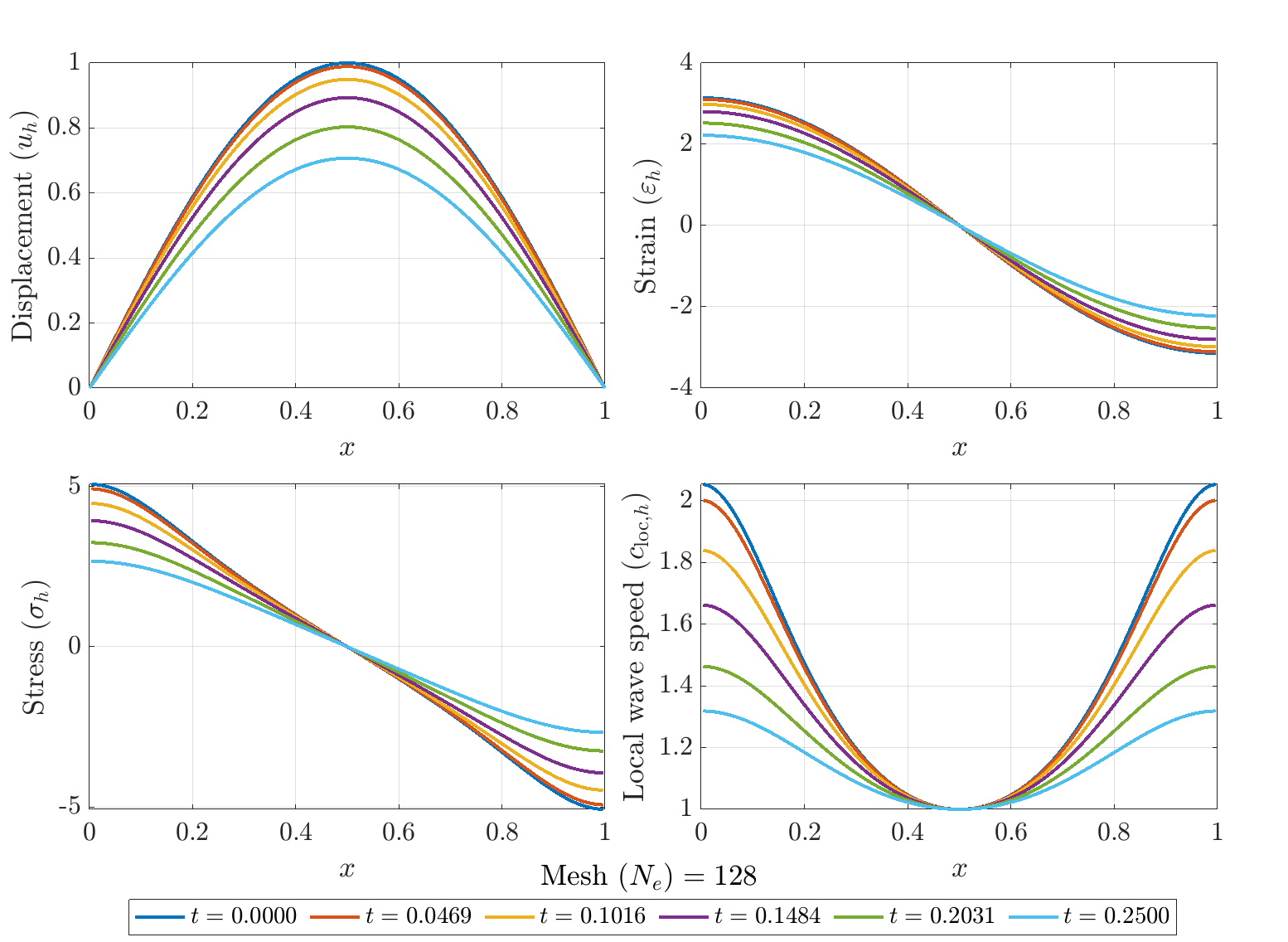} \\
	\vspace{0.3cm}
	\includegraphics[width=\linewidth,height=0.45\textheight]{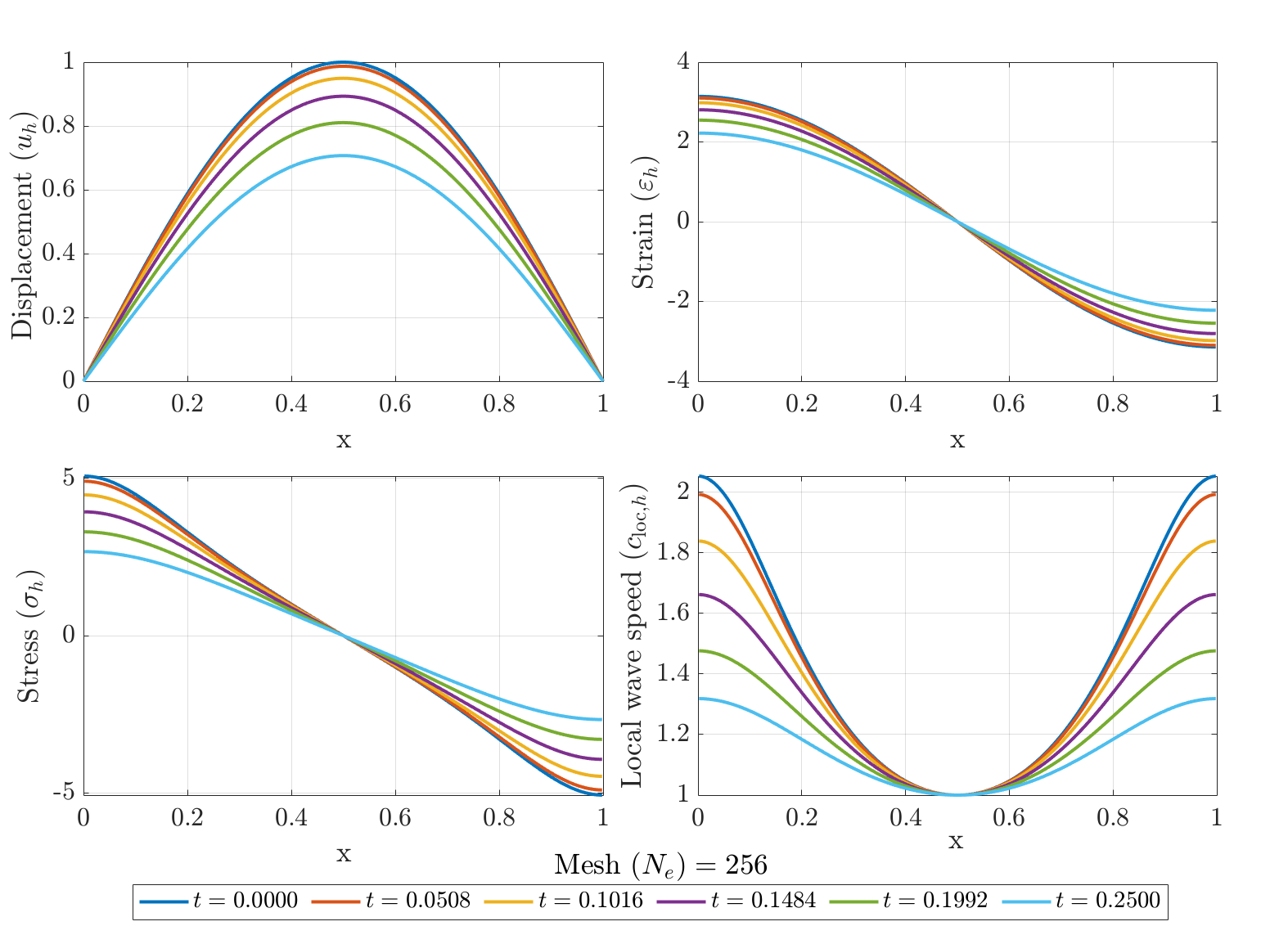}
	\caption{Evolution of the discrete displacement, strain, stress, and local tangent wave speed at selected time instants for meshes with $N_e=128$ and $N_e=256$.}
	\label{mechdatamesh128n256}
\end{figure}
\begin{figure}[htpb]
	\centering
	\includegraphics[width=\linewidth,height=0.45\textheight]{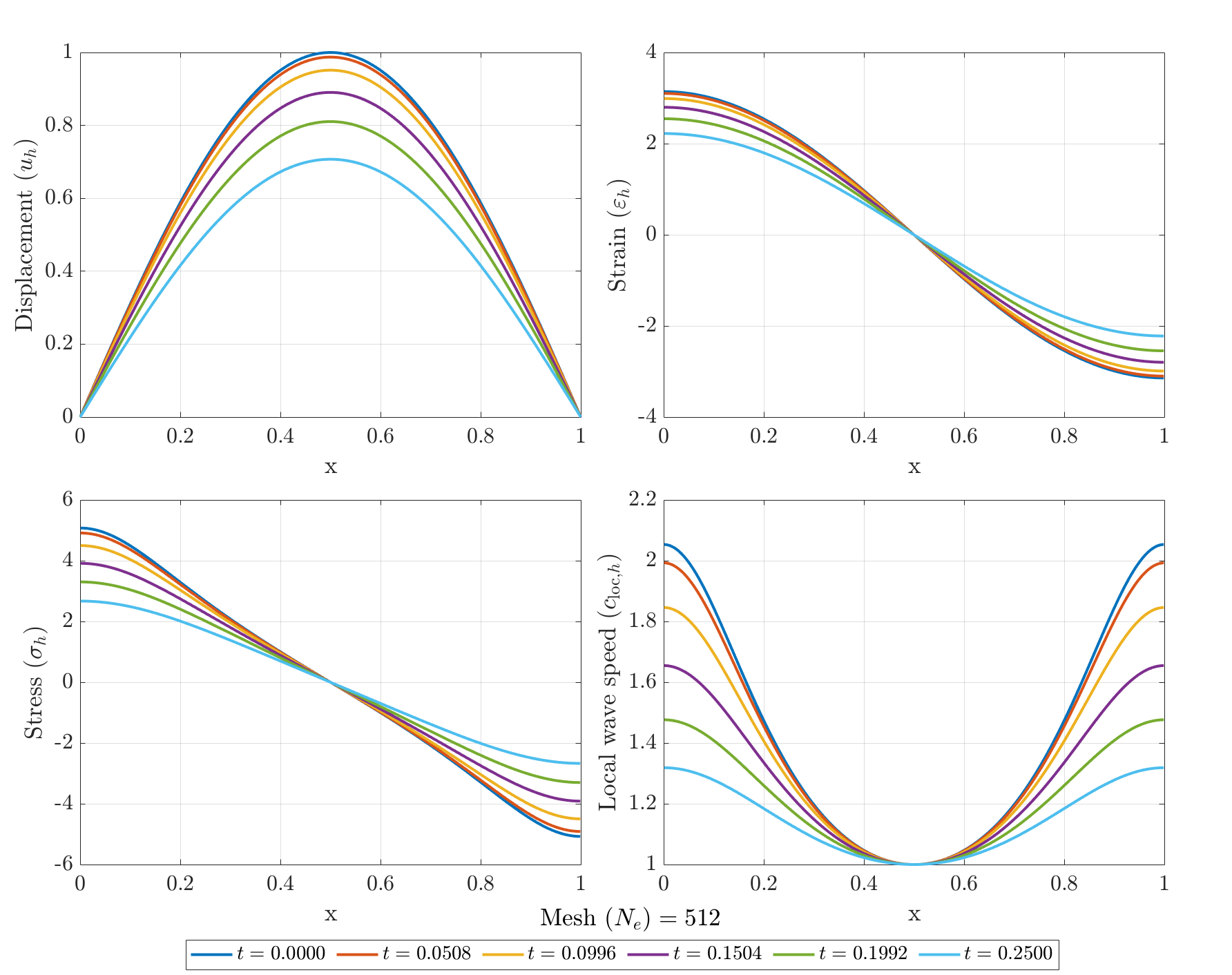} \\
	\vspace{0.3cm}
	\includegraphics[width=\linewidth,height=0.45\textheight]{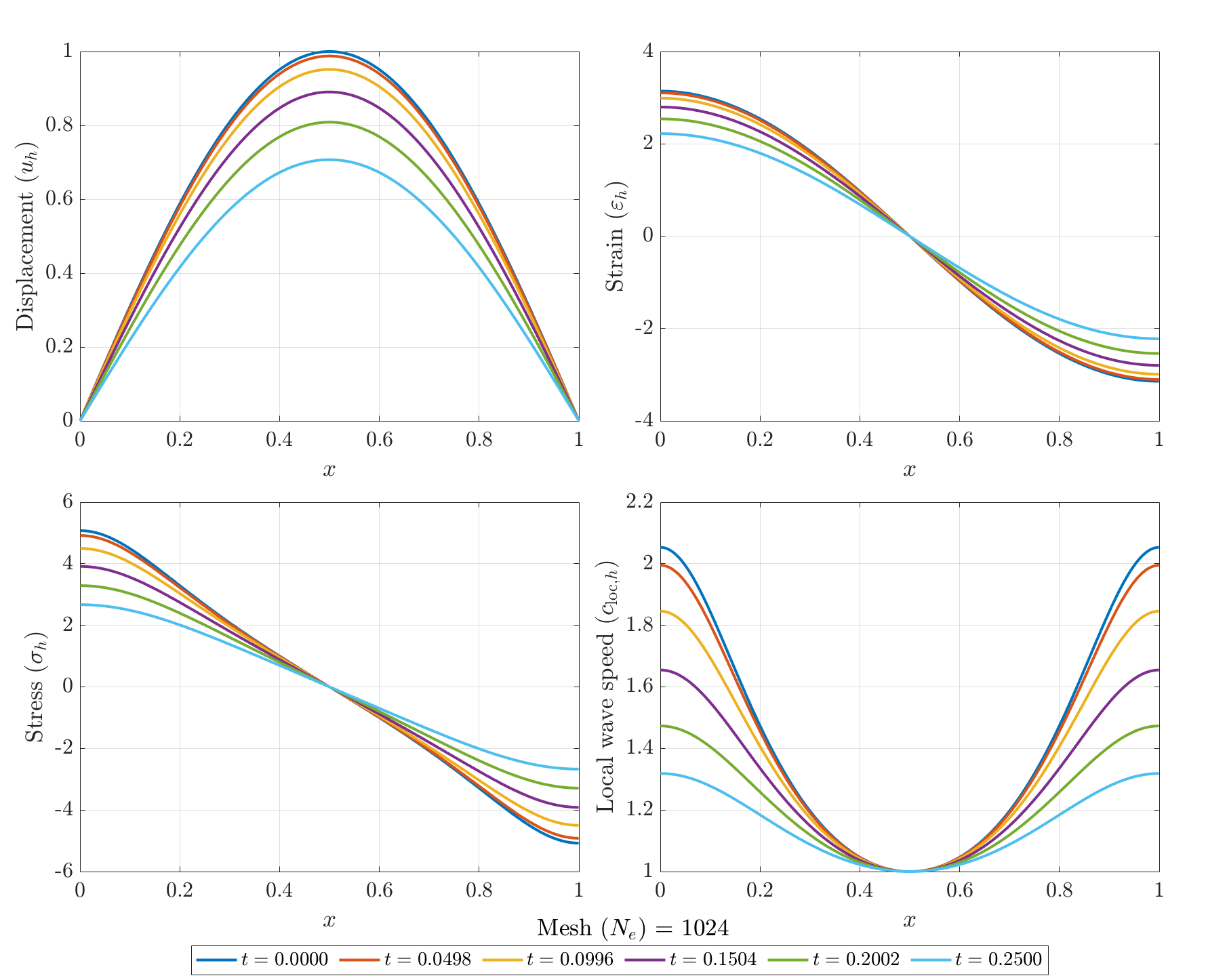}
\caption{Evolution of the discrete displacement, strain, stress, and local tangent wave speed at selected time instants for meshes with $N_e=512$ and $N_e=1024$.}
	\label{mechdatamesh512n1024}
\end{figure}
Figures~\ref{mechdatamesh2n4}--\ref{mechdatamesh512n1024} illustrate the evolution of the discrete displacement, strain, stress, and local tangent wave speed at several representative time instants for a sequence of uniformly refined meshes. On the coarsest meshes ($N_e=2$ and $4$), the displacement field is represented by piecewise linear profiles, resulting in discontinuous strain and stress distributions within individual elements and nearly constant local wave speeds. As the mesh is refined, the discrete solution converges toward a smooth representation of the underlying continuous solution. Hence, the strain and stress fields are progressively well captured, with smooth spatial variations, while preserving the expected symmetry of the problem.  The local tangent wave speed, computed from the nonlinear constitutive law, remains strictly positive throughout the simulation and exhibits its minimum near the center of the domain, where the strain magnitude is smallest, while attaining larger values near the boundaries where the strain magnitude is greatest. 
\begin{figure}[H]
	\centering
	\begin{subfigure}{0.48\textwidth}
    \includegraphics[width=1.05\linewidth,height=5cm]{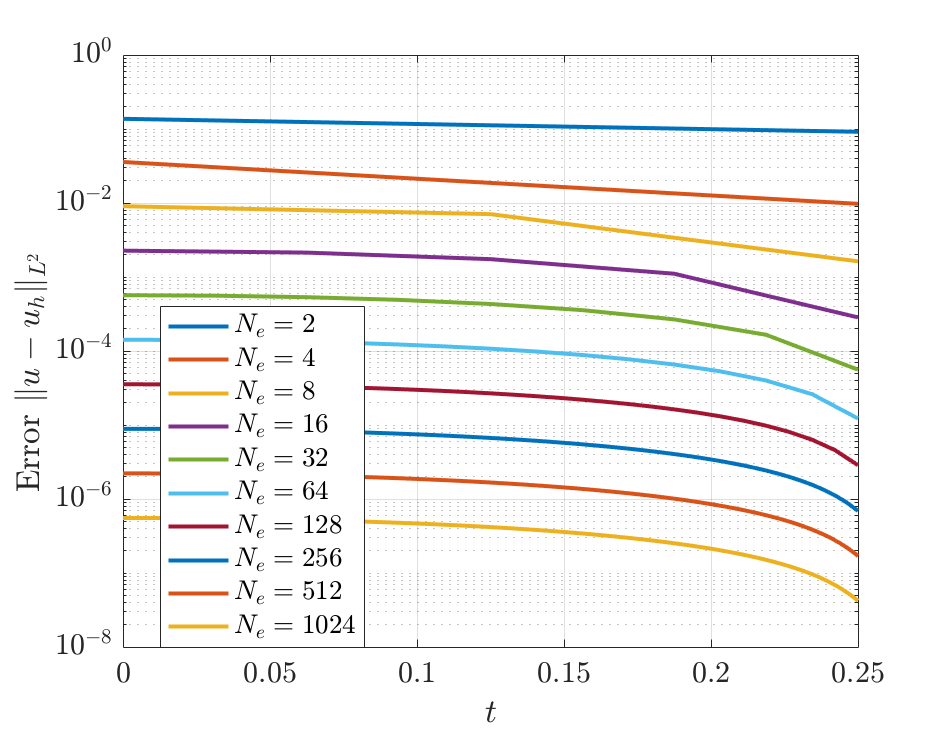}
		\caption{$\mbox{L}^2$-error vs Time ($t$)}
	\end{subfigure}
	\hfill
	\begin{subfigure}{0.48\textwidth}
	    \includegraphics[width=1.05\linewidth,height=5cm]{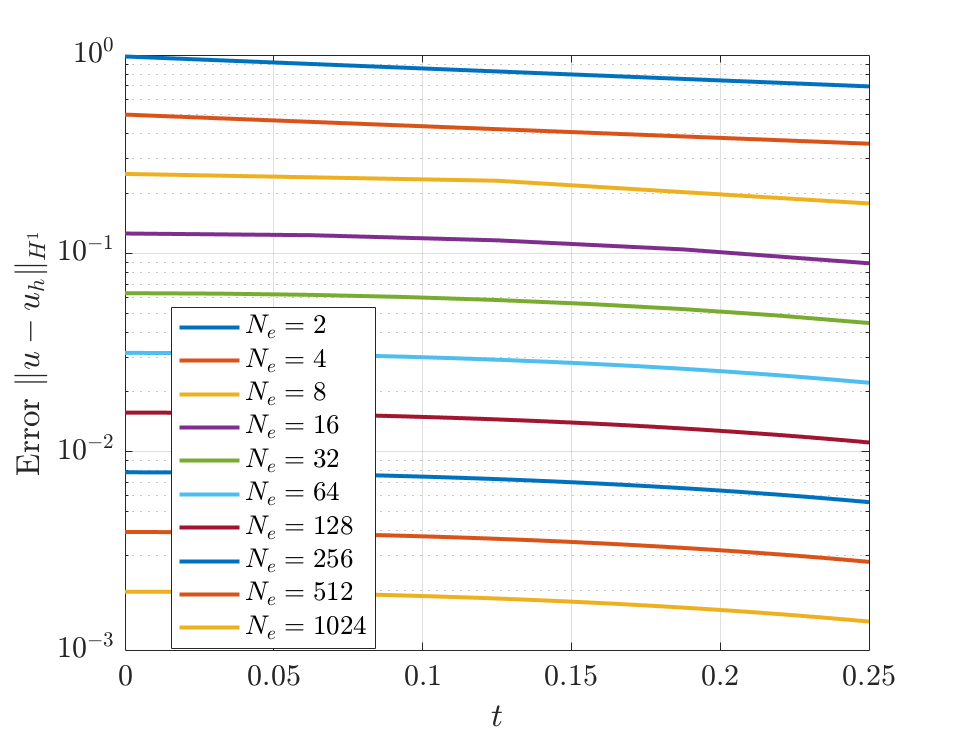}
		\caption{$\mbox{L}^2$-error vs. Time ($t$)}
	\end{subfigure}
	\vspace{0.2cm}
	\begin{subfigure}{0.48\textwidth}
	\includegraphics[width=1.05\linewidth,height=5cm]{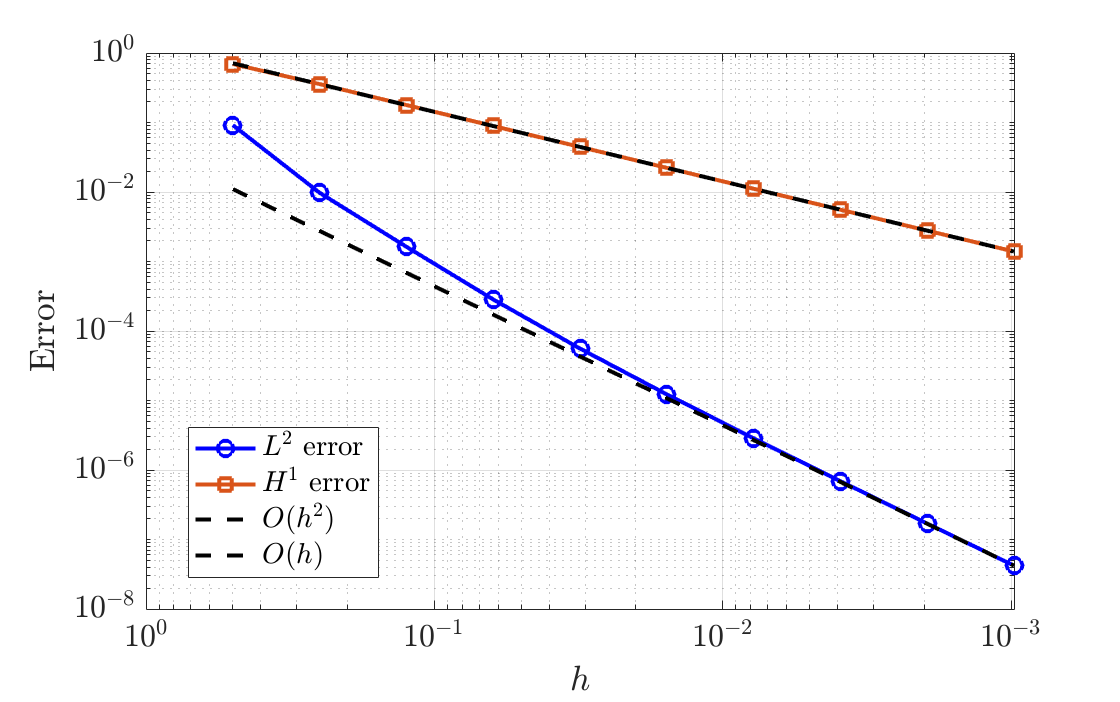}
\caption{$\mbox{L}^2$ and $\mbox{H}^1$-error vs mesh size ($h$)}
	\end{subfigure}
	\hfill
	\begin{subfigure}{0.48\textwidth}
		\includegraphics[width=1.05\linewidth,height=5cm]{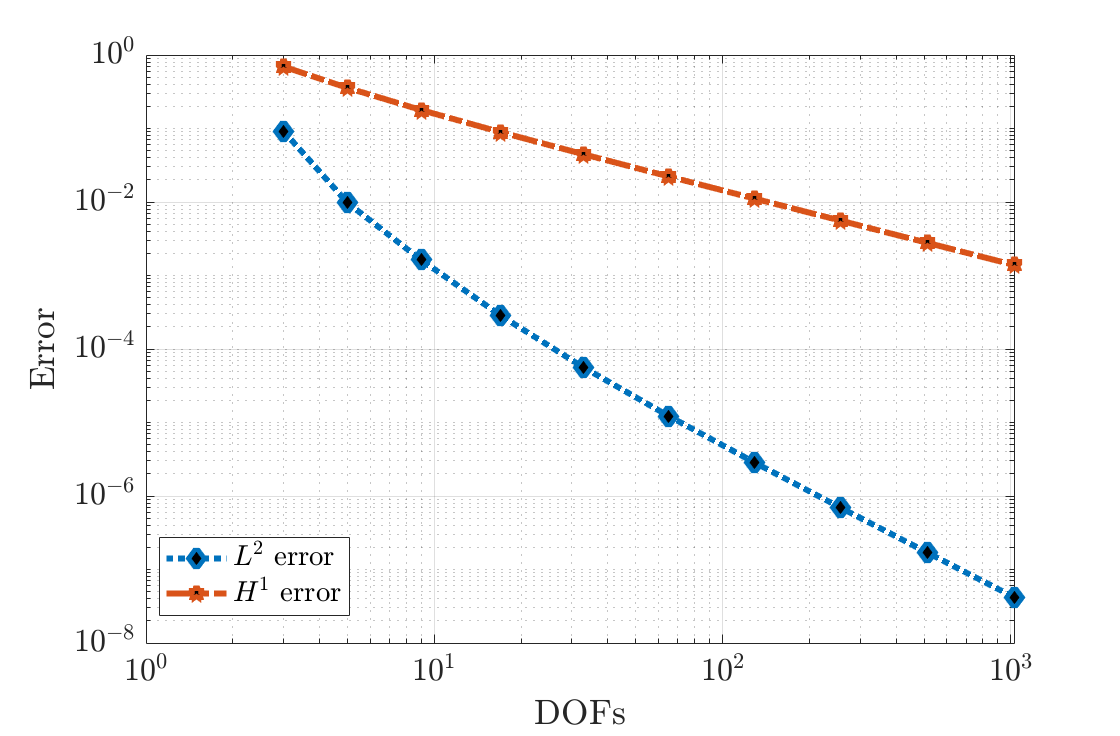}
		\caption{$\mbox{L}^2$ and $\mbox{H}^1$-error vs DOFs}
	\end{subfigure}
	\vspace{0.2cm}
	\begin{subfigure}{0.48\textwidth}
   \includegraphics[width=0.97\linewidth,height=5cm]{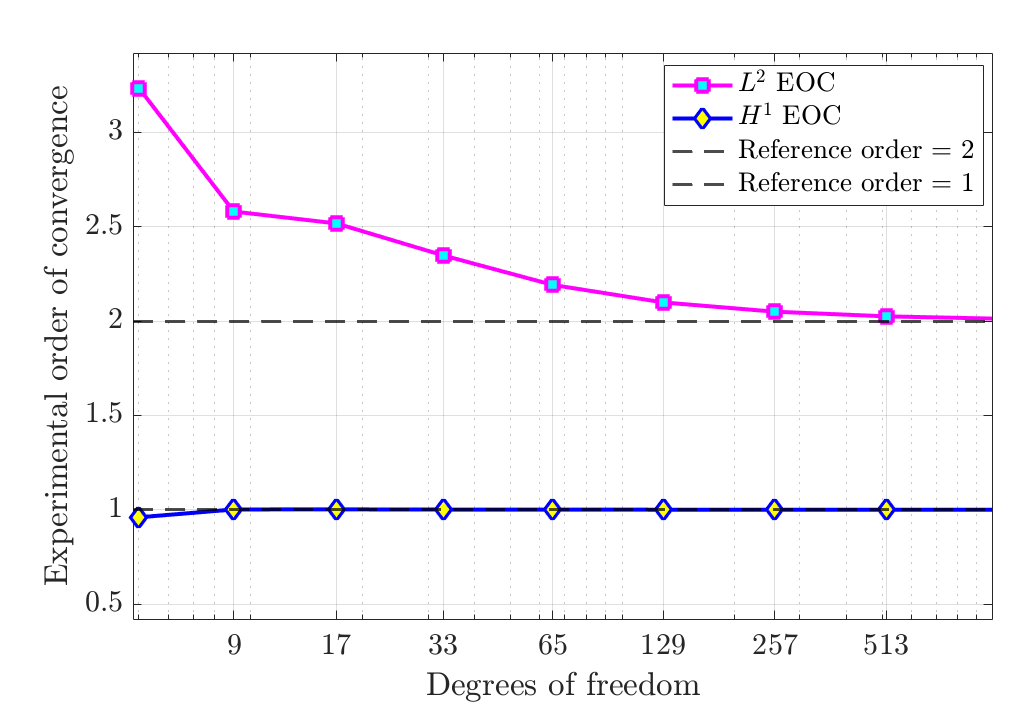}
		\caption{Experimental order of convergence}
	\end{subfigure}
	\hfill
	\begin{subfigure}{0.48\textwidth}
   \includegraphics[width=1.05\linewidth,height=5cm]{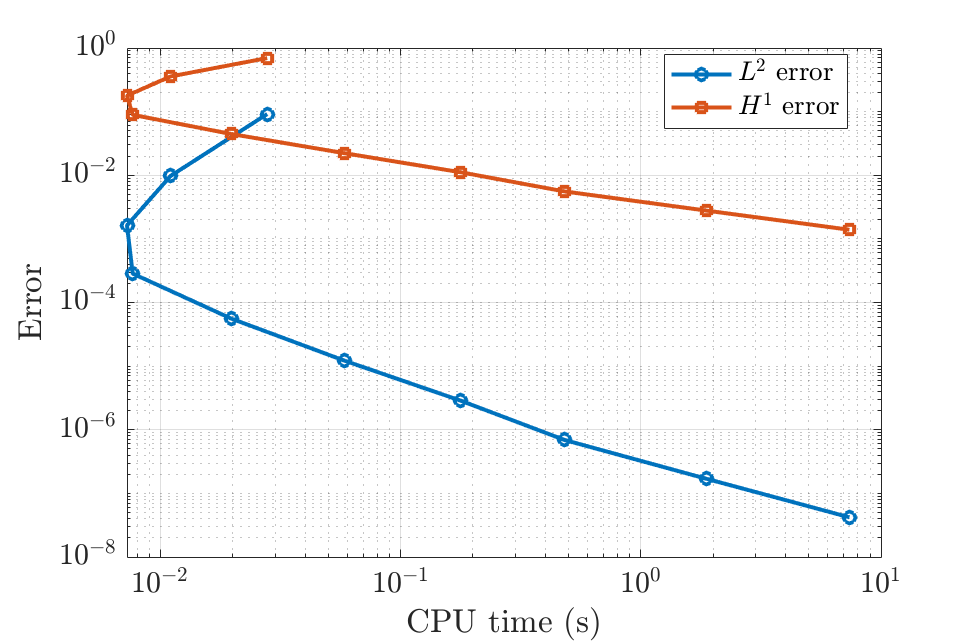}
   \caption{Computational cost vs. accuracy}
	\end{subfigure}
	\caption{$\mbox{L}^2$ and $\mbox{H}^1$-errors, convergence rate, and experimental order of convergence}
	\label{errorshistory}
\end{figure}
Moreover, the displacement amplitude decreases with time owing to the dynamic evolution of the solution, accompanied by corresponding reductions in the strain, stress, and local wave speed. From approximately $N_e=64$ onward, successive mesh refinements produce only negligible changes in all computed fields, demonstrating clear mesh convergence and confirming that the proposed HHT-$\alpha$ finite element formulation accurately captures the nonlinear mechanical response while maintaining numerical stability and consistency over the entire simulation.

Figures~\ref{errorshistory} summarizes the accuracy, convergence behavior, and computational efficiency of the proposed HHT-$\alpha$ finite element method. Figures~\ref{errorshistory} (a) and~\ref{errorshistory} (b) show the temporal evolution of the $\mbox{L}^2$- and $\mbox{H}^1$-errors, respectively. Both scenarios show that errors decrease monotonically with simulation time, and finer meshes produce smaller errors at each time step. The convergence plots in Figures \ref{errorshistory} (c) and (d) show the expected behavior, with numerical errors decreasing with increasing DOFs or mesh refinement. In particular, the $\mbox{L}^2$-error exhibits second-order convergence, whereas the $\mbox{H}^1$-error converges with first-order accuracy. This observation is further corroborated by the experimental orders of convergence shown in Figure~\ref{errorshistory} (e), where the computed EOC approaches the theoretical values of $2$ in the $\mbox{L}^2$-norm and $1$ in the $\mbox{H}^1$-norm as the mesh is refined. Finally, Figure~\ref{errorshistory} (f) shows the computational cost with respect to accuracy. We observe that a higher numerical accuracy is obtained at the expense of more CPU time, while the optimal convergence properties of the proposed method are preserved.
\begin{figure}[H]
	\centering
	\begin{subfigure}{0.48\textwidth}
  \includegraphics[width=\linewidth]{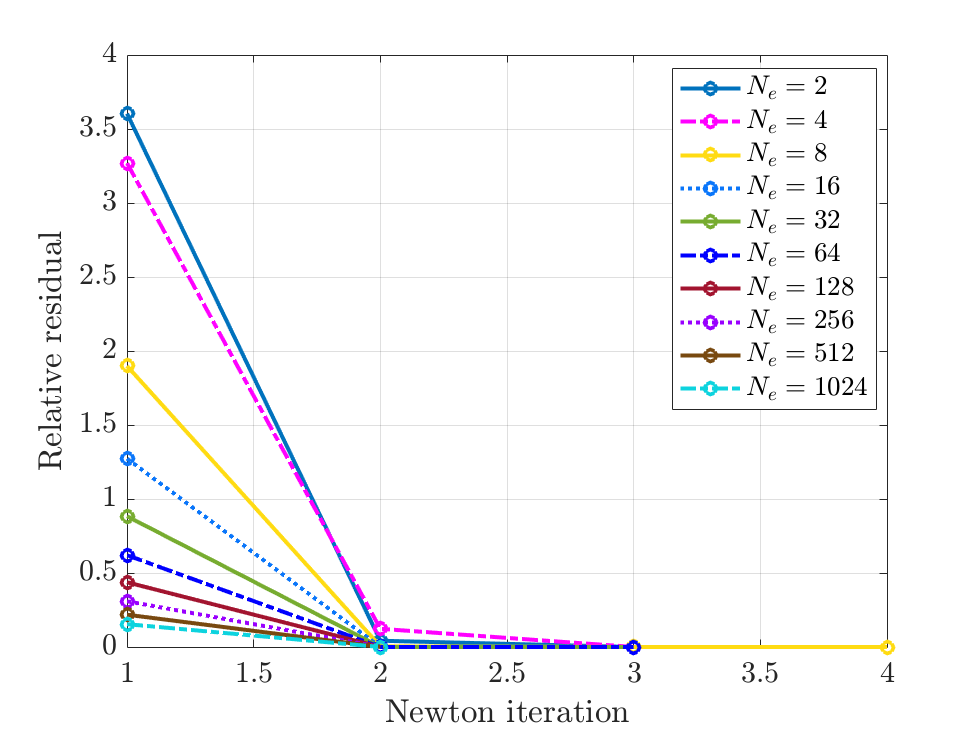}
		\caption{Relative residual vs. Newton iterations}
	\end{subfigure}
	\hfill
	\begin{subfigure}{0.48\textwidth}
		\centering
	\includegraphics[width=\linewidth]{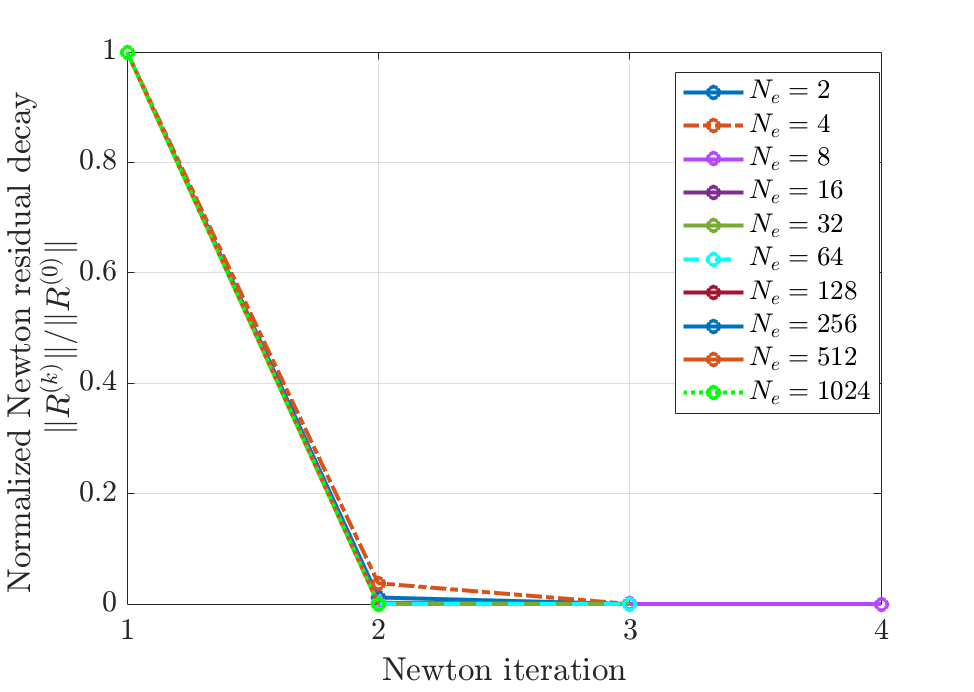}
		\caption{Normalized Newton residual decay}
	\end{subfigure}
	\caption{Computational performance}
	\label{compperform}
\end{figure}
Figure~\ref{compperform} illustrates the nonlinear solver performance of the proposed HHT-$\alpha$ finite element method over the sequence of uniformly refined meshes. Figure~\ref{compperform} (a) presents the relative residual history with respect to the Newton iteration number. For all mesh levels, the residual decreases rapidly, requiring at most four Newton iterations on the coarsest meshes and only two iterations for the finer discretizations, demonstrating the robustness of the nonlinear solver. Figure~\ref{compperform} (b) shows the normalized Newton residual decay, where the residual decreases by several orders of magnitude within just a few iterations. Almost identical decay profiles were obtained for different mesh resolutions, suggesting that the convergence of the Newton method is nearly mesh-independent. These results confirm the effectiveness of the consistent tangent Jacobian and demonstrate that the proposed nonlinear solution strategy converges fast, is stable, and is computationally efficient throughout the simulations.

\begin{figure}[h!]
	\centering
	\begin{subfigure}{0.48\textwidth}
		\includegraphics[width=\linewidth,height=5cm]{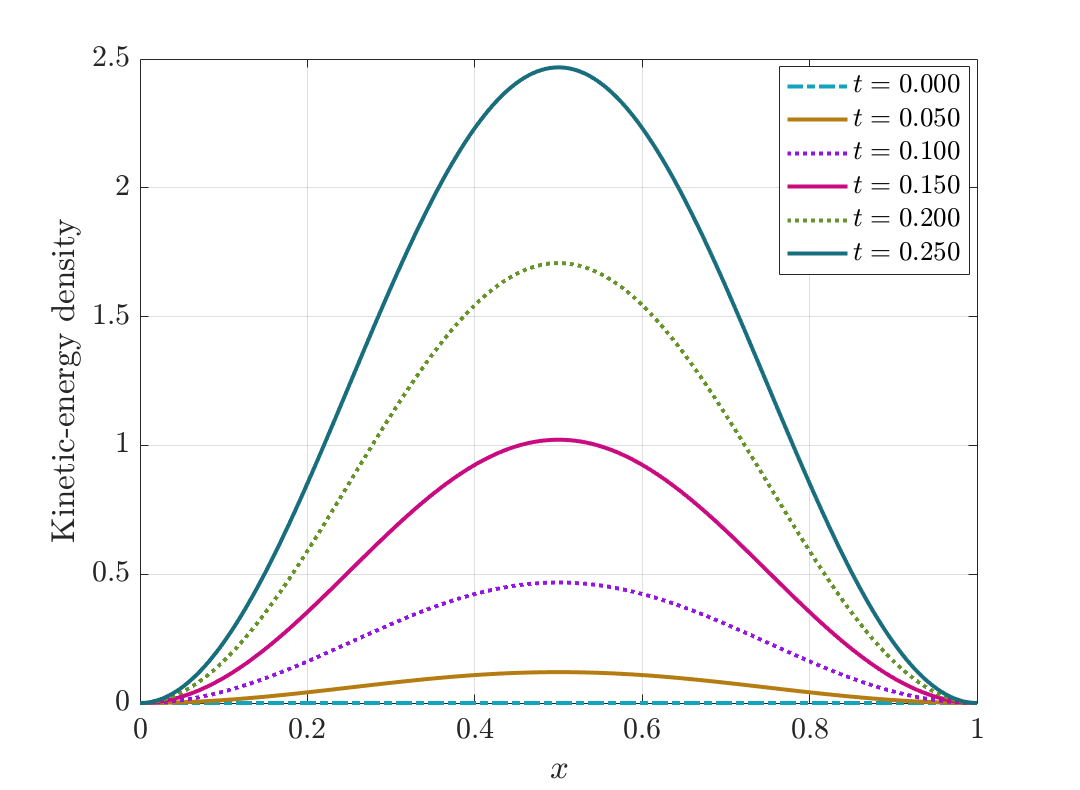}
		\caption{Kinetic energy density}
	\end{subfigure}
	\hfill
	\begin{subfigure}{0.48\textwidth}
	\includegraphics[width=\linewidth,height=5cm]{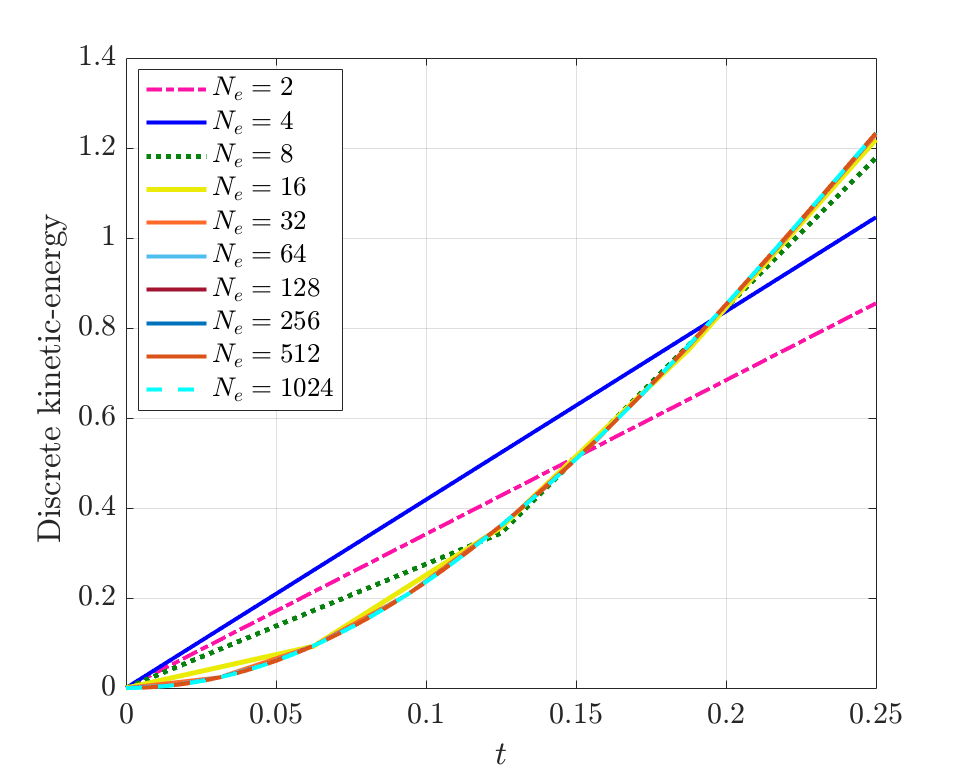}
		\caption{Discrete kinetic energy vs. time}
	\end{subfigure}
	\vspace{0.2cm}
	\begin{subfigure}{0.48\textwidth}
	\includegraphics[width=\linewidth,height=5cm]{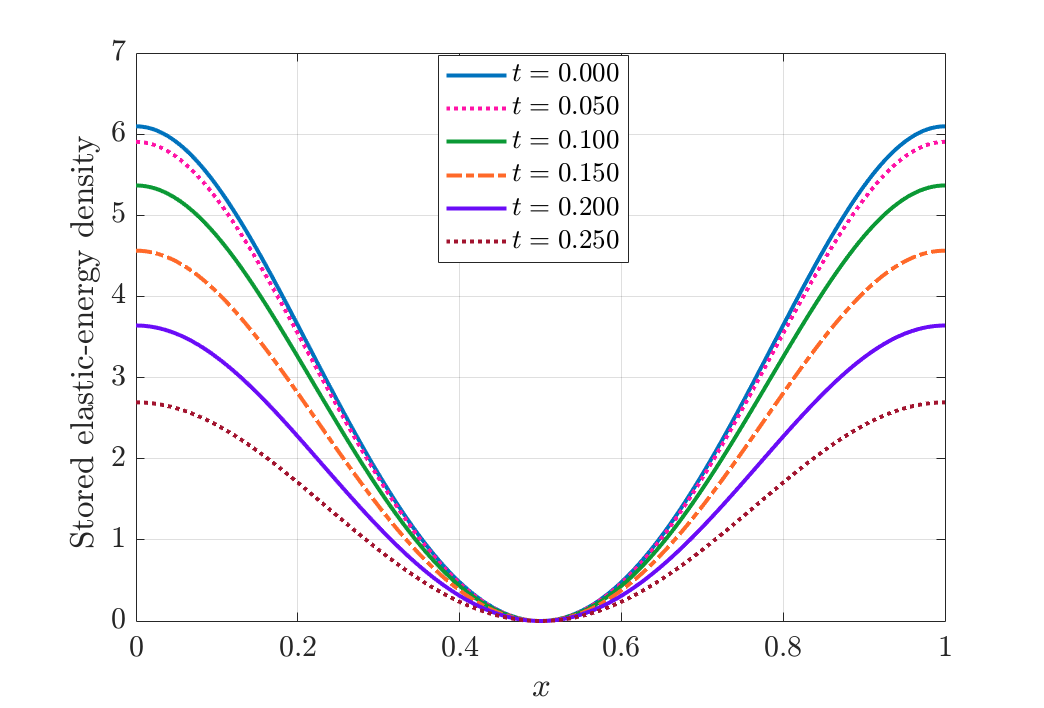}
\caption{Stored elastic-energy density}
	\end{subfigure}
	\hfill
	\begin{subfigure}{0.48\textwidth}
	\includegraphics[width=\linewidth,height=5cm]{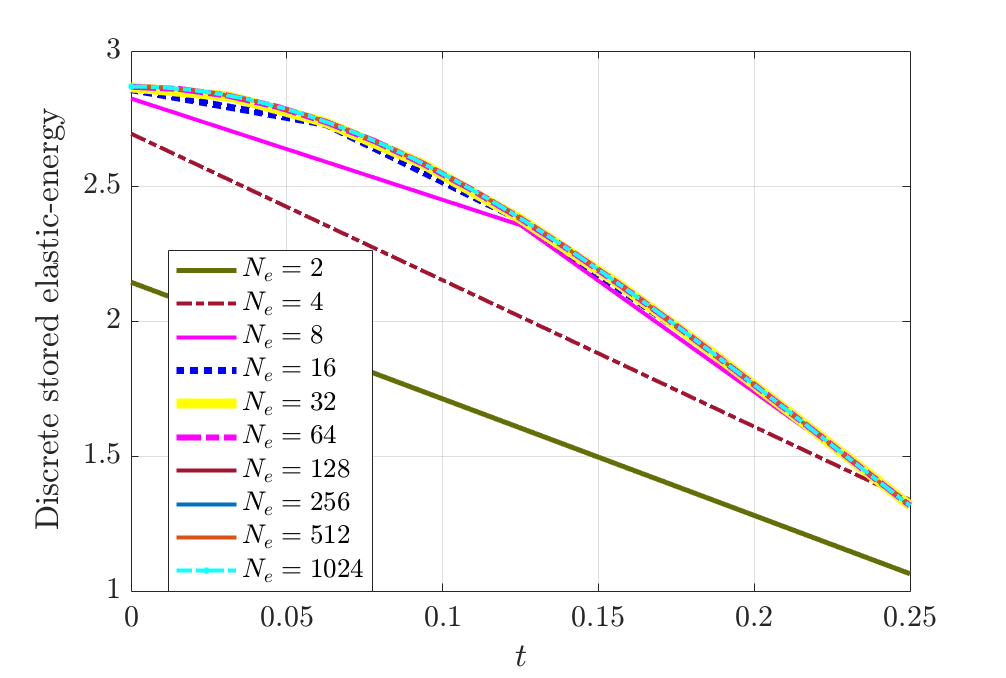}
\caption{Discrete elastic-energy vs time}
	\end{subfigure}
	
	\vspace{0.2cm}
	\begin{subfigure}{0.48\textwidth}
	\includegraphics[width=\linewidth,height=5cm]{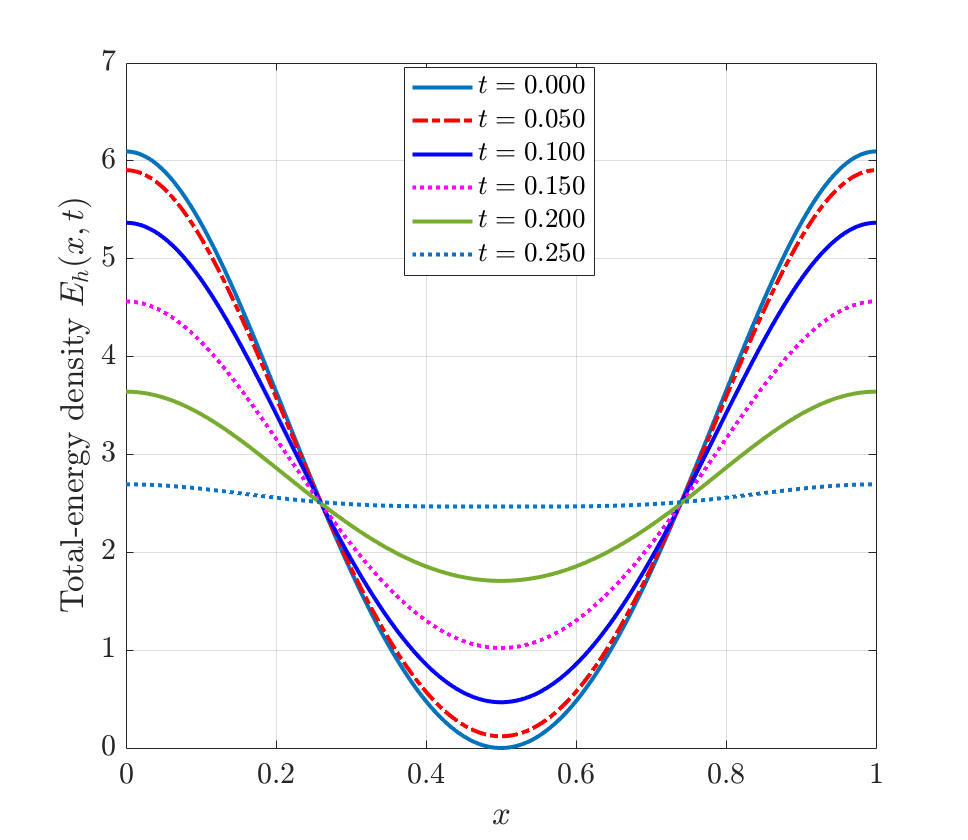}
\caption{Total energy density}
	\end{subfigure}
	\hfill
	\begin{subfigure}{0.48\textwidth}
	\includegraphics[width=\linewidth,height=5cm]{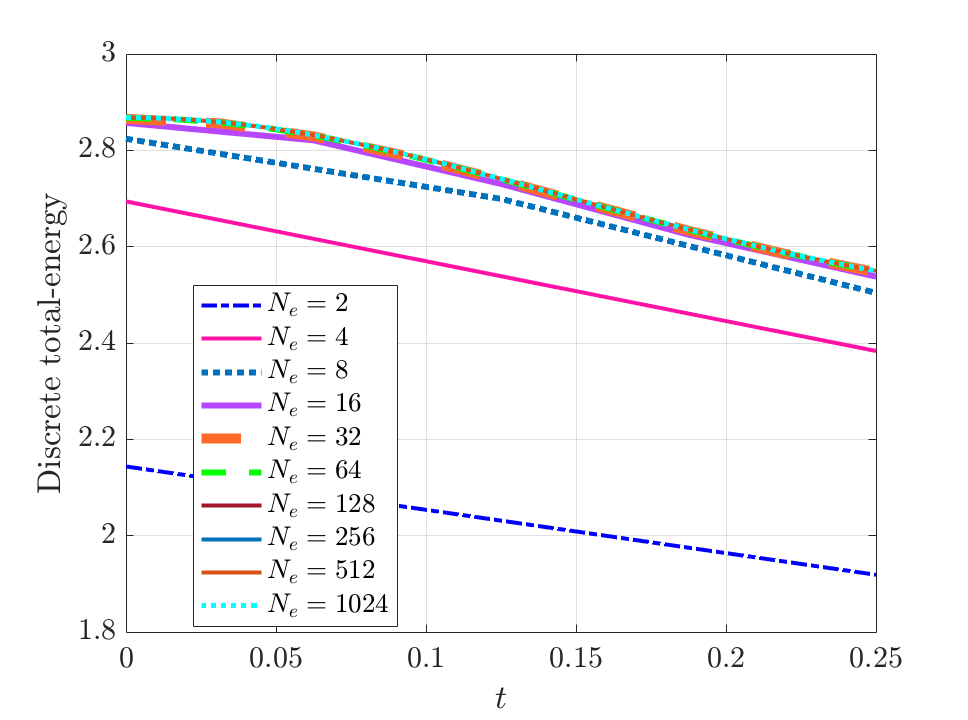}
\caption{Total discrete energy vs time}
	\end{subfigure}
	\caption{Evolution of the kinetic, stored elastic, and total energy during the simulation.
		Figures (a), (c), and (e) show the spatial distributions of the kinetic-energy density, stored elastic-energy density, and total energy density, respectively, at selected time instants. Figures (b), (d), and (f) present the corresponding discrete kinetic, stored elastic, and total energies as functions of time for different mesh resolutions.}
	\label{energydensity}
\end{figure}
Figures~\ref{energydensity} illustrates the evolution of the energy distribution and the corresponding discrete energy quantities during the dynamic response of the nonlinear strain-limiting wave equation. From a physical viewpoint, the total mechanical energy is composed of the kinetic and stored elastic energies, whose relative contributions evolve continuously as the wave propagates. Figures~\ref{energydensity} (a), (c), and (e) present the spatial distributions of the kinetic-energy density, stored elastic-energy density, and total energy density, respectively, at several representative time instants. Initially, the response is dominated by the stored elastic energy because the initial velocity is relatively small. As the solution evolves, elastic energy is progressively converted into kinetic energy, causing the kinetic-energy density to increase and attain its maximum near the center of the domain, where the particle velocity is largest, while the stored elastic-energy density decreases correspondingly. Throughout the simulation, the total energy density remains smooth and physically consistent, indicating that the proposed HHT-$\alpha$ finite element formulation accurately captures the energy transfer mechanism without introducing spurious numerical oscillations.

Figures~\ref{energydensity} (b), (d), and (f) display the temporal evolution of the corresponding discrete kinetic, stored elastic, and total energies for all mesh levels. The discrete kinetic energy increases with time as the wave accelerates, whereas the stored elastic energy decreases due to the release of the initially stored strain energy. Consequently, the total discrete energy exhibits only a slight reduction over the simulation interval, reflecting the controlled numerical dissipation introduced by the HHT-$\alpha$ time integration scheme. Furthermore, the energy histories obtained on successively refined meshes almost coincide, demonstrating excellent mesh independence and confirming the robustness, stability, and consistency of the proposed numerical method.

Overall, the results verify the best convergence rates of the proposed numerical approach as expected for piecewise linear finite element discretizations. For the spatial discretization and for the application of the HHT-$\alpha$ time integration scheme for the one-dimensional nonlinear strain-limiting wave equation, there is a great deal of agreement between the generated solution and the theoretical and empirically observed evolution.

In the following problem, we consider the mix/nonhomogeneous boundary condition and will analyze the boundary effect on the solution at the various time steps and compare with the exact solution.  

\begin{exam}[Non-homogeneous Dirichlet  Boundary]
Let $(\texttt{a}, \texttt{b})$ with $\texttt{a}=0$, $\texttt{b}=1$ and $\texttt{T}>0$. Consider the nonlinear hyperbolic problem
\begin{equation}\label{eq:ex-pde}
\partial_{tt}u(x,t)-\partial_x\!\left(\frac{ \partial_x u}{[1-\beta^2\,|\partial_x u(x,t)|^2]^{\frac{1}{2}}}\right)=f(x,t),
\qquad (x,t)\in(0,1)\times(0, \texttt{T}), 
\end{equation}
with non-zero Dirichlet boundary conditions
\begin{equation}\label{eq:ex-bc}
u(0,t)=0,\qquad u(1,t)=\psi(t), \qquad t\in(0, \texttt{T}),
\end{equation}
and initial conditions
\begin{equation}\label{eq:ex-ic}
u(x,0)=\phi_1(x),\qquad u_t(x,0)=\phi_2(x),\qquad x\in(0,1).
\end{equation}
Further, we have the following exact solution
\begin{equation}\label{eq:ex-uex}
u(x,t)=x\,t+\sin(\pi x)\cos(\pi t).
\end{equation}
Thus, we define the boundary conditions  $\psi(t)=u(1,t)=t,$ so that $u(0,t)=0$ and $u(1,t)=t$, however the initial conditions  are given by $\phi_1(x)=u(x,0)=\sin(\pi x),~ \text{and}~ \phi_2(x)=u_t(x,0)=x,$ respectively. 
Moreover, the load function $f$ is computed using the formula
\begin{equation}\label{eq:ex-f}
f(x,t)=\partial_{tt} u(x,t)-\partial_x\!\left(\frac{\partial_xu(x,t)}{[1-\beta^2\,|\partial_xu(x,t)|^2]^{\frac{1}{2}}}\right).
\end{equation}
Further, we construct the lifting function satisfying \eqref{eq:ex-bc} as
$$\tilde{\psi}(x,t)=x\,\psi(t)=x\,t.$$
Finally, we set the solution $y=u-\tilde{\psi}$, which satifying $y(0,t)=y(1,t)=0$ and, then using \eqref{eq:ex-uex}, we have 
$y(x,t)=\sin(\pi x)\cos(\pi t).$
\end{exam}
Next, we verify the results on ten uniformly refined meshes $N_e=2^j$, $j=1,\,2,\ldots,\, 10$, with finite elements, in the same manner as in Example \ref{exple5.1}. The parameters and abbreviations have been implemented in accordance with the instructions provided in Tables \ref{computationalnotation} and \ref{parameterslist}, respectively. 
\begin{table}[H]
	\centering
	\caption{Convergence history of the finite element approximation.}
	\label{tab:conv}
	\begin{tabular}{ccccc}
		\toprule
		Mesh ($N_e$) &
		$\|u-u_h\|_{\mbox{L}^2}$ &
		EOC &
		$\|u-u_h\|_{\mbox{H}^1}$ &
		EOC \\
		\midrule
		2    & $9.1591\times10^{-2}$ & --    & $6.9455\times10^{-1}$ & -- \\
		4    & $9.7769\times10^{-3}$ & 3.2278 & $3.5715\times10^{-1}$ & 0.9595 \\
		8    & $1.6065\times10^{-3}$ & 2.6054 & $1.7841\times10^{-1}$ & 1.0013 \\
		16   & $2.7816\times10^{-4}$ & 2.5299 & $8.9096\times10^{-2}$ & 1.0018 \\
		32   & $5.4943\times10^{-5}$ & 2.3399 & $4.4526\times10^{-2}$ & 1.0007 \\
		64   & $1.2135\times10^{-5}$ & 2.1787 & $2.2260\times10^{-2}$ & 1.0002 \\
		128  & $2.8532\times10^{-6}$ & 2.0886 & $1.1129\times10^{-2}$ & 1.0001 \\
		256  & $6.9210\times10^{-7}$ & 2.0435 & $5.5647\times10^{-3}$ & 1.0000 \\
		512  & $1.7047\times10^{-7}$ & 2.0215 & $2.7823\times10^{-3}$ & 1.0000 \\
		1024 & $4.2302\times10^{-8}$ & 2.0107 & $1.3912\times10^{-3}$ & 1.0000 \\
		\bottomrule
	\end{tabular}
\end{table}
\begin{table}[htpb]
	\centering
	\caption{Mesh characteristics, nonlinear solver performance, and computational cost.}
	\label{tab:mesh}
	\begin{tabular}{cccccccc}
		\toprule
		$N_e$ &
		$N_{\rm dof}$ &
		$h$ &
		$\Delta t$ &
		$N_t$ &
		$\overline{N}_{\rm Newton}$ &
		$N_{\rm Newton}^{\max}$ &
		All Conv. \\
		\midrule
		2    & 3    & $5.0000\times10^{-1}$ & $2.5000\times10^{-1}$ & 1   & 4.00 & 4 & Yes \\
		4    & 5    & $2.5000\times10^{-1}$ & $2.5000\times10^{-1}$ & 1   & 4.00 & 4 & Yes \\
		8    & 9    & $1.2500\times10^{-1}$ & $1.2500\times10^{-1}$ & 2   & 4.00 & 4 & Yes \\
		16   & 17   & $6.2500\times10^{-2}$ & $6.2500\times10^{-2}$ & 4   & 3.25 & 4 & Yes \\
		32   & 33   & $3.1250\times10^{-2}$ & $3.1250\times10^{-2}$ & 8   & 3.00 & 3 & Yes \\
		64   & 65   & $1.5625\times10^{-2}$ & $1.5625\times10^{-2}$ & 16  & 3.00 & 3 & Yes \\
		128  & 129  & $7.8125\times10^{-3}$ & $7.8125\times10^{-3}$ & 32  & 2.00 & 2 & Yes \\
		256  & 257  & $3.9063\times10^{-3}$ & $3.9063\times10^{-3}$ & 64  & 2.00 & 2 & Yes \\
		512  & 513  & $1.9531\times10^{-3}$ & $1.9531\times10^{-3}$ & 128 & 2.00 & 2 & Yes \\
		1024 & 1025 & $9.7656\times10^{-4}$ & $9.7656\times10^{-4}$ & 256 & 2.00 & 2 & Yes \\
		\bottomrule
	\end{tabular}
\end{table}
It is evident from Table~\ref{tab:conv} that the proposed finite element approximation demonstrates convergence accuracy. The $\mbox{L}^2$-error exhibits nearly second-order convergence, while the $\mbox{H}^1$- error converges with first-order accuracy. These results agree with the theoretical convergence properties of linear finite elements and confirm the consistency and stability of the HHT-$\alpha$ discretization.

The characteristics of the mesh and the performance of the nonlinear solver are summarized in Table~\ref{tab:mesh}. With mesh refinement, the number of degrees of freedom and time steps increase systematically with a consistent refinement strategy. The average Newton iterations decrease on finer meshes, and all simulations converge successfully, demonstrating the robustness and efficiency of the proposed numerical solver.

\begin{figure}[h!]
	\centering
	\begin{subfigure}[b]{0.48\textwidth}
		\centering
		\includegraphics[width=\linewidth,height=0.20\textheight]{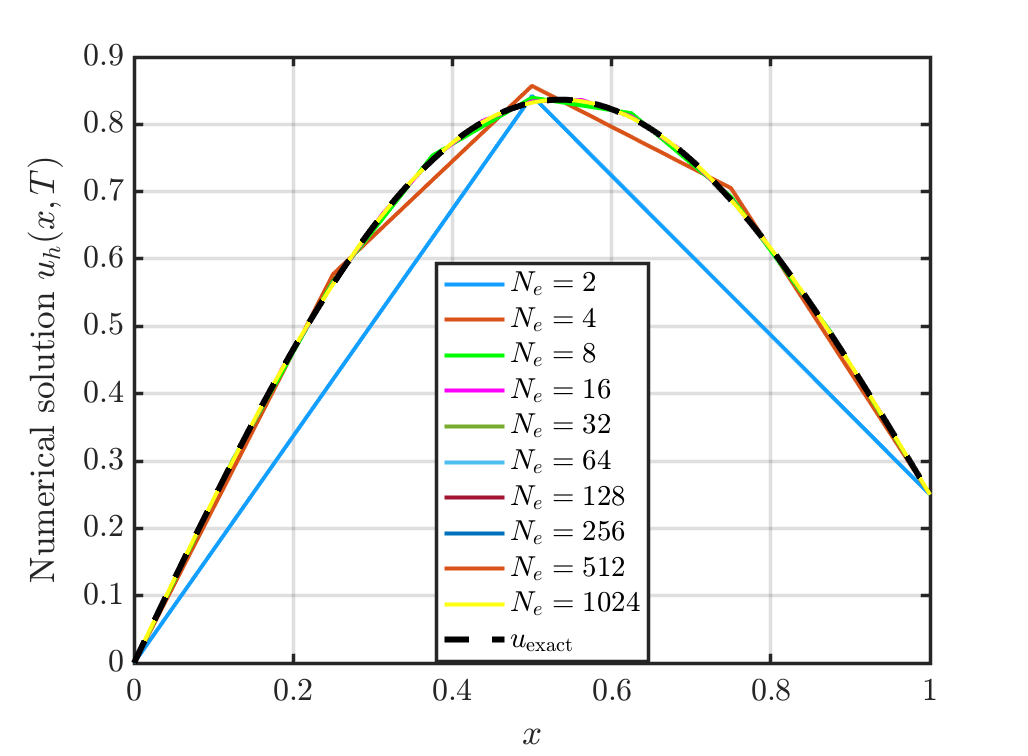}
		\caption{Approximate Solution at Final time}
		\label{ex2solution_profilesa}
	\end{subfigure}
	\hfill
	\begin{subfigure}[b]{0.48\textwidth}
		\centering
		\includegraphics[width=\linewidth,height=0.20\textheight]{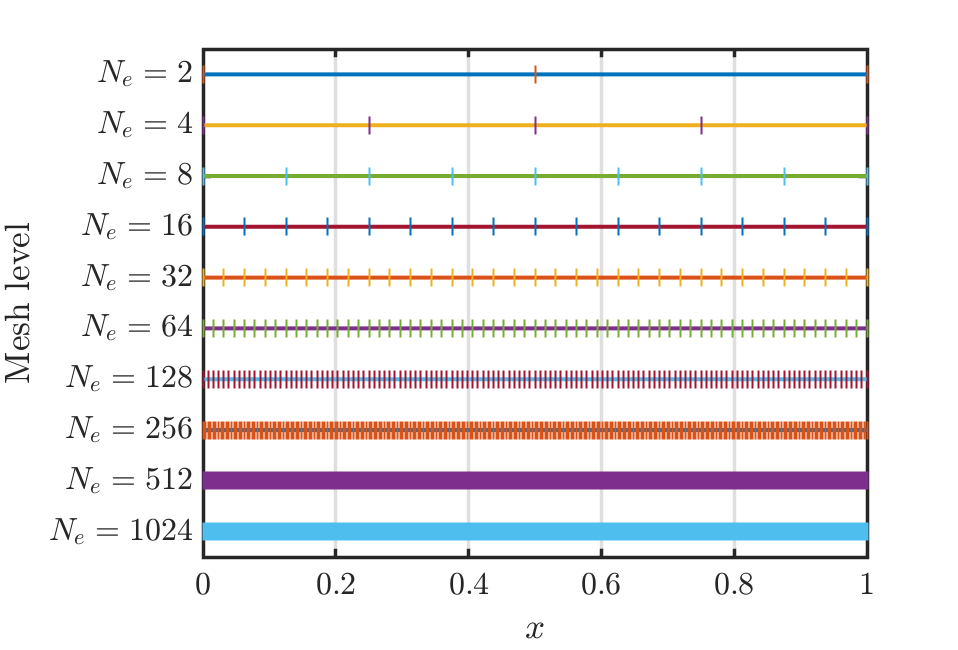}
		\caption{Uniform  finite-element mesh levels}
		\label{ex2meshes}
	\end{subfigure}
	
	\vspace{0.3cm}
	\begin{subfigure}[b]{0.48\textwidth}
		\centering
	\includegraphics[width=\linewidth,height=0.20\textheight]{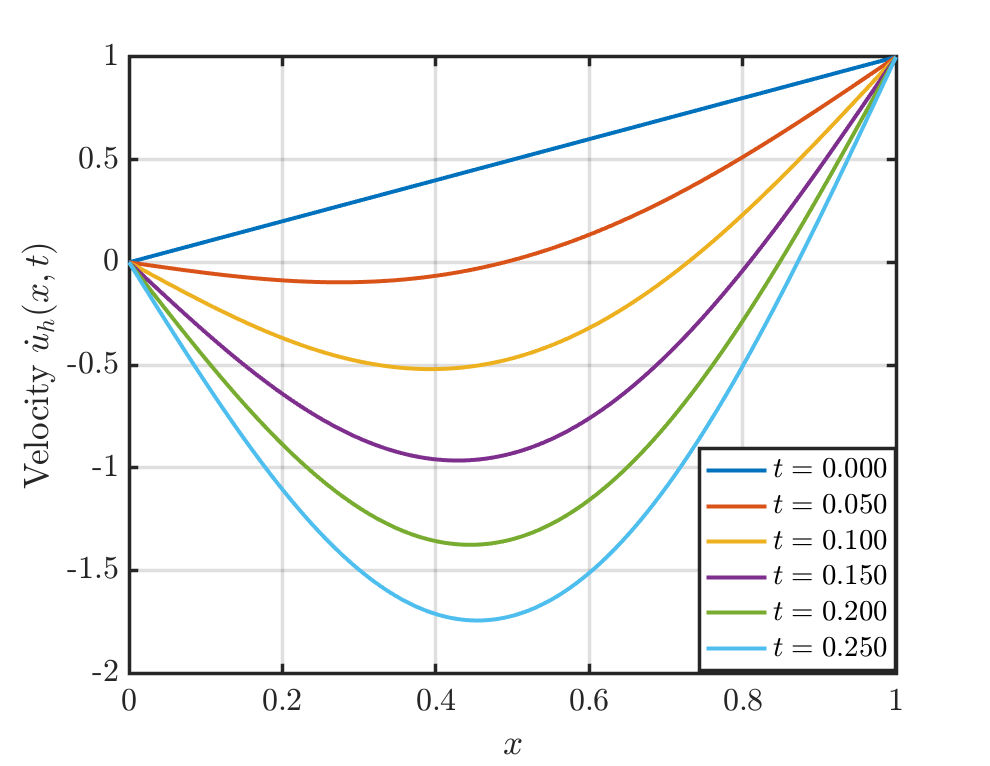}
		\caption{Approximate velocity}
		\label{ex2velocityprofilesc}
	\end{subfigure}
	\hfill
	\begin{subfigure}[b]{0.48\textwidth}
		\centering
\includegraphics[width=\linewidth,height=0.20\textheight]{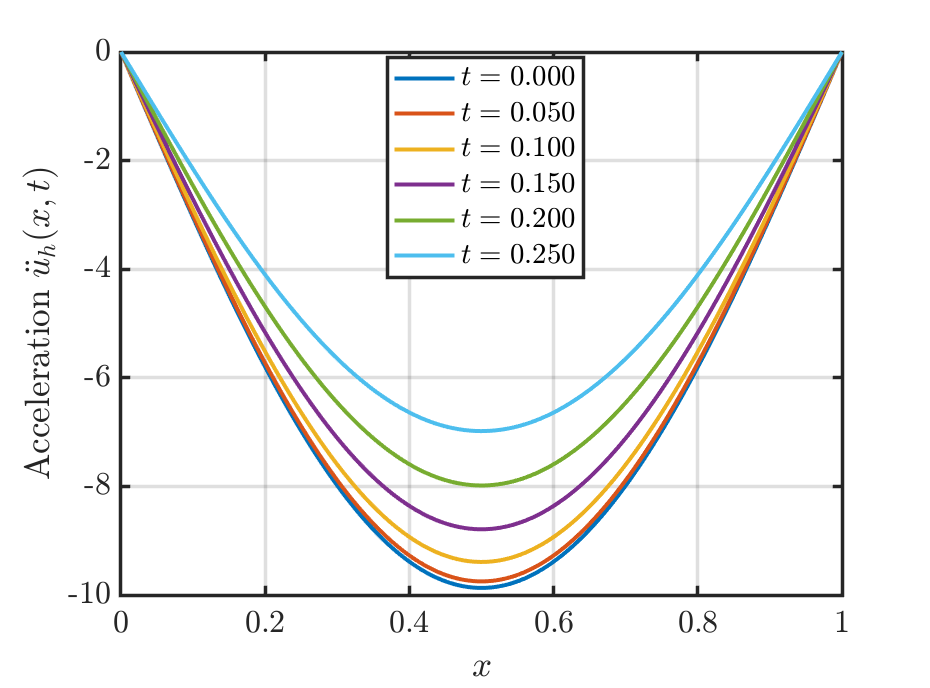}
		\caption{Approximate acceleration}
		\label{ex2accelerationsprofilesd}
	\end{subfigure}
	\caption{Profiles of approximate solution,  displacement, velocity and acceleration}
	\label{ex2profilesdva}
\end{figure}

Table~\ref{tab:wave} reports the nodal error and local wave-speed statistics. The nodal $L^\infty$-error decreases steadily with the mesh refinement, indicating that the numerical accuracy is improving. Furthermore, the local wave speed quantities tend to mesh-independent values, and the constitutive denominator remains strictly positive, confirming that the wave propagation remains stable and the nonlinear constitutive model is admissible for all the simulations.
\begin{table}[H]
	\centering
	\caption{Nodal error and local wave-speed statistics.}
	\label{tab:wave}
	\begin{tabular}{cccccc}
		\toprule
		$N_e$ &
		$\|e\|_{L^\infty_{\rm nodal}}$ &
		$c_{\min}^{T}$ &
		$c_{\max}^{T}$ &
		$c_{\max}^{ST}$ &
		$D_{\min}^{ST}$ \\
		\midrule
		2    & $1.0220\times10^{-2}$ & 1.0713 & 1.1577 & 1.2408 & 0.7500 \\
		4    & $2.5106\times10^{-2}$ & 1.0177 & 1.3554 & 1.6818 & 0.5000 \\
		8    & $7.4817\times10^{-3}$ & 1.0016 & 1.4121 & 1.9368 & 0.4142 \\
		16   & $2.0524\times10^{-3}$ & 1.0000 & 1.4292 & 2.0254 & 0.3902 \\
		32   & $5.3570\times10^{-4}$ & 1.0003 & 1.4331 & 2.0710 & 0.3788 \\
		64   & $1.3656\times10^{-4}$ & 1.0000 & 1.4340 & 2.0772 & 0.3773 \\
		128  & $3.4446\times10^{-5}$ & 1.0000 & 1.4342 & 2.0788 & 0.3769 \\
		256  & $8.6475\times10^{-6}$ & 1.0000 & 1.4342 & 2.0792 & 0.3768 \\
		512  & $2.1662\times10^{-6}$ & 1.0000 & 1.4343 & 2.0793 & 0.3768 \\
		1024 & $5.4209\times10^{-7}$ & 1.0000 & 1.4343 & 2.0793 & 0.3768 \\
		\bottomrule
	\end{tabular}
\end{table}
Figure~\ref{ex2profilesdva} illustrates the numerical behavior of the proposed HHT-$\alpha$ finite element method. Figure~\ref{ex2profilesdva} (a) shows that the numerical solution converges to the manufactured exact solution as the mesh is refined, while Figure~ \ref{ex2profilesdva} (b) presents the sequence of uniformly refined finite-element meshes used in the convergence study. Figures~ \ref{ex2profilesdva}  (c) and~ \ref{ex2profilesdva} (d) depict the temporal evolution of the velocity and acceleration fields, respectively, demonstrating smooth and physically consistent dynamic responses throughout the simulation. Overall, the results confirm the accuracy, stability, and convergence of the proposed numerical scheme.

\begin{figure}[H]
	\centering
	\includegraphics[width=\linewidth,height=0.45\textheight]{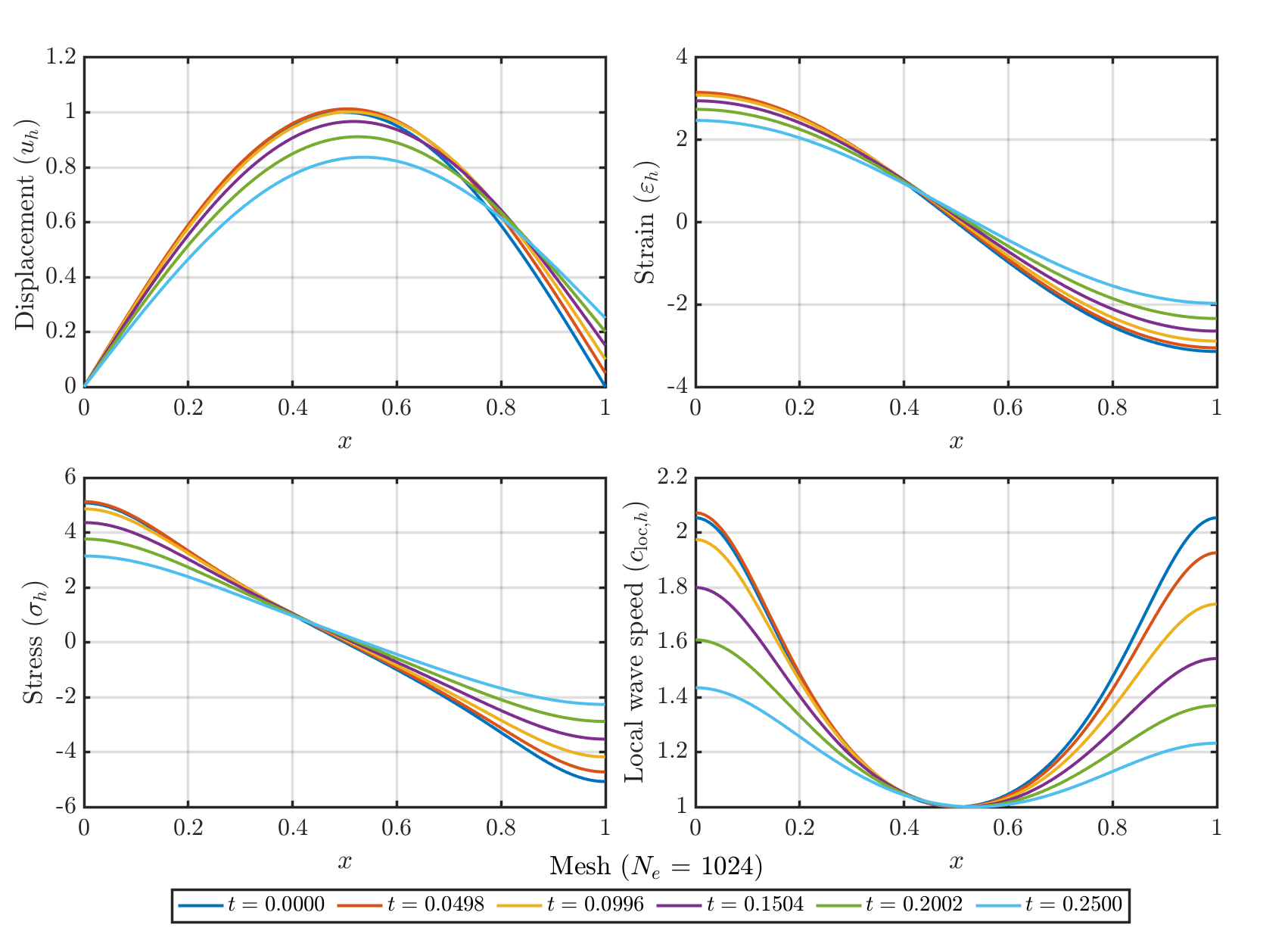}
	\caption{Evolution of the discrete displacement, strain, stress, and local tangent wave speed at selected time instants for meshes with $N_e=512$ and $N_e=1024$.}
	\label{ex2mechdatamesh512n1024}
\end{figure}

The time evolution of the displacement, strain, stress and local tangent wave speed on the finest mesh is shown in Figure~\ref{ex2mechdatamesh512n1024}. As time evolves, the amplitude of the displacement decreases with corresponding variations in the strain and stress fields induced by the nonlinear constitutive response. The local tangent wave speed evolves consistently with the strain distribution, attaining its minimum near the center of the domain while remaining bounded throughout the simulation. These results demonstrate the stable dynamic behavior and physical consistency of the proposed HHT-$\alpha$ finite element formulation.
\begin{figure}[H]
	\centering
	\begin{subfigure}{0.48\textwidth}
	\includegraphics[width=1.05\linewidth,height=5cm]{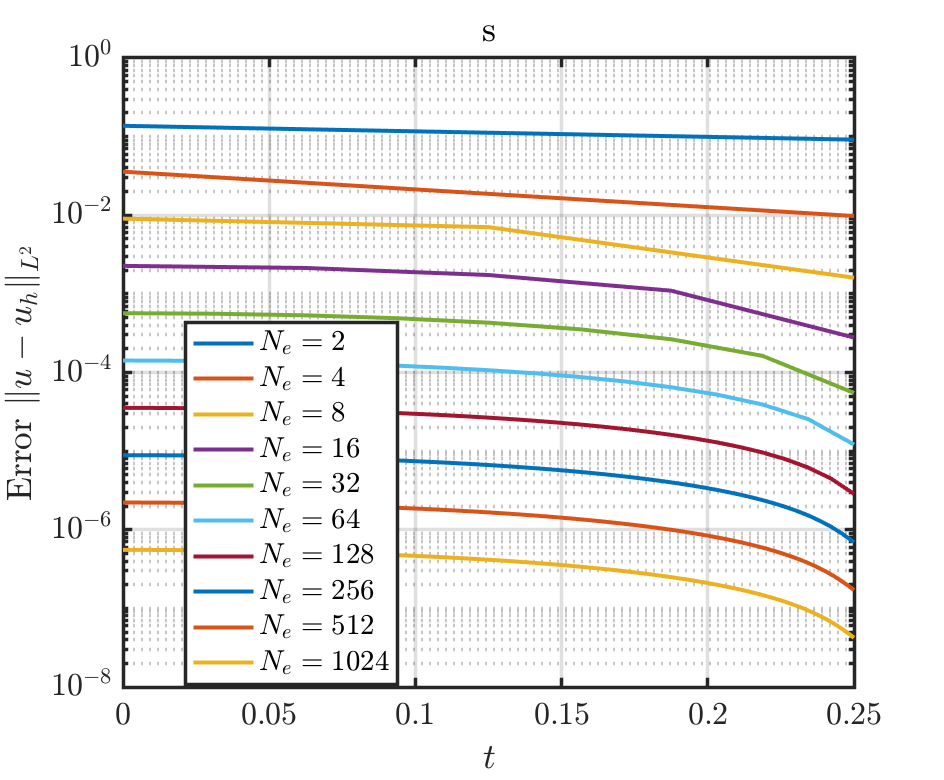}
		\caption{$\mbox{L}^2$-error vs Time ($t$)}
	\end{subfigure}
	\hfill
	\begin{subfigure}{0.48\textwidth}
		\includegraphics[width=1.05\linewidth,height=5cm]{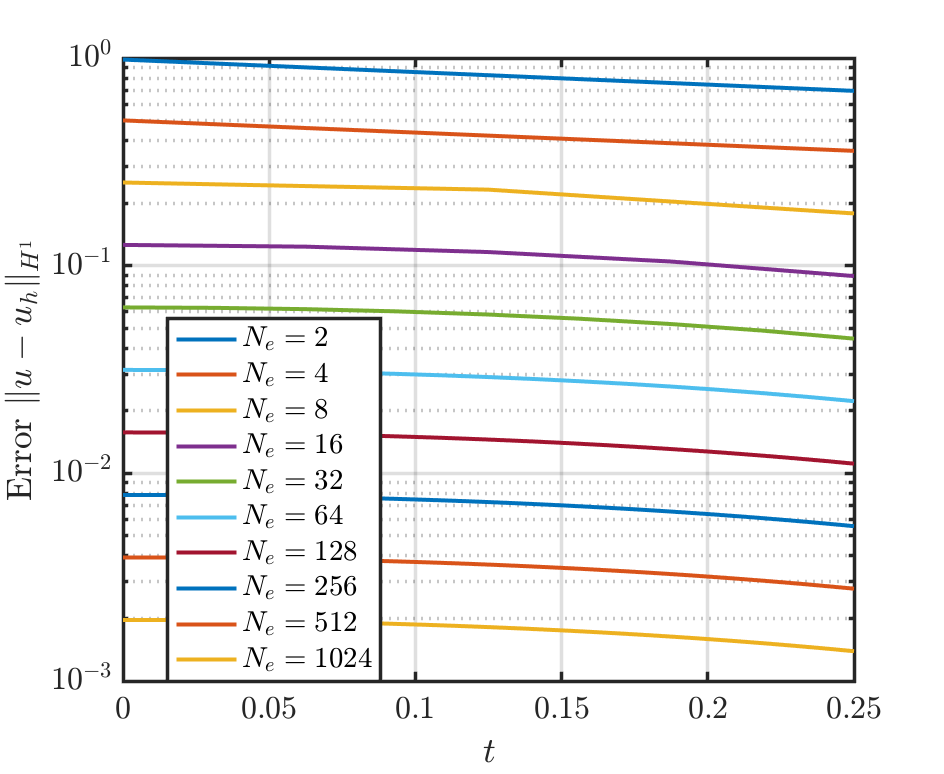}
		\caption{$\mbox{L}^2$-error vs. Time ($t$)}
	\end{subfigure}
	\vspace{0.2cm}
	\begin{subfigure}{0.48\textwidth}
	\includegraphics[width=1.05\linewidth,height=5cm]{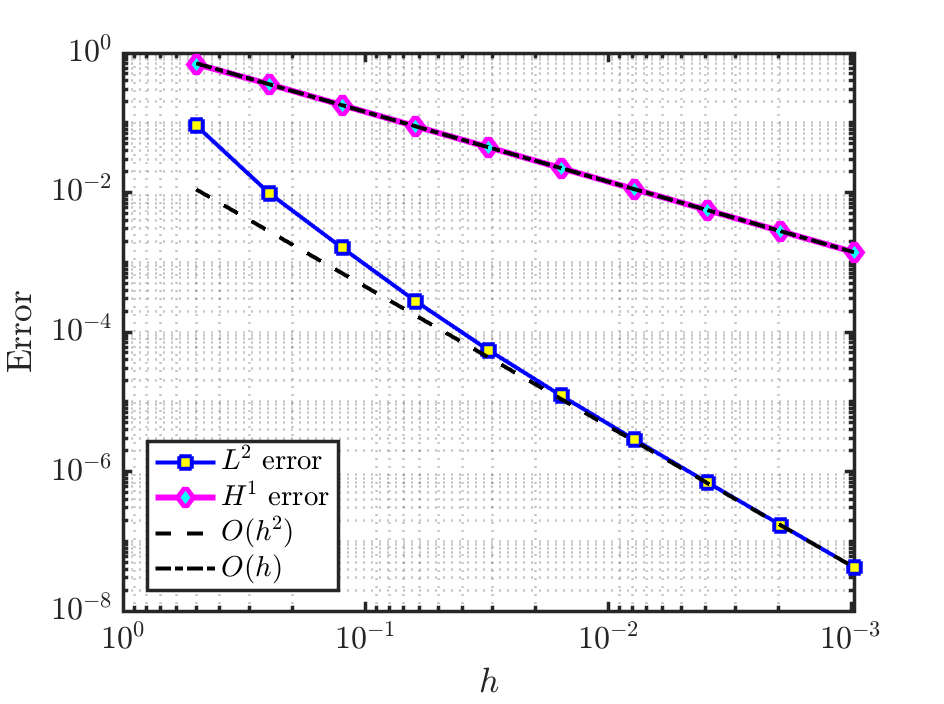}
		\caption{$\mbox{L}^2$ and $\mbox{H}^1$-error vs mesh size ($h$)}
	\end{subfigure}
	\hfill
	\begin{subfigure}{0.48\textwidth}
		\includegraphics[width=1.05\linewidth,height=5cm]{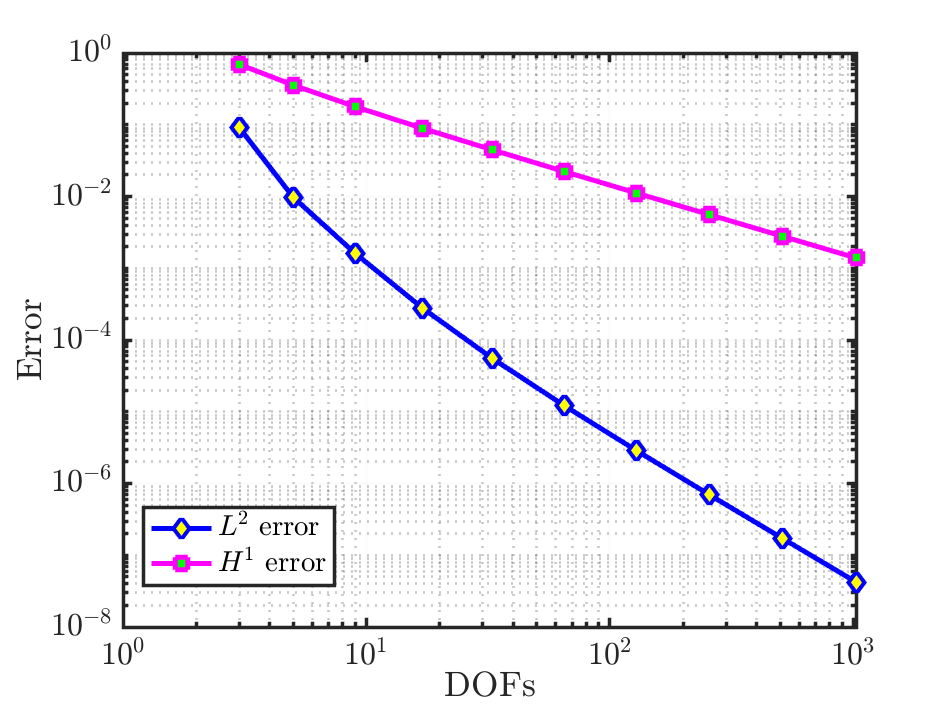}
		\caption{$\mbox{L}^2$ and $\mbox{H}^1$-error vs DOFs}
	\end{subfigure}
	\vspace{0.2cm}
	\begin{subfigure}{0.48\textwidth}
    \includegraphics[width=0.97\linewidth,height=5cm]{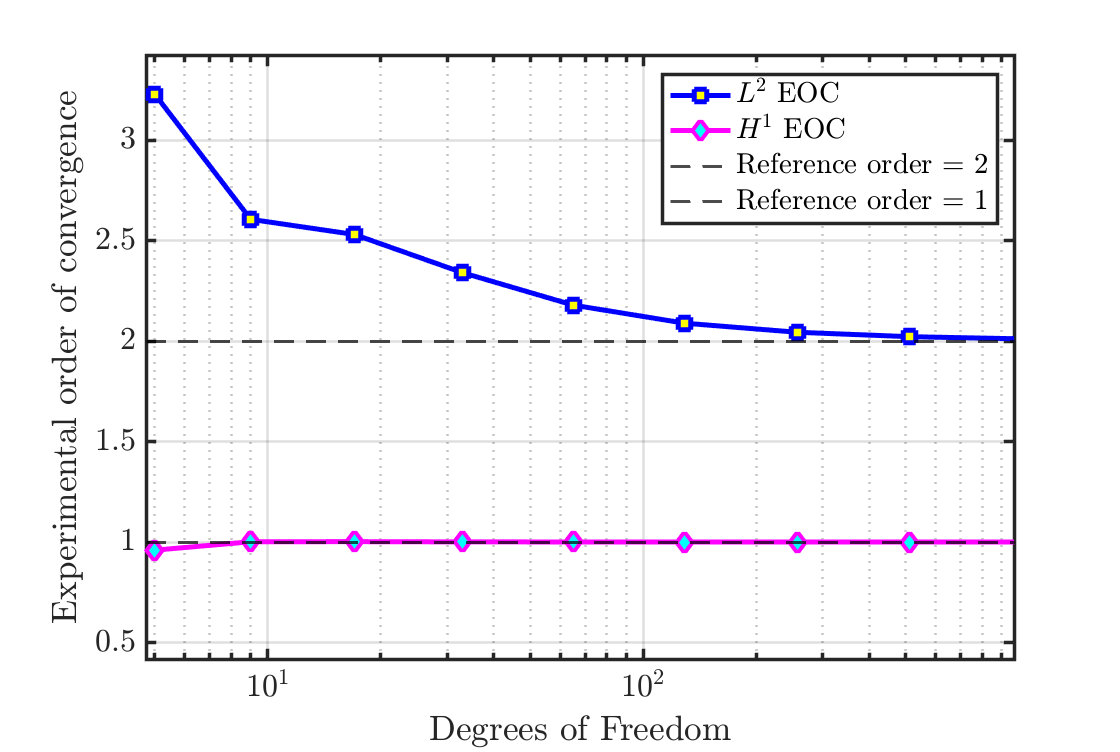}
		\caption{Experimental order of convergence}
	\end{subfigure}
	\hfill
	\begin{subfigure}{0.48\textwidth}
\includegraphics[width=1.05\linewidth,height=5cm]{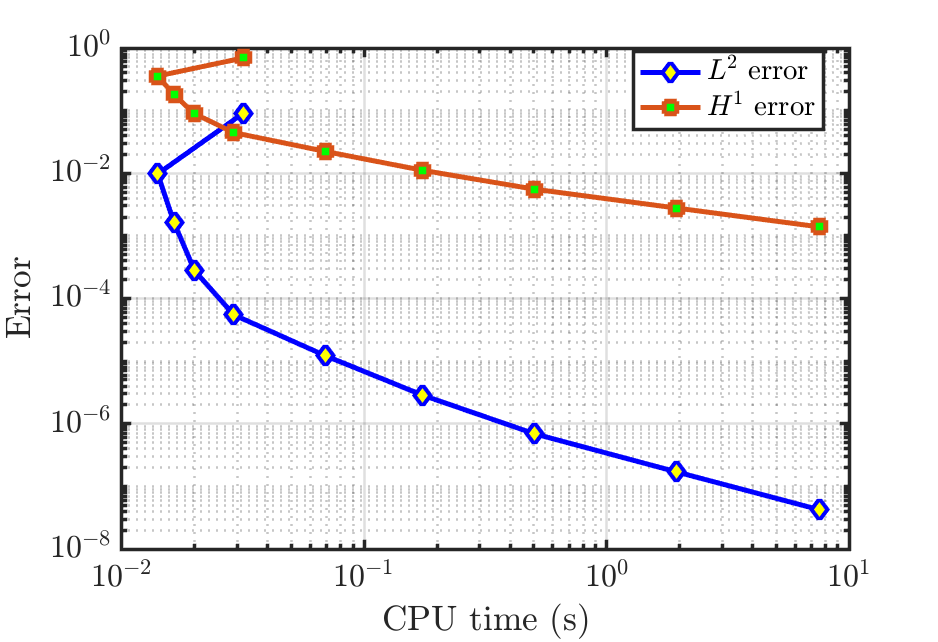}
		\caption{Computational cost vs. accuracy}
	\end{subfigure}
	\caption{$\mbox{L}^2$ and $\mbox{H}^1$-errors, convergence rate, and experimental order of convergence vs DOFs and }
	\label{ex2errorshistory}
\end{figure}
Figure~\ref{ex2errorshistory} shows the history of the $\mbox{L}^2$ and $\mbox{H}^1$ errors, the convergence plots and the computational efficiency of the proposed HHT-$\alpha$ finite element technique. The error histories are stable in time and systematically decrease with the mesh refinement, the convergence plots confirm the expected second-order convergence in the $\mbox{L}^2$-norm and first-order convergence in the $\mbox{H}^1$-norm. The experimental orders of convergence approach their theoretical orders for a large number of degrees of freedom and the analysis of computational cost shows that the higher accuracy is obtained with the predictable increase of CPU time. All these results verify the accuracy and efficiency of the proposed numerical scheme in a robust way.
\begin{figure}[H]
	\centering
	\begin{subfigure}{0.48\textwidth}
	\includegraphics[width=1.05\linewidth,height=5cm]{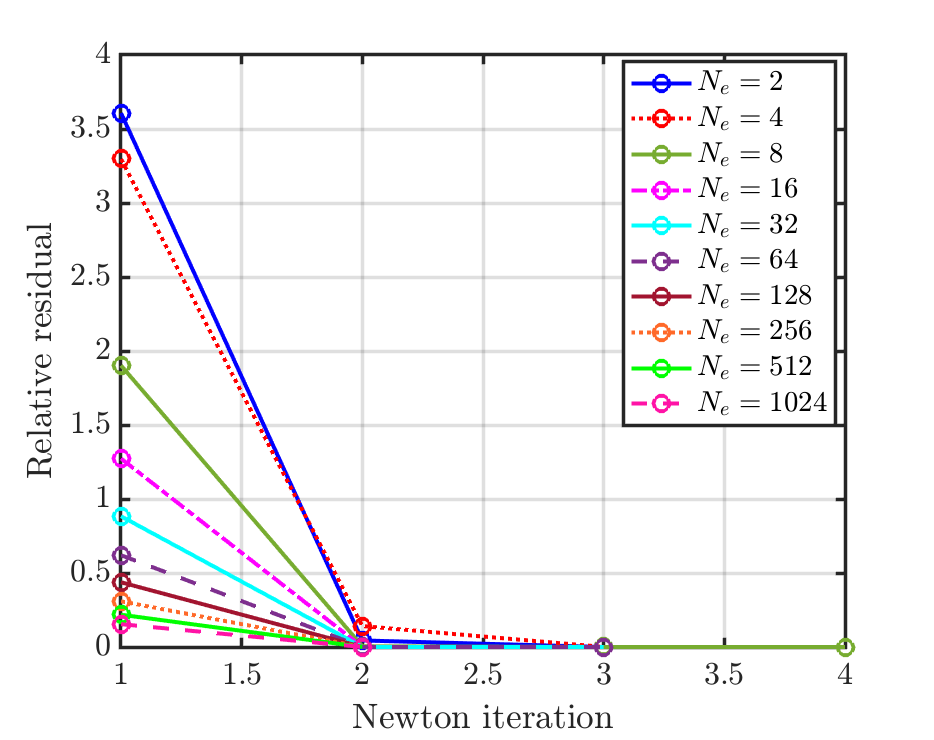}
		\caption{Relative residual vs. Newton iterations}
	\end{subfigure}
	\hfill
	\begin{subfigure}{0.48\textwidth}
		\centering
		\includegraphics[width=1.05\linewidth,height=5cm]{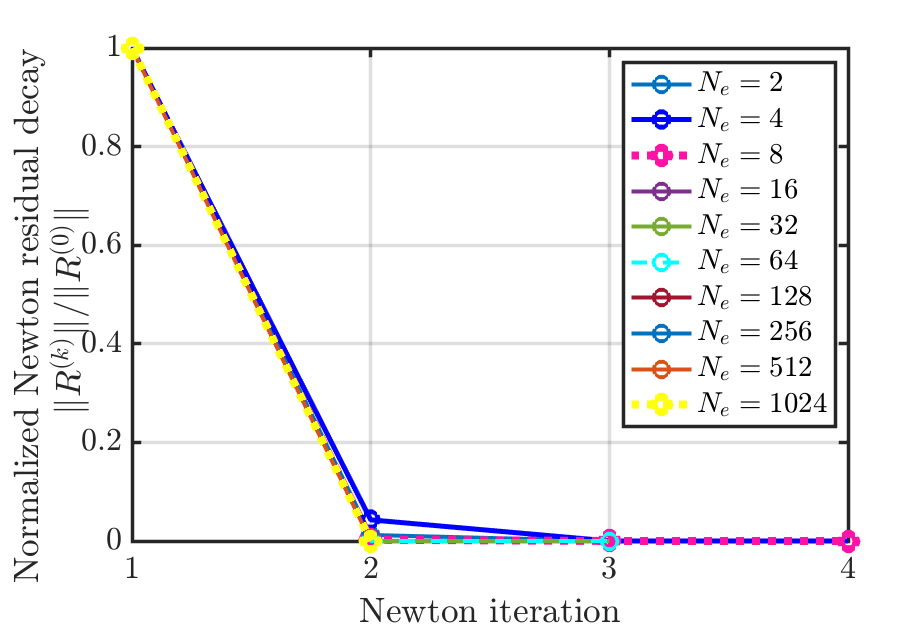}
		\caption{Normalized Newton residual decay}
	\end{subfigure}
	\caption{Computational performance}
	\label{ex2compperform}
\end{figure}
The convergence behavior of the Newton solver for various mesh refinements is shown in Figure~\ref{ex2compperform}. From Figure~(a), we observe that for all meshes the relative residual decreases rapidly within only a few Newton iterations, while in Figure~(b) the normalized residual decay is almost identical, indicating mesh-independent convergence. The fast reduction of the residual confirms the quadratic convergence characteristics of Newton's method and the effectiveness of the consistent tangent Jacobian. In summary, the results demonstrate the robustness, efficiency and numerical stability of the proposed nonlinear HHT-$\alpha$ finite element solver.

\begin{figure}[h!]
	\centering
	\begin{subfigure}{0.48\textwidth}
   \includegraphics[width=\linewidth,height=5cm]{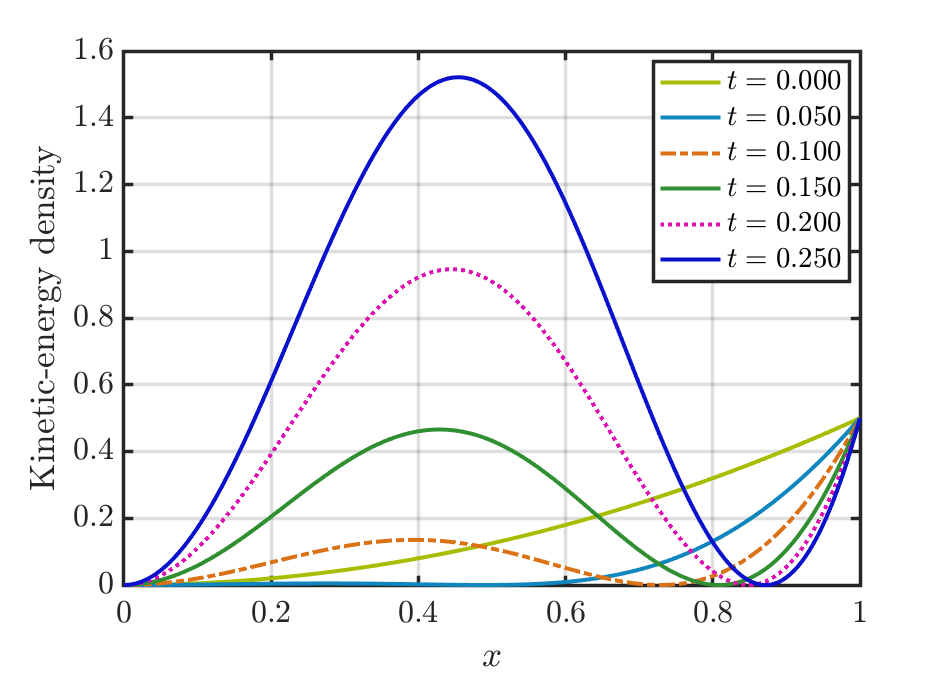}
		\caption{Kinetic energy density}
	\end{subfigure}
	\hfill
	\begin{subfigure}{0.48\textwidth}
		\includegraphics[width=\linewidth,height=5cm]{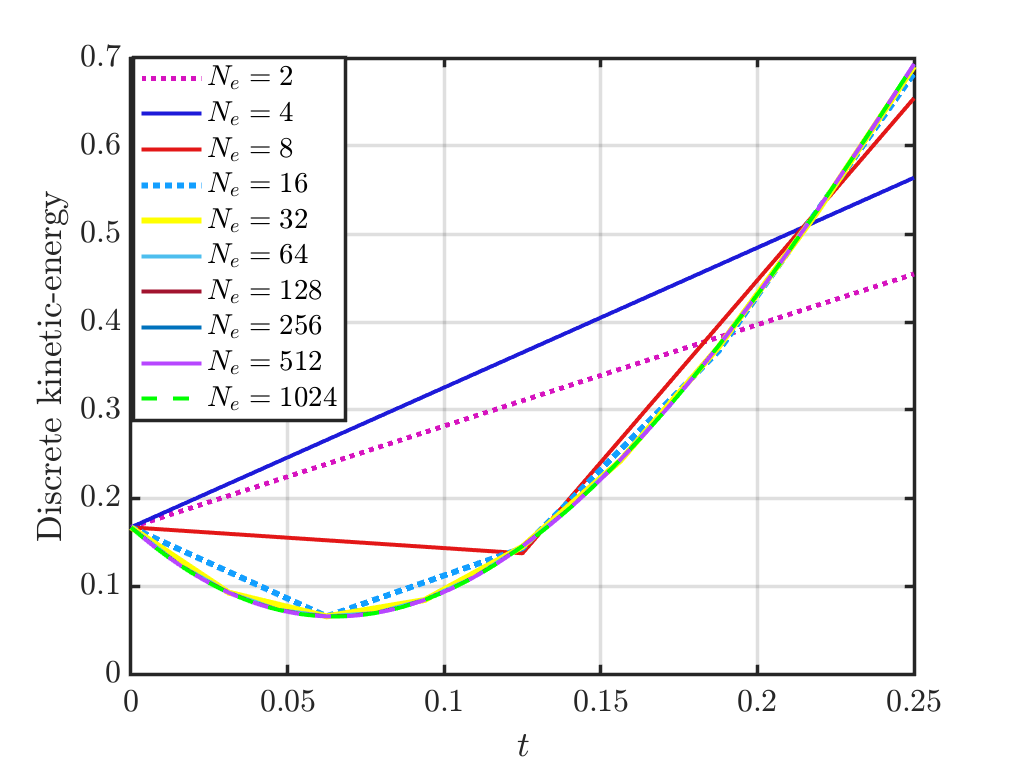}
		\caption{Discrete kinetic energy vs. time}
	\end{subfigure}
	\vspace{0.2cm}
	\begin{subfigure}{0.48\textwidth}
\includegraphics[width=\linewidth,height=5cm]{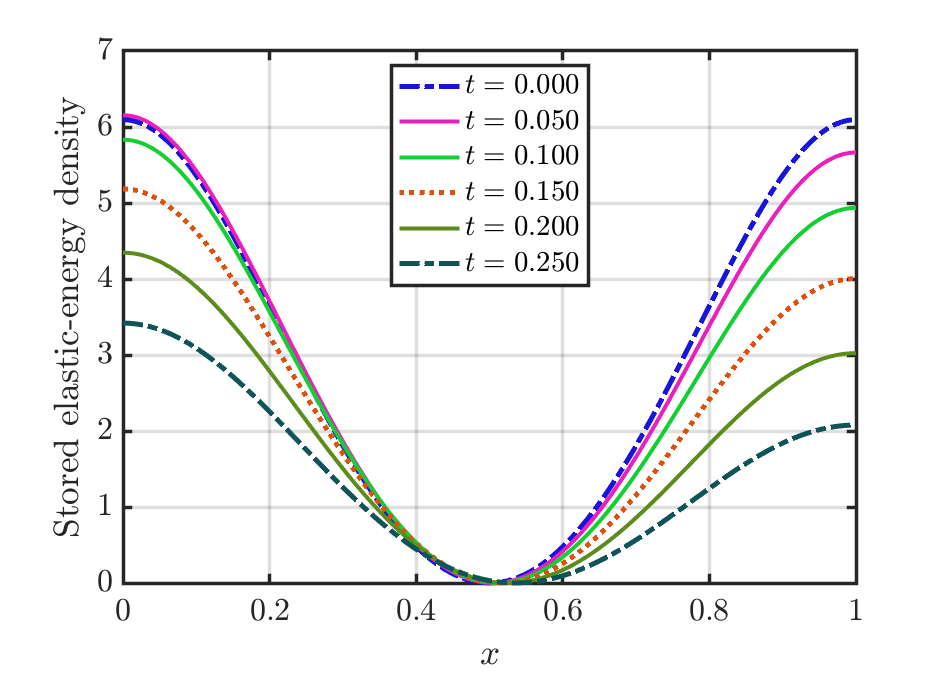}
		\caption{Stored elastic-energy density}
	\end{subfigure}
	\hfill
	\begin{subfigure}{0.48\textwidth}
		\includegraphics[width=\linewidth,height=5cm]{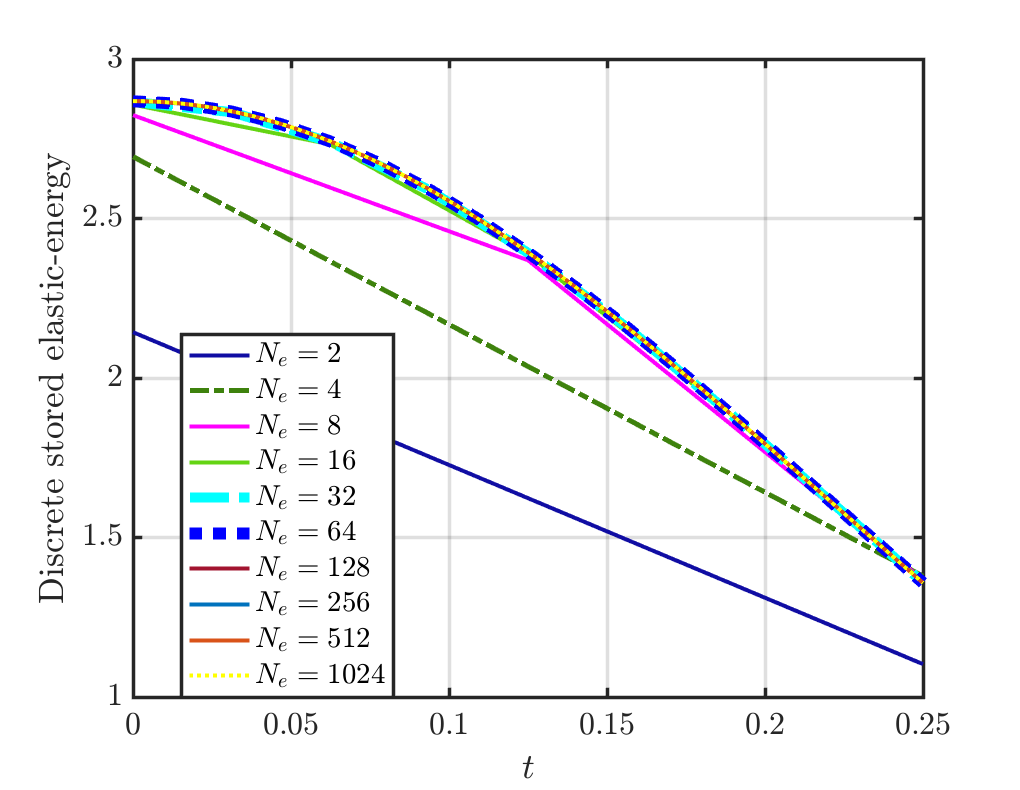}
		\caption{Discrete elastic energy vs time}
	\end{subfigure}
	
	\vspace{0.2cm}
	\begin{subfigure}{0.48\textwidth}
	\includegraphics[width=\linewidth,height=5cm]{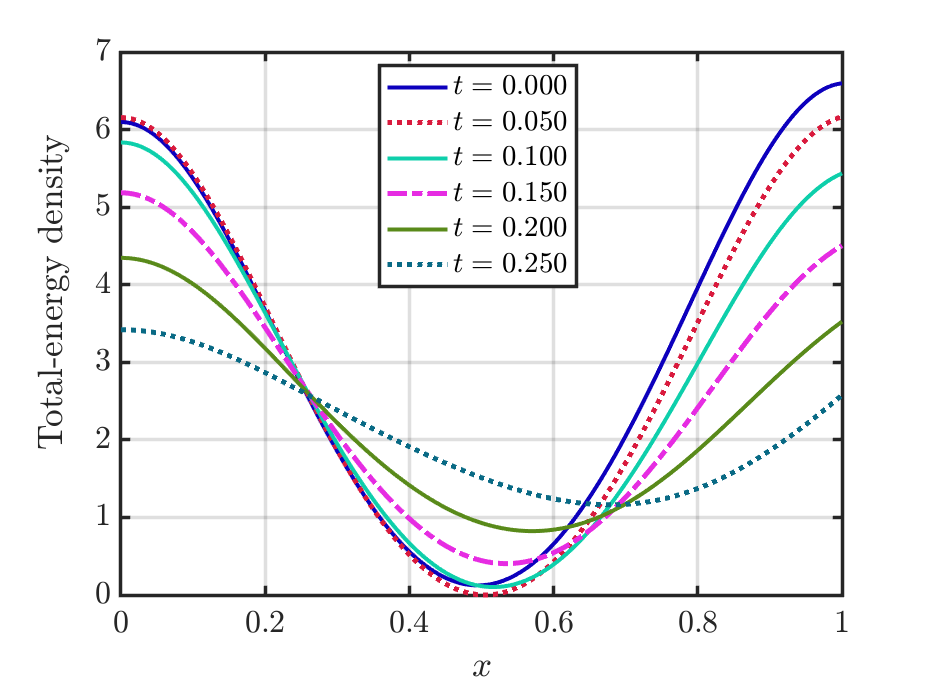}
		\caption{Total energy density}
	\end{subfigure}
	\hfill
	\begin{subfigure}{0.48\textwidth}
	\includegraphics[width=\linewidth,height=5cm]{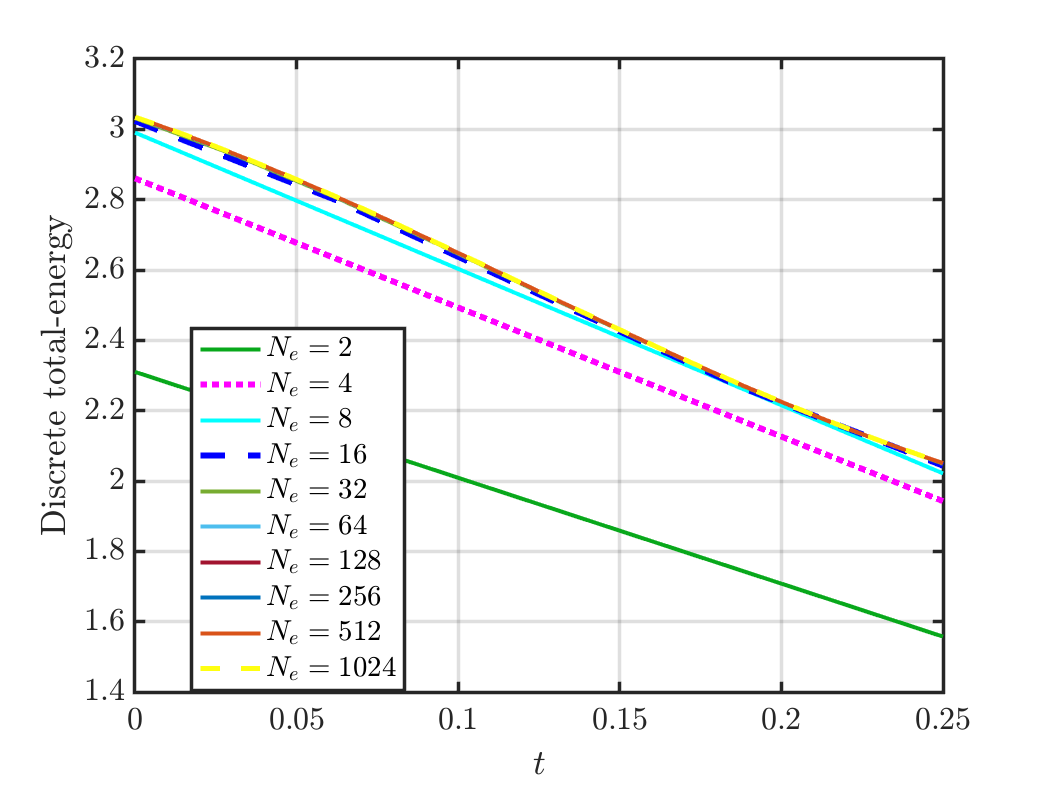}
		\caption{Total discrete energy vs time}
	\end{subfigure}
	\caption{Evolution of the kinetic, stored elastic, and total energy during the simulation.
		Figures (a), (c), and (e) show the spatial distributions of the kinetic-energy density, stored elastic-energy density, and total energy density, respectively, at selected time instants. Figures (b), (d), and (f) present the corresponding discrete kinetic, stored elastic, and total energies as functions of time for different mesh resolutions.}
	\label{ex2energydensity}
\end{figure}
Figure~\ref{ex2energydensity} shows the time evolution of the kinetic, stored elastic, and total energy calculated by the proposed HHT-$\alpha$ finite element method. The spatial distributions of the corresponding energy density at selected time instants are depicted in Figure~\ref{ex2energydensity} (a), (c) and (e), which reveal continuous energy transfer between kinetic and elastic parts during the wave propagation. The discrete kinetic, stored elastic, and total energy for different mesh resolutions are shown in Figure~\ref{ex2energydensity} (b), (d), and (f), which quantitatively assess the global energy evolution during the simulation.

The numerical results demonstrate that the energy histories are almost mesh independent when the mesh is refined, which suggests good numerical convergence. The total discrete energy decreases gradually, which indicates the algorithmic dissipation introduced by the HHT-$\alpha$ scheme, and the smooth evolution of the kinetic and elastic energies confirms the stability and robustness of the proposed nonlinear finite element formulation.

\section{Conclusion} \label{sec:6}
In this article,  we develop and analyze a continuous Galerkin finite element method for a class of quasilinear hyperbolic equations. The proposed fully discrete formulation combines continuous linear finite elements in space with the HHT-$\alpha$ time integration scheme and a consistent Newton solver, providing a stable and robust computational framework for nonlinear and degenerate regimes. The algorithmic dissipation introduced by the parameter $\alpha\in[0,1/3]$ effectively suppresses spurious high-frequency oscillations near degeneracy while preserving consistency and accuracy. Comprehensive numerical experiments demonstrate the expected second-order convergence in the $\mbox{L}^2$-norm and first-order convergence in the $\mbox{H}^1$-norm, rapid and nearly mesh-independent Newton convergence, physically consistent evolution of strain, stress, and local tangent wave speed, and smooth energy dissipation characteristic of the \texttt{HHT}-$\alpha$ method. The results demonstrate that the proposed formulation can accurately capture strain-limiting wave propagation with excellent stability, robustness, and computational efficiency.

The proposed nonlinear \texttt{HHT}-$\alpha$ finite element framework offers a reliable and efficient method for solving dynamic strain-limiting hyperbolic problems with algebraic nonlinear relationships between stress and strain. The robustness, accuracy, and good nonlinear convergence properties of the proposed method make it a promising basis for future developments such as multidimensional nonlinear elastodynamics, adaptive finite element methods, heterogeneous media, and nonlinear fracture and damage mechanics.

\section*{Acknoledgement}
RM, AGH, and HL gratefully acknowledge the financial support provided by the Dean of the College of Science at Texas A\&M University–Corpus Christi. RM and SMM acknowledge the support of the National Science Foundation under Grant No.\ 2316905.

\bibliographystyle{plain} 
\bibliography{references}    
\end{document}